\documentclass[CJK,10pt]{article}
\usepackage{geometry}
\usepackage{amsmath,amsfonts,amssymb,amsthm}
\usepackage{graphics}
\usepackage{dsfont}
\usepackage{epsfig}
\usepackage{esint}
\usepackage[colorlinks,linkcolor=blue]{hyperref}
\usepackage{stmaryrd}
\usepackage{tikz}
\usetikzlibrary{decorations.pathmorphing}
\newcommand{\lollipop}{%
\mathord{\raisebox{-0.30ex}{%
\begin{tikzpicture}[scale=0.4]
  \coordinate (root) at (0,0);
  \node[circle, fill=black, inner sep=1pt] (leaf) at (0,0.6) {};
  \draw[thick] (root) -- (leaf);
\end{tikzpicture}%
}}}
\newcommand{\cherry}{%
\mathord{\raisebox{-0.30ex}{%
\begin{tikzpicture}[scale=0.4]
  \coordinate (root) at (0,0);
  \node[circle, fill=black, inner sep=1pt] (left) at (-0.3,0.6) {};
  \node[circle, fill=black, inner sep=1pt] (right) at (0.3,0.6) {};
  \draw[thick] (root) -- (left);
  \draw[thick] (root) -- (right);
\end{tikzpicture}%
}}}

\newcommand{\chickenfoot}{%
\mathord{\raisebox{-0.30ex}{%
\begin{tikzpicture}[scale=0.4]
  \coordinate (root) at (0,0);
  \node[circle, fill=black, inner sep=1pt] (left) at (-0.4,0.6) {};
  \node[circle, fill=black, inner sep=1pt] (middle) at (0,0.6) {};
  \node[circle, fill=black, inner sep=1pt] (right) at (0.4,0.6) {};
  \draw[thick] (root) -- (left);
  \draw[thick] (root) -- (middle);
  \draw[thick] (root) -- (right);
\end{tikzpicture}%
}}}

\newcommand{\cherryS}{%
\mathord{\raisebox{-0.30ex}{%
\begin{tikzpicture}[scale=0.3]
  \coordinate (root) at (0,0);
  \node[circle, fill=black, inner sep=1pt] (left) at (-0.3,0.6) {};
  \node[circle, fill=black, inner sep=1pt] (right) at (0.3,0.6) {};
  \draw[thick] (root) -- (left);
  \draw[thick] (root) -- (right);
\end{tikzpicture}%
}}}

\newcommand{\chickenfootS}{%
\mathord{\raisebox{-0.30ex}{%
\begin{tikzpicture}[scale=0.35]
  \coordinate (root) at (0,0);
  \node[circle, fill=black, inner sep=1pt] (left) at (-0.4,0.6) {};
  \node[circle, fill=black, inner sep=1pt] (middle) at (0,0.6) {};
  \node[circle, fill=black, inner sep=1pt] (right) at (0.4,0.6) {};
  \draw[thick] (root) -- (left);
  \draw[thick] (root) -- (middle);
  \draw[thick] (root) -- (right);
\end{tikzpicture}%
}}}

\newcommand{\ignore}[1]{}

\newtheorem{definition}{Definition}
\newtheorem{proposition}{Proposition}
\newtheorem{theorem}{Theorem}
\newtheorem{remark}{Remark}
\newtheorem{lemma}{Lemma}

\newtheorem{assumption}{Assumption}

\newcommand{\cc}{\mathrm{C}}

\newcommand{\dd}{\mathrm{d}}

\newcommand{\bb}{\mathrm{B}}
\newcommand{\LL}{\mathrm{L}}

\usepackage{color}
\definecolor{darkred}{rgb}{0.9,0.1,0.1}
\definecolor{darkblue}{rgb}{0,0,0.7}
\definecolor{darkgreen}{rgb}{0,0.5,0}

\author{Nicolas Clozeau\thanks{Laboratoire IMATH, Université de Toulon, email: nicolas.clozeau@univ-tln.fr}}
\title{An inductive approach to stochastic estimates for singular SPDEs with correlated coefficient fields}
\begin{document}
\maketitle
\begin{abstract}
We develop an inductive approach to obtain stochastic estimates for singular SPDEs when the coefficient field is correlated with the driving noise. We implement our method, in the whole space, on the example of the $\varphi^{4}_2$-equation. Our approach is based on (infinite-dimensional) integration by parts with respect to Wick products of Gaussian random variables (more precisely, mollifications of space–time white noise) together with an inductive scheme over the moment's exponent.
\end{abstract}
\tableofcontents

\section{Introduction}

Sub-critical singular stochastic partial differential equations (SPDEs) in homogeneous media are by now well understood, and general techniques have been developed to perform renormalizations in a systematic manner; see, for instance, the theory of regularity structures by Hairer \cite{hairer2014theory}, paracontrolled calculus by Gubinelli \& Imkeller \& Perkowski \cite{gubinelli2015paracontrolled} or renormalization group techniques by Kupiainen \& Duch \cite{kupiainen2016renormalization,duch2025flow}. More recently, singular SPDEs in inhomogeneous but deterministic media have been investigated in \cite{singh2025canonical,broux2025renormalised}, see also \cite{hairer2025periodic,chen2024periodic,chen2023periodic} for periodic homogenization results. In this work, we aim to advance the theory of sub-critical singular SPDEs to the setting of random inhomogeneous media that may be correlated with the driving noise (see Assumption \ref{AssumptionCoef}). As motivating examples, we consider the $\varphi^4_2$-equation in correlated environment, a reaction-diffusion equation symbolically given by 
\begin{equation}\label{QuantizationEquation}
\partial_t u-\nabla\cdot a\nabla u=-u^3+\xi\quad\text{in $(0,\infty)\times \mathbb{R}^2$},
\end{equation}
where $\xi$ denotes a space-time white noise and $a : \mathbb{R}^{1+2}\rightarrow\mathbb{R}^{2\times 2}$ takes values in the space of elliptic and bounbed matrices. 

\medskip

The $\varphi_2^{4}$-equation \eqref{QuantizationEquation} with constant coefficients (i. e. $a=\mathrm{Id}$) is a very classical toy-model and has been largely studied in the physics and mathematics literature. It has been first introduced by Parisi \& Wu in \cite{parisi1981perturbation} as the Langevin dynamic of the associated $\varphi^4_2$-Gibbs measure (see also \cite{damgaard1987stochastic, glimm1987quantum}). Mathematically, it has been first solved by Da Prato \& Debussche in \cite{da2003strong} via a fixed point approach in suitable Besov spaces, showing local in time existence. See also the work of Mourrat \& Weber \cite{mourrat2017global} where they show global in time existence as well as its generalization by Tsatsoulis \& Weber \cite{tsatsoulis2019spectral}. Beyond its Langevin dynamic interpretation, the $\varphi^{4}_2$-equation also describes, for instance, the dynamic of thermodynamic fluctuations in two-dimensional ferromagnet near the critical temperature (see for instance \cite{zinn2002quantum,mourrat2017convergence,brydges1982random}).
The case of non-constant coefficients can be seen as local impurities or random bonds in the material, which we allow to be correlated to the external random forcing. The mathematical study of \eqref{QuantizationEquation} in the non-constant coefficient setting is more recent, see for instance the work of Singh \cite{singh2025canonical} and Broux \& Singh \& Steele \cite{broux2025renormalised} for a general theory (which covers much more than \eqref{QuantizationEquation}); as well as Singh \& Hairer \cite{hairer2025periodic,hairer2025homogenisation} and Chen \& Xu \cite{chen2023periodic} for periodic homogenization results.

\medskip

In this contribution, we introduce an inductive approach to obtain stochastic estimates for the renormalized model associated with \eqref{QuantizationEquation} in correlated random media (i. e. when $a$ is correlated to $\xi$). The case of random coefficients for \eqref{QuantizationEquation} has been recently studied by the author \& Singh in \cite{clozeau2025renormalisation} where we established stochastic estimates for the renormalised model using tools from probability (in the form of finite-dimensional Gaussian integration by parts) and the Hairer-Quastel criterion \cite[Appendix A]{hairer2018class} for bounds on general convolutions. In this contribution, we revisit the results of \cite{clozeau2025renormalisation} and present an alternative proof with a more elementary approach which does not rely on the machinery of \cite{hairer2018class} (and thus diagram-free as well). In addition, it has the advantage of relaxing the regularity assumption on the correlations, which we only require to be in $\cc^\alpha(\mathbb{R}^{1+2})$ (for some $\alpha\in (0,1)$, together with integrability conditions, see Assumption \ref{AssumptionCoef}) instead of $\cc^3(\mathbb{R}^{1+2})$. More precisely, we prove the stochastic estimates for the renormalized model by induction on the moment's exponent where the main novelty is to use infinite-dimensional Gaussian calculus (we recall the framework in Appendix \ref{AppendixGIBP}) in the spirit of recent works applying Malliavin calculus to singular SPDEs (see \cite{linares2024diagram,tempelmayr2024characterizing,hairer2024bphz}) but with a very different methodology due to the correlations (see the discussion in Section \ref{SecStrat}). In this paper we only focus on establishing the stochastic estimates, which is the non-trivial part of the analysis; convergence of the model can be obtained following a very similar strategy and a solution theory can be derived following the arguments in the literature \cite[Proposition $4.4$]{da2003strong}, \cite[Theorem $6.2$]{mourrat2017global}, \cite[Theorem,$3.3$]{tsatsoulis2019spectral} or \cite[Proposition $3.2$]{hairer2025periodic}. 
Finally, this work represents a first step towards understanding the stochastic homogenization of \eqref{QuantizationEquation}, a question we leave for future investigations with the present result providing a general strategy for obtaining stochastic estimates.


%
\section{Notations}
Throughout the article, we use the following notational conventions:

\medskip

\textbf{Abbreviations. }We generically write $G_{\sigma}$ (with $\sigma>0$) for a scalar centred Gaussian with variance $\sigma^2$, i. e. 
$$G_{\sigma}(x):=\frac{1}{2\pi\,\sigma^2}\exp\Big(-\tfrac{\vert x\vert^2}{2\sigma^2}\Big)\quad\text{for any $x\in \mathbb{R}^2$,}$$
and for any $C>0$, we define
\begin{equation}\label{DeriveGaussianGamma}
H_{C}(z):=\mathds{1}_{t > 0}\,G_{\sqrt{t}}\Big(\tfrac{1}{\sqrt{C}}x\Big)\quad\text{for any $z=(t,x)\in \mathbb{R}^{1+d}$.}
\end{equation}

We generically denote an element of $\mathbb{R}^{1+2}$ by $z=(t,x)\in \mathbb{R}\times \mathbb{R}^2$, and keep the suggestive notation $z'=(t',x')$, $\tilde{z}=(\tilde{t},\tilde{x})$, etc.

\medskip

In order to keep long integral expressions lighter, we will also write $\int \dd z_1\cdots \dd z_n$ for $\int_{\mathbb{R}^{1+2}}\dd z_1\cdots\int_{\mathbb{R}^{1+2}}\dd z_n$. 

\medskip

\textbf{Differential calculus. }For $X$ a Banach space and $F : X\rightarrow \mathbb{R}$ we write, for any $n\geq 1$, $\dd^n F$ for the $n$-th Fréchet derivative of $F$. 

\medskip

\textbf{Geometry. }For any $r>0$ and $z=(t,x)\in \mathbb{R}^{1+2}$, we write $\cc_r(z)=[t,t+r^2)\times \bb_r(x)$ for the parabolic cylinder of radius $r$ and centered at $z$. We denote by $\vert \cdot\vert$ the parabolic distance in $\mathbb{R}^{1+2}$ given by $\vert z\vert^2:=\sqrt{t}+\vert x\vert^2$ for any $z=(t,x)$.

\medskip

\textbf{Functional spaces. }For any $\alpha\in (0,1)$, $r>0$ and $y\in \mathbb{R}^{1+2}$, we denote by $\cc^\alpha(\cc_r(y))$ the space of H\"older continuous functions in $\cc_r(y)$  with respect to the parabolic distance. We denote by $|\cdot|_{\cc^{\alpha}(\cc_r(y))}$ the parabolic H\"older semi-norm defined as for $f : \mathbb{R}^{1+2}\rightarrow \mathbb{R}^k$ (for any $k\geq 1$)
$$|f|_{\cc^{\alpha}(\cc_r(y))}:=\sup_{z,z'\in \cc_r(y)}\frac{|\vert f(z)-f(z')|\vert}{\vert z-z'\vert^\alpha}.$$

\medskip

We denote by $\bb^\alpha_{2,\infty}(\mathbb{R}^{1+2})$ the Besov space endowed with the norm for any $f : \mathbb{R}^{1+2}\rightarrow \mathbb{R}$
$$\|f\|_{\bb^\alpha_{2,\infty}(\mathbb{R}^{1+2})}:=\sup_{z\neq 0}\frac{\|f(\cdot+z)-f\|_{\LL^2(\mathbb{R}^{1+2})}}{\vert z\vert^\alpha}+\|f\|_{\LL^2(\mathbb{R}^{1+2})}.$$
\section{Assumptions and statement of the main result}
We first state the assumptions on our coefficient field $a$. 
\begin{assumption}[Assumptions on the coefficient field]\label{AssumptionCoef}
Let $\xi\in \mathcal{S}'(\mathbb{R}^{1+2})$ be a scalar space-time white noise. We assume that there exist $m\in \cc^{\alpha}(\mathbb{R}^{1+2})\cap \bb^\alpha_{2,\infty}(\mathbb{R}^{1+2})$ satisfying
\begin{equation}\label{AssIntM}
\int_{\mathbb{R}^{1+2}}\dd z\, \big(1+\vert z\vert^2\big)^{\sigma}\vert m(z)\vert<\infty,
\end{equation}
for given $\alpha \in (0,1)$ and $\sigma>\frac{1}{2}$, as well as $A : \mathbb{R}\rightarrow \mathbb{R}^{2\times 2}$ such that 

\medskip

\begin{itemize}
\item[(i)]There exists $\lambda\in (0,1)$ such that $A$ takes values in the set of smooth and symmetric $\lambda$-uniformly elliptic matrices, i. e. for any $x\in \mathbb{R}$, $k\geq 1$ and $\zeta\in \mathbb{R}^2$,
$$\lambda\vert \zeta\vert^2\leq A(x)\zeta\cdot \zeta\leq \vert \zeta\vert^2\quad \text{and}\quad \sup_{x\in \mathbb{R}}\Big\vert \frac{\dd^k}{\dd x^k}A(x)\Big\vert<\infty;$$
%
%
%
%
%

%
\medskip

\item[(ii)]$a$ is correlated to the driving noise $\xi$, i. e. 
\begin{equation}\label{CoefField}
a:=A(g)\quad\text{with $g=m\star \xi$.}
\end{equation}
\end{itemize}
\end{assumption}
\begin{remark}
Note that the assumption $m\in \bb^\alpha_{2,\infty}(\mathbb{R}^{1+2})$ implies that the coefficient field $a$ satisfies for any $\alpha'<\alpha$, $k\geq 0$, $r>0$ and $\varepsilon>0$
\begin{equation}\label{RegCoefAHolder}
\sup_{z\in \mathbb{R}^{1+2}}\mathbb{E}^\frac{1}{p}\Big[|a^{(k)}|^p_{\cc^{\alpha'}(\cc_r(z))}\Big]\lesssim \max\big\{1,r^\varepsilon\big\}\quad\text{for any $p<\infty$},
\end{equation}
where $a^{(k)}:=A^{(k)}(g)$. Estimate \eqref{RegCoefAHolder} follows from the covariance sturcture: for any $z,z'\in \mathbb{R}^{1+2}$
$$\mathbb{E}\Big[\vert g(z)-g(z')\vert^2\Big]=\|m(z-\cdot)-m(z'-\cdot)\|^2_{\LL^2(\mathbb{R}^{1+2})}\leq \|m\|^2_{\bb^{\alpha}_{2,\infty}(\mathbb{R}^{1+2})}\vert z-z'\vert^{2\alpha},$$
which we can upgrade by Gaussianity into arbitrary moments
$$\mathbb{E}^\frac{1}{p}\Big[\vert g(z)-g(z')\vert^p\Big]\lesssim_p \|m\|_{\bb^{\alpha}_{2,\infty}(\mathbb{R}^{1+2})} \vert z-z'\vert^\alpha\quad \text{for any $p<\infty$}.$$
We then conclude by the regularity assumption on $A$ and a standard Kolmogorov continuity criterion. 
%
\end{remark}
We now turn to the statement of our main result. To renormalize the equation \eqref{QuantizationEquation}, we follow the so-called \textit{Da Prato-Debussche trick} in \cite{da2003strong} that we briefly reproduce in the following. The idea is to define the auxiliary quantity $\lollipop$ which solves the linear part of the equation
$$\partial_t \lollipop-\nabla\cdot a\nabla\,\lollipop=\xi\quad\text{in $(0,\infty)\times \mathbb{R}^2$},$$
and declare that $u$ is a solution of \eqref{QuantizationEquation} if and only if $u=w+\lollipop$ where $w$ symbolically solves 
\begin{equation}\label{DDEquation}
\partial_t w-\nabla\cdot a\nabla w=-w^3-3 w^2\, \lollipop-3 w\,\cherry -\chickenfoot\quad\text{in $(0,\infty)\times \mathbb{R}^2$},
\end{equation}
where we have set $\cherry:=\lollipop^2$ and $\chickenfoot:=\lollipop^3$. Note that \eqref{DDEquation} is more regular than \eqref{QuantizationEquation} since, provided we can make sense of $\cherry$ and $\chickenfoot$ as distributions in $\cc^{-0^{-}}$, it now holds that $(\partial_t-\nabla\cdot a\nabla)^{-1}\chickenfoot\in \cc^{2^{-}}$ which makes the product $w\cherry$ well-defined by Young's multiplication and thus \eqref{DDEquation} solvable in $\cc^{0^{+}}$.

\medskip

Rigourously, we regularise the noise $\xi$ with a smooth mollifier: for any $\delta>0$ and for some even function $\psi\in \cc^{\infty}_c(\mathbb{R})$, we set
\begin{equation}\label{RegKernel}
\rho_\delta(t,x):=\psi_\delta(t)G_{\delta}(x)\quad\text{ with $\psi_\delta(t):=\delta^{-2}\psi(\tfrac{t}{\delta^2})$ for any $(t,x)\in \mathbb{R}^{1+2}$},
\end{equation}
and $\xi_\delta:=\rho_\delta\star \xi$. We then define the regularised $\lollipop_\delta$ which solves 
$$
\left\{
    \begin{array}{ll}
        \partial_t \lollipop_\delta-\nabla \cdot a\nabla\,\lollipop_\delta=\xi_\delta & \text{in $(0,\infty)\times \mathbb{R}^2$,} \\
        \lollipop_\delta(0)=0. & 
    \end{array}
\right.
$$
We can now safely take the products $\cherry_\delta:=\lollipop^2_\delta$ and $\chickenfoot_\delta:=\lollipop^3_\delta$ and perform a renormalization. 

\medskip

Following \cite{singh2025canonical}, we decompose the kernel $\Gamma$ of the parabolic operator $\partial_t-\nabla\cdot a\nabla$ into a singular part and a regular part: for any $z,z'\in \mathbb{R}^{1+2}$, we write 
\begin{equation}\label{RegPartDef}
\Gamma(z,z')=K_{z}(z-z')+R(z,z'),
\end{equation}
where $K_{z}$ is the heat-kernel associated to the linear operator (with \textit{frozen coefficients}) $-\nabla\cdot a(z)\nabla$ given by 
\begin{equation}\label{DefKernelFrozen}
K_{z}(s,y):=\mathds{1}_{\mathbb{R}^+}(s)\frac{1}{4\pi s\sqrt{\text{det}(a(z))}}\exp\bigg(-\frac{y\cdot a^{-1}(z)y}{4s}\bigg)\quad \text{for any $(s,y)\in (0,\infty)\times\mathbb{R}^2$}.
\end{equation}
We then decompose $\cherry_\delta$ via Green's function representation and by expanding the square as
\begin{align*}
\cherry_\delta(z)=&\Big(\underbrace{\int_{\mathbb{R}^{1+2}} \dd z'\, R(z,z')\xi_\delta(z')}_{:=\mathcal{R}_\delta(z)}\Big)^2+2\Big(\int_{\mathbb{R}^{1+2}} \dd z'\, R(z,z')\xi_\delta(z')\Big)\underbrace{\Big(\int_{\mathbb{R}^{1+2}} \dd z'\, K_{z}(z-z')\xi_\delta(z')\Big)}_{:=\widehat{\lollipop}_\delta(z)}\\
&+\underbrace{\Big(\int_{\mathbb{R}^{1+2}} \dd z'\, K_{z}(z-z')\xi_\delta(z')\Big)^2}_{:=\widehat{\cherry}_\delta(z)},
\end{align*}
and likewise 
$$\chickenfoot_\delta=\big(\mathcal{R}_\delta(z)\big)^3+3\big(\mathcal{R}_\delta(z)\big)^2\,\widehat{\lollipop}_\delta(z)+2\mathcal{R}_\delta(z)\widehat{\cherry}_\delta(z)+\widehat{\chickenfoot}_\delta(z).$$
As observed in \cite{singh2025canonical}, the regular parts $\{(\mathcal{R}_\delta)^m\}_{m\in \{1,2,3\}}$ are uniformly bounded in $\delta>0$ (we provide a proof in Proposition \ref{PropRegPart}) and $\widehat{\lollipop}_\delta$ is as well controlled uniformly in $\delta>0$ as a distribution in $\cc^{0^-}(\mathbb{R}^{1+2})$. Therefore, only the parts $\widehat{\cherry}_\delta$ and $\widehat{\chickenfoot}_\delta$ need a renormalization. Despite the fact that the law of the coefficient field $a$ is stationary, which may suggest to use constants for the renormalization procedure, it has been observed in \cite[Proposition $1.2$]{clozeau2025renormalisation} that the use of deterministic renormalization functions generically lead to variance blow-up. We thus use random renormalization functions and define:
$$c^{\cherryS}_\delta(z):=\int_{\mathbb{R}^{1+2}}\dd z'\, \int_{\mathbb{R}^{1+2}}\dd z''\,K_z(z-z')K_z(z-z'')\rho_\delta\star \rho_\delta(z'-z''),$$
and 
$$c^{\chickenfootS}_\delta(z):=3c^{\cherryS}_\delta(z)\,\widehat{\lollipop}_\delta(z);$$
as well as the renormalized products
\begin{align*}
\overline{\cherry}_\delta(z):=\widehat{\cherry}_\delta(z)-c^{\cherryS}_\delta(z)=\int_{\mathbb{R}^{1+2}}\dd z'\int_{\mathbb{R}^+\times \mathbb{R}^2}\dd z''\, K_{z}(z-z')K_{z}(z-z'')\xi_\delta(z')\diamond \xi_\delta(z''),
\end{align*}
likewise
\begin{align*}
\overline{\chickenfoot}_\delta(z)&:=\widehat{\chickenfoot}_\delta(z)-c^{\chickenfootS}_\delta(z)\\
&=\int_{\mathbb{R}^{1+2}}\dd z'\int_{\mathbb{R}^{1+2}}\dd z''\, \int_{\mathbb{R}^{1+2}}\dd z'''\,K_{z}(z-z')K_{z}(z-z'')K_{z}(z-z''')\xi_\delta(z')\diamond \xi_\delta(z'')\diamond \xi_\delta(z'''),
\end{align*}
where we recall that the Wick product $\diamond$ is defined in Appendix \ref{AppendixGIBP}. We now state our main result that establishes the stochastic estimates for $\{\widehat{\lollipop}_\delta\}_{\delta>0}$, $\{\overline{\cherry}_\delta\}_{\delta>0}$ and $\{\overline{\chickenfoot}_\delta\}_{\delta>0}$.
\begin{theorem}[Stochastic estimates of the renormalized model]\label{MainResult}
For any $\lambda\in (0,e^{-1}]$ and test functions $\psi^\lambda:=\lambda^{-4}\psi(\frac{\cdot}{\lambda^2},\frac{\cdot}{\lambda})$ with $\mathrm{supp }\,\psi\subset \cc_1$ and $\|\psi\|_{\cc^0(\mathbb{R}^{1+2})}\leq 1$, it holds uniformly in $\delta\in (0,1)$: for any $p<\infty$
\begin{equation}\label{MainResultStoEsti}
\mathbb{E}^\frac{1}{p}\Big[\big\vert\big(\,\widehat{\lollipop}_\delta,\psi^\lambda\big)\big\vert^p\Big]\lesssim \vert \log(\lambda)\vert^\frac{1}{2},\quad\mathbb{E}^\frac{1}{p}\Big[\big\vert\big(\overline{\cherry}_\delta,\psi^\lambda\big)\big\vert^p\Big]\lesssim \vert \log(\lambda)\vert\quad\text{and}\quad \mathbb{E}^\frac{1}{p}\Big[\big\vert\big(\overline{\chickenfoot}_\delta,\psi^\lambda\big)\big\vert^p\Big]\lesssim \vert \log(\lambda)\vert^{\frac{3}{2}}.
\end{equation}
\end{theorem}
\begin{remark}
Note that $\widehat{\lollipop}_\delta$, $\overline{\cherry}_\delta$ and $\overline{\chickenfoot}_\delta$ are stationary random fields. Indeed, for any $z,\star \in \mathbb{R}^{1+2}$ it holds 
\begin{align*}
\overline{\cherry}_\delta(\xi,z+\star)&=\int \dd z'\int \dd z''\, K_{z+\star}(z+\star-z')K_{z+\star}(z+\star-z'')\xi_\delta(z')\diamond\xi_\delta(z'')\\
&=\int \dd z'\int \dd z''\, K_{z+\star}(z-z')K_{z+\star}(z-z'')\xi_\delta(z'+\star)\diamond\xi_\delta(z''+\star)\\
&\stackrel{\eqref{DefKernelFrozen}}{=}\overline{\cherry}_\delta\big(\xi(\cdot+\star),z\big),
\end{align*}
likewise for $\overline{\chickenfoot}_\delta$ and $\widehat{\lollipop}_\delta$. Therefore, by stationarity, \eqref{MainResultStoEsti} also holds for any base point $\star\in \mathbb{R}^{1+2}$ and test functions $\psi^\lambda(\cdot-\star)$. As a consequence of a Komolgorov continuity type argument, Theorem \ref{MainResult} implies tightness of $\{\widehat{\lollipop}_\delta\}_{\delta>0}$, $\{\overline{\cherry}_\delta\}_{\delta>0}$ and $\{\overline{\chickenfoot}_\delta\}_{\delta>0}$ in $\cc^{-0^{-}}$ (a proof can be done following the lines of \cite[Theorem $10.7$]{hairer2014theory}). 
\end{remark}
\section{Strategy of the proof}\label{SecStrat}
In this section, we present and show how to implement the inductive approach to the specific model \eqref{QuantizationEquation}. At first glance, one might be tempted to apply the spectral gap inequality approach developed in \cite{linares2024diagram,hairer2024bphz} to get the stochastic estimates in Theorem \ref{MainResult}. Unfortunately, since the kernel $K_z$ is correlated with $\xi_\delta$ through $a$, the spectral gap inequality is not fine enough to capture the fluctuations of $\big(\widehat{\lollipop}_\delta,\psi^\lambda\big)$, $\big(\overline{\cherry}_\delta,\psi^\lambda\big)$ and $\big(\overline{\chickenfoot}_\delta,\psi^\lambda\big)$: indeed, when the Malliavin derivative acts on one of the kernel $K_z$ it produces a term that is typically of the same type which, unless being in a weakly correlated setting, prevents us from closing the estimate. Instead, we need finer expressions that account for these correlations more carefully. For the exposition, we use the framework of Gaussian calculus recalled in Appendix \ref{AppendixGIBP} and the notation defined in Definition \ref{QuantityNorm}. 

\medskip

First, we identify three main quantities on which to run the induction that we denote by $\{\mathbb{X}_n\}_{n\leq 3}$. Looking at the structure of $\big(\widehat{\lollipop}_\delta,\psi^\lambda\big)$, $\big(\overline{\cherry}_\delta,\psi^\lambda\big)$ and $\big(\overline{\chickenfoot}_\delta,\psi^\lambda\big)$, a natural choice is 
\begin{equation}\label{QuantityInducStratProof}
\mathbb{X}_n:\Gamma\mapsto \int \dd z_1\cdots\dd z_{n+1}\, \Gamma(z_1,\cdots,z_{n+1})\, \mathop{\diamond}_{i=2}^{n+1}\xi_\delta(z_i),
\end{equation}
where $\Gamma$ is a general kernel satisfying heat-kernel bounds, see Definition \ref{DefClassKernel} for more details. Typically, for any $z_1,\cdots,z_{n+1}\in \mathbb{R}^{n+1}$
\begin{equation}\label{TypiclKernel}
\Gamma(z_1,\cdots,z_{n+1})\sim \prod_{i=2}^{n+1}K_{z_1}(z_1-z_i).
\end{equation}
Second, the idea is to use the integration by parts formula in Lemma \ref{IntWRTWick} to establish a recursive inequality between moments of $\mathbb{X}_i$ and that of $\{\mathbb{X}_{n}\}_{n\leq i}$ for any $i\leq 3$. To do so, we use that: for any $p\geq 1$ and any kernel $\Gamma$
$$\mathbb{E}\big[(\mathbb{X}_i(\Gamma))^{2p}\big] =\int \dd z_1\cdots \dd z_{i+1}\, \mathbb{E}\Big[(\mathbb{X}_i(\Gamma))^{2p-1}\Gamma(z_1,\cdots,z_{i+1})\mathop{\diamond}_{j=2}^{i+1}\xi_\delta(z_j)\Big].$$
and use integration by parts \eqref{GaussianInteNnoises} to obtain 
$$\mathbb{E}\big[(\mathbb{X}_i(\Gamma))^{2p}\big]=\int \dd z_1\cdots \dd z_{i+1}\, \mathbb{E}\Big[\mathrm{D}_{\delta,z_1,\cdots,z_{i+1}}\big((\mathbb{X}_n(\Gamma))^{2p-1}\Gamma(z_1,\cdots,z_{i+1})\big)\Big].$$
Expanding the derivative, we can deduce a recursive inequality for $\mathbb{E}\big[(\mathbb{X}_n(\Gamma))^{2p}\big]$ involving $\{\mathbb{X}_n\}_{n\leq 3}$, each of them evaluated at possibly different kernels. The different kernels involve natturally derivatives of $\Gamma$ for which heat-kernel bounds should be stable (see the typical kernel \eqref{TypiclKernel}), that is made precise via the class of kernels in Definition \ref{DefClassKernel}. The recursive inequalities are shown in details in Lemma \ref{LinkMalliavinDerive}. Finally, once recursive inequalities are established, we bound the quantities $\{\mathbb{X}_n\}_{n\leq 3}$ in increasing order where we can run an induction on the moment's exponent and derive the estimate of Theorem \ref{MainResult}, see Proposition \ref{EstimateForOneNoise}, Proposition \ref{X2Estimates} and Proposition \ref{EstiX3}. 

\medskip

In principle, the strategy presented here is fully automated and can be applied to more singular SPDEs requiring more advanced renormalisation techniques, such as the theory of regularity structures \cite{hairer2014theory} or paracontrolled calculus \cite{gubinelli2015paracontrolled}. The key step is to identify the key quantities $\{\mathbb{X}_n\}_{n\leq N}$, where $N\geq 1$ denotes the number of noises involved, which can be revealed via Gaussian integration by parts. In this paper, we chose to focus on the simple case of \eqref{QuantizationEquation} to make the computations explicit and the strategy easier to reproduce. Generalizations of the method is left for future investigations.

\section{Auxiliary lemmas and definitions}
In this section, we state auxiliary lemmas that will be useful to prove the stochastic bounds of Theorem \ref{MainResult}. Our first result is a well known integration by parts formula with respect to Wick products of Gaussian random variables. Such formulas are classical in the literature,  see for instance \cite{Nualart2006}, and here we state them in a very particular form that will be convenient for us for establishing the stochastic estimates. For completeness, we also provide a proof in Appendix \ref{AppendixGIBP}.
\begin{lemma}[Gaussian integration by parts]\label{IntWRTWick}
With the notations of the Appendix \ref{AppendixGIBP}, we consider $F : X \rightarrow \mathbb{R}$ infinitely many Fréchet differentiable in $X$ and which satisfies the following assumptions: 

\medskip

\begin{itemize}
\item[(i)]There exists $p\geq 1$ and $C>0$ such that for any $T,T'\in X$, 
\begin{equation}\label{AssumptionInteOnF}
\vert F(T)\vert\leq C\big(\|T\|_X+1\big)^p\quad\text{and}\quad \vert F(T)-F(T')\vert\leq C \big(\|T\|_{X}+\|T'\|_{X}+1\big)^p \|T-T'\|_{X};
\end{equation}
\item[(ii)]For any $n\geq 1$, there exists $p_n\geq 1$ and $C_n>0$ such that for any $T,T'\in X$,
\begin{equation}\label{AssumptionInteOnDNF}
|||\dd^n F(T)||| \leq C_n \big(\|T\|_{X}+1\big)^{p_n}\quad\text{and}\quad |||\dd^n F(T)-\dd^n F(T')||| \leq C_n \big(\|T\|_{X}+\|T'\|_{X}+1\big)^{p_n} \|T-T'\|_{X},
\end{equation}
\end{itemize} 
where $|||\cdot|||$ denotes the operator norm. Then for any $n\geq 1$ and $\{\varphi_i\}_{i\leq n}\subset\cc^\infty_c(\mathbb{R}^{1+2})$, it holds
\begin{equation}\label{GaussianInteNnoises}
\mathbb{E}\bigg[F(\xi)\, \mathop{\diamond}_{i=1}^n (\xi,\varphi_i)\bigg]=\int \dd y_1\cdots \dd y_n \, \mathbb{E}\big[\mathrm{D}_{y_1}\cdots \mathrm{D}_{y_n} F(\xi)\big]\varphi_1(y_1)\cdots \varphi_n(y_n),
\end{equation}
where $\mathrm{D}_{y_1}\cdots \mathrm{D}_{y_n} F\in \LL^2((\mathbb{R}^{1+d})^n)$ denotes the kernel of $\dd^n F$ in $\LL^2((\mathbb{R}^{1+d})^n)$. 
\end{lemma}
\begin{remark}
The fact that, for any Fréchet differentiable $F : X\rightarrow \mathbb{R}$, $\dd^n F$ (with $n\geq 1$) admits a square integrable kernel $\mathrm{D}_{y_1}\cdots \mathrm{D}_{y_n} F$ in $\LL^2((\mathbb{R}^{1+2})^n)$ comes from the structure of $X$. We recall that $\|\cdot\|_X=\|S\cdot\|_{\LL^{2}(\mathbb{R}^{1+2})}$ where $S : \LL^{2}(\mathbb{R}^{1+2})\rightarrow \LL^{2}(\mathbb{R}^{1+2})$ is an Hilbert-Schmidt operator. Thus, given an orthonormal basis $\{h_j\}_{j\geq 1}$ of $\LL^{2}(\mathbb{R}^{1+2})$, we can check that for any $T\in X$ and $y_1,\cdots,y_n\in \mathbb{R}^{1+2}$
\begin{equation}\label{FormulaKernelL2}
\mathrm{D}_{y_1}\cdots \mathrm{D}_{y_n} F(T)=\sum_{j_1,\cdots,j_n=1}^\infty \dd^n F(T).(h_{j_1},\cdots,h_{j_n})h_{j_1}(y_1)\cdots h_{j_n}(y_n).
\end{equation} 
Indeed, on the one hand \eqref{FormulaKernelL2} is a well defined function of $\LL^{2}((\mathbb{R}^{1+2})^n)$: by orthogonality of $\{h_j\}_{j\geq 1}$ and since $S$ is Hilbert-Schmidt,
\begin{align*}
\int\dd y_1\cdots\dd y_n\, \big\vert \mathrm{D}_{y_1}\cdots \mathrm{D}_{y_n} F(T)\big\vert^2 &= \sum_{j_1,\cdots,j_n=1}^\infty\int\dd y_1\cdots\dd y_n\, \big\vert  \dd^n F(T).(h_{j_1},\cdots,h_{j_n})h_{j_1}(y_1)\cdots h_{j_n}(y_n)\big\vert^2\\
&\leq |||\dd^n F(T)|||\Big(\sum_{j=1}^\infty \|h_j\|^2_{X}\Big)^n=|||\dd^n F(T)|||\Big(\sum_{j=1}^\infty \|S h_j\|^2_{\LL^2(\mathbb{R}^{1+2})}\Big)^n<\infty.
\end{align*}
On the other hand, it is now immediate from the $\LL^2((\mathbb{R}^{1+2})^n)$-convergence of \eqref{FormulaKernelL2} and the multi-linearity of $\dd^n F(T)$ that for any $T_1,\cdots,T_n\in \LL^2(\mathbb{R}^{1+2})$,
\begin{align*}
\int\dd y_1\cdots\dd y_n\,  \mathrm{D}_{y_1}\cdots \mathrm{D}_{y_n} F(T)\,T_1(y_1)\cdots T_n(y_n)&=\dd^n F(T).\Big( \sum_{j=1}^\infty (T_1,h_j)_{\LL^2(\mathbb{R}^{1+2})}h_j,\cdots,\sum_{j=1}^\infty (T_n,h_j)_{\LL^2(\mathbb{R}^{1+2})}h_j\Big)\\
&=\dd^n F(T).(T_1,\cdots,T_n).
\end{align*}
\end{remark}
Next, we show a general bound on a convolution type quantity that appears several times in the proof of Theorem \ref{MainResult}.
\begin{lemma}[Convolution bound]\label{LemmaConvoBound}
For any $C>0$, $z,\tilde{z}\in \mathbb{R}^{1+2}$ and $\delta>0$, we define 
$$\mathcal{C}_{\delta,C}(z,\tilde{z}):=\int \dd z_1\, \Big(\int \dd z_2\, H_C(z-z_2)\rho_\delta\star\rho_\delta(z_2-z_1)\Big)H_C(\tilde{z}-z_1),$$
where we recall that $H_C$ is defined in \eqref{DeriveGaussianGamma}. Then, for any $T>0$ and $z,\tilde{z}\in [0,T]\times \mathbb{R}^2$, it holds uniformly in $\delta\in (0,1)$
\begin{equation}\label{BoundRandomTestFunction}
\mathcal{C}_{\delta,C}(z,\tilde{z})\lesssim_T 1+\vert \log(\vert x-\tilde{x} \vert)\vert.
\end{equation}

\end{lemma}
\begin{proof}
The proof of \eqref{BoundRandomTestFunction} is based on the semi-group property for Gaussians: for any $\sigma_1,\sigma_2>0$, $G_{\sigma_1}\star G_{\sigma_2}=G_{\sqrt{\sigma^2_1+\sigma^2_2}}$. Using our choice of regularization kernel $\rho_\delta$ in \eqref{RegKernel}, we have for any $z_1,z_2\in \mathbb{R}^{1+2}$
$$\rho_\delta\star \rho_\delta(z_2-z_1)=\psi_\delta\star \psi_\delta(t_2-t_1)G_{\delta\sqrt{2}}(x_2-x_1).$$
Now, using in addition the definition \eqref{DeriveGaussianGamma} of $H_C$, we have
\begin{align*}
\int \dd z_2\, H_C(z-z_2)\rho_\delta\star \rho_\delta(z_2-z_1)& = C\int_{0}^t \dd t_2\,\psi_\delta\star\psi_\delta (t_2-t_1)G_{\sqrt{C(t-t_2)}}\star G_{\delta\sqrt{2}}(x-x_1)\\
&=C\int_0^t \dd t_2\, \psi_\delta\star \psi_\delta(t_2-t_1)G_{\sqrt{C(t-t_2)+2\delta^2}}(x-x_1).
\end{align*}
Integrating against $\int \dd z_1\, H_C(\tilde{z}-z_1)$ and applying again the semi-group property, we obtain
\begin{align*}
&\int_{0}^{t}\dd t_2\int \dd z_1\, \psi_\delta\star \psi_\delta(t_2-t_1)H_C(\tilde{z},z_1)G_{\sqrt{C(t-t_2)+2\delta^2}}(x-x_1)\\
&\lesssim \int_{0}^t \dd t_2\int_{0}^{\tilde{t}}\dd t_1\, \psi_\delta\star \psi_\delta(t_2-t_1)G_{\sqrt{C(\tilde{t}-t_1)}}\star G_{\sqrt{C(t-t_2)+2\delta^2}}(x-\tilde{x})\\
&= \int_{0}^t \dd t_2\int_{0}^{\tilde{t}}\dd t_1\, \psi_\delta\star \psi_\delta(t_2-t_1)G_{\sqrt{C(\tilde{t}-t_1)+C(t-t_2)+2\delta^2}}(x-\tilde{x}).
\end{align*}
Finally, using the change of variable $t_1\mapsto C(\tilde{t}-t_1)+C(t-t_2)+2\delta^2$ together with $C(t-t_2)+2\delta^2\leq CT+2$ and $\int_{\mathbb{R}}\dd t_2\, \psi_\delta\star\psi_\delta(t_2)=(\int \psi)^2$, we have
$$ \int_{0}^t \dd t_2\int_{0}^{\tilde{t}}\dd t_1\, \psi_\delta\star \psi_\delta(t_2-t_1)G_{\sqrt{C(\tilde{t}-t_1)+C(t-t_2)+2\delta^2}}(x-\tilde{x})\lesssim \int_{0}^{CT+2}\dd t_1\, G_{\sqrt{t_1}}(x-\tilde{x})\lesssim_T 1+\vert \log(\vert x-\tilde{x}\vert)\vert,$$
which gives \eqref{BoundRandomTestFunction}.
\end{proof}
%
%
%
%
%
%
%
Then, we state and prove a general bound on convolutions of logarithms that will be useful for the proof of Proposition \ref{X2Estimates}.
\begin{lemma}[Bound on general convolutions of logarithms]\label{GeneLog}
Let $\tilde{x}\in \mathbb{R}^2$ and, for $n\geq 1$, $\{\alpha_{ij}\}_{1\leq i\neq j\leq n}$ a finite set of integers. Define for any $\lambda\in (0,e^{-1}]$
$$\ell_n:=\fint_{\bb_\lambda}\dd x_1\cdots \dd x_n\, \prod_{i=1}^n \big(1+\vert \log(\vert x_i-\tilde{x}\vert)\vert\big)\prod_{i\neq j}\big(1+\vert \log(\vert x_i-x_j\vert)\vert\big)^{\alpha_{ij}}.$$
Then, it holds
$$\ell_n\lesssim_n \vert \log(\lambda)\vert^{\sum_{i\neq j}\alpha_{ij}}\Big(1+\fint_{\bb_\lambda}\dd x\,\vert \log(\vert x-\tilde{x}\vert)\vert^2\Big)^{\frac{n}{2}}.$$
\end{lemma}
\begin{proof}
We apply the Cauchy-Schwarz inequality in the form of
\begin{align*}
\ell_n\lesssim_n \sqrt{\fint_{\bb_\lambda}\dd x_1\cdots \dd x_n\, \prod_{i\neq j}\big(1+\vert \log(\vert x_i-x_j\vert)\vert\big)^{2\alpha_{ij}}}\sqrt{\fint_{\bb_\lambda}\dd x_1\cdots \dd x_n\, \prod_{i=1}^n \big(1+\vert \log(\vert x_i-\tilde{x}\vert)\vert\big)^2},
\end{align*}
where, with the change of variables $x_i\mapsto \frac{x_i}{\lambda}$, we have 
\begin{align*}
\sqrt{\fint_{\bb_\lambda}\dd x_1\cdots \dd x_n\, \prod_{i\neq j}\big(1+\vert \log(\vert x_i-x_j\vert)\vert\big)^{2\alpha_{ij}}}&=\sqrt{\fint_{\bb_1}\dd x_1\cdots \dd x_n\, \prod_{i\neq j}\big(1+\vert \log(\lambda \vert x_i-x_j\vert)\vert\big)^{2\alpha_{ij}}}\\
&\lesssim \vert \log(\lambda)\vert^{\sum_{i\neq j}\alpha_{ij}}.
\end{align*}
Finally, from Fubini's theorem
$$\sqrt{\fint_{\bb_\lambda}\dd x_1\cdots \dd x_n\, \prod_{i=1}^n \big(1+\vert \log(\vert x_i-\tilde{x}\vert)\vert\big)^2}=\Big(\fint_{\bb_\lambda}\dd x\,\big(1+\vert \log(\vert x-\tilde{x}\vert)\vert\big)^2\Big)^{\frac{n}{2}}\lesssim_n \Big(1+\fint_{\bb_\lambda}\dd x\,\vert \log(\vert x-\tilde{x}\vert)\vert^2\Big)^{\frac{n}{2}}.$$
\end{proof}
For notational convenience, as it will appear several times in the statements, we define a random derivative on Malliavin differentiable random variables.
\begin{definition}[Random derivative]\label{QuantityNorm}
For any $k\geq 1$ and $z_1,\cdots,z_k\in \mathbb{R}^{1+2}$, we define the differential operator $\mathrm{D}_{\delta,z_1,\cdots,z_k}$ on $k$-times Malliavin differentiable (in the sense of Appendix \ref{AppendixGIBP}) $F :  X\rightarrow \mathbb{R}$ as
$$\mathrm{D}_{\delta,z_1,\cdots,z_k}F:=\int \dd y_1\cdots \dd y_k\, \rho_\delta(y_1-z_1)\cdots \rho_\delta(y_k-z_k) \mathrm{D}_{y_1}\cdots \mathrm{D}_{y_k} F.$$
\end{definition}
Finally, as we will refer to Lemma \ref{IntWRTWick} several times in the stochastic estimates, we introduce the following class of random kernels. For this, we use the notations in Appendix \ref{AppendixGIBP}.
\begin{definition}[Class of random kernels]\label{DefClassKernel}
For any $k\geq 0$ and measurable function $\eta : \mathbb{R}^{1+2}\rightarrow (0,\infty)$, we define the class $\mathcal{K}_k(\eta)$ of random kernels as the set of $\Gamma : X\times (\mathbb{R}^{1+d})^{k+1} \rightarrow (0,\infty)$ such that 
\begin{itemize}
\item[(i)]There exists $\bar{\Gamma} : X\times (\mathbb{R}^{1+2})^{k+1}\rightarrow (0,\infty)$ such that $\Gamma=\bar{K}(\xi,\cdot)$ and, for any $z_1,\cdots,z_k\in \mathbb{R}^{1+2}$, $\bar{K}(\cdot,z_1,\cdots,z_k)$ satisfies the assumptions of Lemma \ref{IntWRTWick} with constants $C(z_1,\cdots,z_{k+1})$ and $C_n(z_1,\cdots,z_k)$ satisfying for some $C_n,\sigma>0$
\begin{equation}\label{ConstantDefClass}
\max\big\{C(z_1,\cdots,z_{k+1}),C_n(z_1,\cdots,z_{k+1})\big\}\lesssim_n \eta(z_1)(1+\vert z_1\vert^2)^\sigma \prod_{i=2}^{k+1} H_{C_n}(z_1-z_i).
\end{equation}
\item[(ii)]For any $n\geq 1$, there exist $C_n>0$ such that for any $\delta>0$ and $z_1,\cdots,z_k,z'_1,\cdots,z'_n\in \mathbb{R}^{1+2}$, it holds
\begin{equation}\label{AssumptionKernelFunction1Noise}
\begin{aligned}
\vert \Gamma(z_1,\cdots,z_{k+1})\vert+\vert \mathrm{D}_{\delta,z'_1,\cdots,z'_n}\Gamma(z_1,\cdots,z_{k+1})\vert\lesssim \eta(z_1)\prod_{i=2}^{k+1} H_{C_n}(z_1-z_i),
\end{aligned}
\end{equation}
where for convention we have set $\prod_{i=2}^1H_{C_n}(z_1-z_i)=1$.
\end{itemize}
\end{definition}
\begin{remark}
The point $(i)$ in Definition \ref{DefClassKernel} ensures that we can apply Lemma \ref{IntWRTWick} for $\mathbb{X}_1$, $\mathbb{X}_2$ and $\mathbb{X}_3$ defined in \eqref{QuantityInducStratProof}, whereas $(ii)$ is imposed to establish the stochastic estimates.
\end{remark}
\section{Proof of the stochastic estimates}
In this section, we give the proofs of the stochastic estimates leading to Theorem \ref{MainResult}. As already explained in Section \ref{SecStrat}, our main task is to control the quantities $\{\mathbb{X}_i\}_{i\leq 3}$ in increasing order which is the purpose of Lemma \ref{LinkMalliavinDerive}, Proposition \ref{EstimateForOneNoise}, Proposition \ref{X2Estimates} and Proposition \ref{EstiX3}. Finally, we give in Proposition \ref{PropRegPart} the estimates on the regular part $\{\mathcal{R}_\delta\}_{\delta>0}$. 
\subsection{A technical lemma for Malliavin derivatives of $\{\mathbb{X}_n\}_{n\geq 1}$}
The proof of Theorem \ref{MainResult} is based on a recursive argument and exploits the connection between the Malliavin derivatives (up to three orders) of $\mathbb{X}_n$ and the lower-order processes $\{\mathbb{X}_i\}_{i\leq n-1}$, for $n\in \{2,3\}$. This connection is established in the following lemma.
\begin{lemma}\label{LinkMalliavinDerive}
Let $n\geq 3$ and $\eta : \mathbb{R}^{1+2} \rightarrow (0,\infty)$ be measurable. Let $\Gamma_1,\Gamma_2 \in \mathcal{K}_n(\eta)$ symmetric w. r. t. their $n$-lasts variables and define
$$\mathbb{X}_n(\Gamma_1):=\int \dd z_1\cdots \dd z_{n+1} \, \Gamma_1(z_1,\cdots,z_{n+1})\mathop{\diamond}_{i=2}^{n+1} \xi_\delta(z_i),$$
where we recall that $\displaystyle\mathop{\diamond}_{i=2}^{n+1}$ is defined in \eqref{WickProductDef}, as well as for any $z'\in \mathbb{R}^{1+2}$ and $i\in \{1,2\}$
$$ \eta_{i}(z'):= \eta(z')\displaystyle\int \dd z\, \eta(z)\big(1+\vert \log(\vert x'-x\vert)\vert\big)^i.$$
Then, it holds for some $C>0$
\begin{itemize}
\item[(i)] For the first-order Malliavin derivative: there exist $\Gamma^{\mathrm{1st}}_{\delta,1}\in \mathcal{K}_n(\eta\int \eta)$ and $\Gamma^{\mathrm{1st}}_{\delta,2}\in \mathcal{K}_{n-1}(\eta_{1})$ such that
\begin{equation}\label{EquationFirstOrderDeriv}
\begin{aligned}
&\int \dd z_2 \,\mathrm{D}_{\delta,z_2}\big(\mathbb{X}_{n}(\Gamma_1)\big)^{2p-1}\int \dd z_1\, \dd z_3\cdots \dd z_{n+1}\, \Gamma_2(z_1,\cdots, z_{n+1})\\
&=(2p-1)\big(\mathbb{X}_n(\Gamma_1)\big)^{2(p-1)}\big(\mathbb{X}_n(\Gamma^{\mathrm{1st}}_{\delta,1})+\mathbb{X}_{n-1}(\Gamma^{\mathrm{1st}}_{\delta,2})\big).
\end{aligned}
\end{equation}
\item[(ii)] For the second-order Malliavin derivative: there exist $\Gamma^{\mathrm{2nd}}_{\delta,1}\in \mathcal{K}_n(\eta\int \eta)$, $\Gamma^{\mathrm{2nd}}_{\delta,2}\in \mathcal{K}_{n-1}(\eta_{1})$, $\Gamma^{\mathrm{2nd}}_{\delta,3}\in \mathcal{K}_{n-2}(\eta_{2})$ as well as $\mathbb{R}^{1+2}\ni z\mapsto (\Gamma^{\mathrm{2nd}}_{\delta,4}(z),\Gamma^{\mathrm{2nd}}_{\delta,5}(z))\in \mathcal{K}_{n}(\eta)\times \mathcal{K}_{n-1}(\eta \,\mathcal{C}_{\delta,C}(\cdot,z))$ such that
\begin{equation}\label{EquationSecondOrderDeriv}
\begin{aligned}
&\int \dd z_2\, \dd z_3\, \mathrm{D}_{\delta,z_2,z_3}\big(\mathbb{X}_n(\Gamma_1)\big)^{2p-1}\int \dd z_1\, \dd z_4\cdots \dd z_{n+1}\, \Gamma_2(z_1,\cdots, z_{n+1})\\
&\lesssim_p \big(\mathbb{X}_n(\Gamma_1)\big)^{2(p-1)}\big(\mathbb{X}_n(\Gamma^{\mathrm{2nd}}_{\delta,1})+\mathbb{X}_{n-1}(\Gamma^{\mathrm{2nd}}_{\delta,2})\big)+\big(\mathbb{X}_n(\Gamma_1)\big)^{2(p-1)}|\mathbb{X}_{n-2}(\Gamma^{\mathrm{2nd}}_{\delta,3})|\\
&\quad +(p-1)\vert\mathbb{X}_n(\Gamma_1)\vert^{2p-3}\int \dd z\, \eta(z)\, \big(\vert \mathbb{X}_n(\Gamma^{\mathrm{2nd}}_{\delta,4}(z))\vert +\vert \mathbb{X}_{n-1}(\Gamma^{\mathrm{2nd}}_{\delta,5}(z))\vert\big)^2,
\end{aligned}
\end{equation}
where we recall that $\mathcal{C}_\delta$ was defined in Lemma \ref{LemmaConvoBound}.

\medskip

\item[(iii)]For the third Malliavin derivative: there exist $\Gamma^{\mathrm{3rd}}_{\delta,1}\in \mathcal{K}_n(\eta\int \eta)$, $\Gamma^{\mathrm{3rd}}_{\delta,2}\in \mathcal{K}_{n-1}(\eta_{1})$, $\Gamma^{\mathrm{3rd}}_{\delta,3}\in \mathcal{K}_{n-2}(\eta_{2})$, $\mathbb{R}^{1+2}\ni z\mapsto (\Gamma^{\mathrm{3rd}}_{\delta,4}(z),\Gamma^{\mathrm{3rd}}_{\delta,5}(z))\in \mathcal{K}_{n}(\eta)\times \mathcal{K}_{n-1}(\eta \,\mathcal{C}_\delta(\cdot,z))$, $\Gamma^{\mathrm{3rd}}_{\delta,6}\in \mathcal{K}_{n-3}(\eta_3)$, $\mathbb{R}^{1+2}\ni z\mapsto (\Gamma^{\mathrm{3rd}}_{\delta,7}(z),\Gamma^{\mathrm{3rd}}_{\delta,8}(z))\in \mathcal{K}_{n}(\eta)\times \mathcal{K}_{n-1}(\eta \,\mathcal{C}_{\delta,C}(\cdot,z))$, $z\mapsto(\Gamma^{\mathrm{3rd}}_{\delta,9}(z),\Gamma^{\mathrm{3rd}}_{\delta,11}(z))\in \mathcal{K}_n(\eta)\times \mathcal{K}_n(\eta)$, $z\mapsto (\Gamma^{\mathrm{3rd}}_{\delta,10}(z),\Gamma^{\mathrm{3rd}}_{\delta,12}(z),\Gamma^{\mathrm{3rd}}_{\delta,13}(z))\in \mathcal{K}_{n-1}(\eta\,\mathcal{C}_{\delta,C}(\cdot,z))\times \mathcal{K}_{n-1}(\eta\, \mathcal{C}_{\delta,C}(\cdot,z))\times \mathcal{K}_{n-2}(\eta \,\mathcal{C}^2_{\delta,C}(\cdot,z))$ such that
\begin{equation}\label{EquationThirdOrderDeriv}
\begin{aligned}
&\int \dd z_2\, \dd z_3\, \dd z_4\, \mathrm{D}_{\delta,z_2,z_3,z_4}\big(\mathbb{X}_n(\Gamma_1)\big)^{2p-1}\int \dd z_1\, \dd z_5\cdots \dd z_{n+1}\Gamma_2(z_1,\cdots,z_{n+1})\\
&\lesssim_p \big(\mathbb{X}_n(\Gamma_1)\big)^{2(p-1)}\big(\mathbb{X}_n(\Gamma^{\mathrm{3rd}}_{\delta,1})+\mathbb{X}_{n-1}(\Gamma^{\mathrm{3rd}}_{\delta,2})\big)+\big(\mathbb{X}_n(\Gamma_1)\big)^{2(p-1)}|\mathbb{X}_{n-2}(\Gamma^{\mathrm{3rd}}_{\delta,3})|\\
&\quad +(p-1)\vert\mathbb{X}_n(\Gamma_1)\vert^{2p-3}\int \dd z\, \eta(z)\, \big(\vert \mathbb{X}_n(\Gamma^{\mathrm{3rd}}_{\delta,4}(z))\vert +\vert \mathbb{X}_{n-1}(\Gamma^{\mathrm{3rd}}_{\delta,5}(z))\vert\big)^2\\
&\quad + \big(\mathbb{X}_n(\Gamma_1)\big)^{2(p-1)}|\mathbb{X}_{n-3}(\Gamma^{\mathrm{3rd}}_{\delta,6})|\\
&\quad + (p-1)(2p-3)\vert\mathbb{X}_n(\Gamma_1)\vert^{2p-4}\int \dd z\, \eta(z)\, \big(\vert \mathbb{X}_n(\Gamma^{\mathrm{3rd}}_{\delta,7}(z))\vert +\vert \mathbb{X}_{n-1}(\Gamma^{\mathrm{3rd}}_{\delta,8}(z))\vert\big)^3\\
&\quad + (p-1)|\mathbb{X}_n(\Gamma_1)|^{2p-3}\int \dd z\, \eta(z)\, \big(|\mathbb{X}_n(\Gamma^{\mathrm{3rd}}_{\delta,9}(z))|+|\mathbb{X}_{n-1}(\Gamma^{\mathrm{3rd}}_{\delta,10}(z))|\big)\\
&\qquad\qquad\qquad\qquad\qquad\qquad\qquad\qquad\big(|\mathbb{X}_n(\Gamma^{\mathrm{3rd}}_{\delta,11}(z))|+|\mathbb{X}_{n-1}(\Gamma^{\mathrm{3rd}}_{\delta,12}(z))|+|\mathbb{X}_{n-2}(\Gamma^{\mathrm{3rd}}_{\delta,13}(z))|\big).
\end{aligned}
\end{equation}

\end{itemize}
\end{lemma}
\begin{proof}
For notational convenience, we write throughout the proof $\mathbb{X}_n$ for $\mathbb{X}_n(\Gamma_1)$. We split the proof into three steps, proving separately \eqref{EquationFirstOrderDeriv}, \eqref{EquationSecondOrderDeriv} and \eqref{EquationThirdOrderDeriv}.

\medskip

{\sc Step $1$. Proof of \eqref{EquationFirstOrderDeriv}. }First, we compute the Malliavin derivative $\mathrm{D}_{\delta,z_2}\mathbb{X}^{2p-1}_n$, using properties of the Wick product \eqref{PermutationWick} and \eqref{DiffEqWick}, that we express in terms of an $\mathbb{X}_{n-1}$-type quantity: on the one hand
\begin{equation}\label{InteMalliavinDerivTwoNoisesPowerLess}
\begin{aligned}
\mathrm{D}_{\delta,z_2}\mathbb{X}_n=&\mathbb{X}_n\big(\mathrm{D}_{\delta,z_2}\Gamma_1\big)+\sum_{\sigma\in \mathfrak{S}_n}\int \dd z'_1\cdots \dd z'_{n+1}\, \Gamma_1(z'_1,z'_{\sigma(2)},\cdots,z'_{\sigma(n+1)})\rho_\delta\star\rho_\delta(z'_2-z_2)\mathop{\diamond}_{i=3}^{n+1}\xi_\delta(z'_i)\\
=&\mathbb{X}_n\big(\mathrm{D}_{\delta,z_2}\Gamma_1\big)+\mathbb{X}_{n-1}(\Gamma_{\delta,z_2}),
\end{aligned}
\end{equation}
where we have set 
\begin{equation}\label{KernelXn}
\Gamma_{\delta,z_2}(z'_1,z'_3,\cdots ,z'_{n+1}):=\int \dd z'_2\,\sum_{\sigma\in \mathfrak{S}_n}\Gamma_1(z'_1,z'_{\sigma(2)},\cdots,z'_{\sigma(n+1)})\rho_\delta\star\rho_\delta(z'_2-z_2)\quad\text{for any $z'_1,z'_3,\cdots,z'_{n+1}\in \mathbb{R}^{1+2}$},
\end{equation}
and $\mathfrak{S}_n$ denotes the permutations of $\{2,\cdots,n+1\}$.  On the other hand, using \eqref{InteMalliavinDerivTwoNoisesPowerLess} we obtain
\begin{align}
\mathrm{D}_{\delta,z_2}\mathbb{X}_n^{2p-1}=(2p-1)\mathbb{X}_n^{2(p-1)}\big(\mathbb{X}_n\big(\mathrm{D}_{\delta,z_2}\Gamma_1\big)+\mathbb{X}_{n-1}(\Gamma_{\delta,z_2})\big).\label{InteMalliavinDerivTwoNoises}
\end{align}
Second, we set 
$$\Gamma^{\mathrm{1st}}_{\delta,1}:=\int \dd z_1\cdots\dd z_{n+1}\, \Gamma_2(z_1,\cdots,z_{n+1})\mathrm{D}_{\delta,z_2}\Gamma_1\quad\text{and}\quad \Gamma^{\mathrm{1st}}_{\delta,2}:=\int \dd z_1\cdots\dd z_{n+1}\, \Gamma_2(z_1,\cdots,z_{n+1})\Gamma_{\delta,z_2},$$
and \eqref{InteMalliavinDerivTwoNoises} leads to \eqref{EquationFirstOrderDeriv}. Finally, since $\Gamma^{\mathrm{1st}}_{\delta,1}\in \mathcal{K}_{n}(\eta\int \eta)$ follows directly from the assumption $\Gamma_1,\Gamma_2\in \mathcal{K}_{n}(\eta)$, it remains to show that $\Gamma^{\mathrm{1st}}_{\delta,2}\in \mathcal{K}_{n-1}(\eta_{1})$. Using Lemma \ref{LemmaConvoBound} and further the assumption $\Gamma_1,\Gamma_2\in \mathcal{K}_n(\eta)$, we have for any $\ell\geq 1$, $z''_1,\cdots,z''_{\ell}\in \mathbb{R}^{1+2}$ and $z'_1,z'_3,\cdots, z'_{n+1}\in \mathbb{R}^{1+2}$
\begin{align}
&\big\vert \Gamma^{\mathrm{1st}}_{\delta,2}(z'_1,z'_3,\dots,z'_{n+1})\big\vert +\big\vert \mathrm{D}_{\delta,z''_1,\cdots, z''_{\ell}} \Gamma^{\mathrm{1st}}_{\delta,2}(z'_1,z'_3,\cdots,z'_{n+1})\big\vert\nonumber\\
& \lesssim \eta(z'_1)\int \dd z_1\, \eta(z_1)\bigg(\int \dd z_2\, H_{C_\ell}(z_1-z_2)\int \dd z'_2\, H_{C_\ell}(z'_1-z'_2)\rho_\delta\star \rho_\delta(z'_2-z_2)\bigg)\prod_{i=3}^{n+1} H_{C_\ell}(z'_1-z'_i)\nonumber\\
&\stackrel{\eqref{BoundRandomTestFunction}}{\lesssim} \eta_{1}(z'_1)\prod_{i=3}^{n+1} H_{C_\ell}(z'_1-z'_i),\label{ShowInKappan1}
\end{align}
which means $\Gamma^{\mathrm{1st}}_{\delta,2}\in \mathcal{K}_{n-1}(\eta_{1})$.

\medskip

{\sc Step $2$. Proof of \eqref{EquationSecondOrderDeriv}. }First, we compute the second-order Malliavin derivative $\mathrm{D}_{\delta,z_2,z_3}\mathbb{X}^{2p-1}_n$ by using \eqref{InteMalliavinDerivTwoNoisesPowerLess}: on the one hand 
\begin{equation}\label{Malliavin2ndOrderTwoNoisesPowerLess}
\begin{aligned}
\mathrm{D}_{\delta,z_2,z_3}\mathbb{X}_n=&\mathrm{D}_{\delta,z_3}\mathbb{X}_n(\mathrm{D}_{\delta,z_2}\Gamma_1)+\mathrm{D}_{\delta,z_3}\mathbb{X}_{n-1}(\Gamma_{\delta,z_2})\\
=&\mathrm{D}_{\delta,z_3}\mathbb{X}_n(\mathrm{D}_{\delta,z_2}\Gamma_1)+\mathbb{X}_{n-1}(\mathrm{D}_{\delta,z_3}\Gamma_{\delta,z_2})\\
&+\sum_{\sigma\in \mathfrak{S}_{n-1}}\int \dd z'_1\, \dd z'_3\cdots \dd z'_{n+1}\, \Gamma_{\delta,z_2}(z'_1,z'_{\sigma(3)},\cdots,z'_{\sigma(n+1)})\rho_\delta\star\rho_\delta(z'_3-z_3)\mathop{\diamond}_{i=4}^{n+1} \xi_\delta(z'_i)\\
=&\mathrm{D}_{\delta,z_3}\mathbb{X}_n(\mathrm{D}_{\delta,z_2}\Gamma_1)+\mathbb{X}_{n-1}(\mathrm{D}_{\delta,z_3}\Gamma_{\delta,z_2})+\mathbb{X}_{n-2}(\Gamma_{\delta,z_2,z_3}),
\end{aligned}
\end{equation}
where we have set 
\begin{equation}\label{DefGammaTwoIndices}
\Gamma_{\delta,z_2,z_3}(z'_1,z'_4,\cdots ,z'_{n+1}):=\int \dd z'_3\,\sum_{\sigma\in \mathfrak{S}_{n-1}}\Gamma_{\delta,z_2}(z'_1,z'_{\sigma(3)},\cdots,z'_{\sigma(n+1)})\rho_\delta\star\rho_\delta(z'_3-z_3)\quad\text{for any $z'_1,z'_3,\cdots,z'_{n+1}\in \mathbb{R}^{1+2}$},
\end{equation}
and $\mathfrak{S}_{n-1}$ denotes the permutations of $\{3,\cdots,n+1\}$. On the other hand, using \eqref{Malliavin2ndOrderTwoNoisesPowerLess} we obtain 
\begin{equation}\label{Malliavin2ndOrderTwoNoisesP}
\begin{aligned}
\mathrm{D}_{\delta,z_2,z_3}\mathbb{X}^{2p-1}_n=&(2p-1)\mathbb{X}^{2(p-1)}_n\big(\mathrm{D}_{\delta,z_3}\mathbb{X}_n(\mathrm{D}_{\delta,z_2}\Gamma_1)+\mathbb{X}_{n-1}(\mathrm{D}_{\delta,z_3}\Gamma_{\delta,z_2})\big)\\
& +(2p-1)\mathbb{X}^{2(p-1)}_n\mathbb{X}_{n-2}(\Gamma_{\delta,z_2,z_3})\\
&+2(2p-1)(p-1)\mathbb{X}^{2p-3}_n\big(\mathbb{X}_n\big(\mathrm{D}_{\delta,z_2}\Gamma_1\big)+\mathbb{X}_{n-1}(\Gamma_{\delta,z_2})\big)\big(\mathbb{X}_n\big(\mathrm{D}_{\delta,z_3}\Gamma_1\big)+\mathbb{X}_{n-1}(\Gamma_{\delta,z_3})\big).
\end{aligned}
\end{equation}
We then treat each term in the r. h. s. of \eqref{Malliavin2ndOrderTwoNoisesP} separately:
\begin{itemize}
\item The first r. h. s. term of \eqref{Malliavin2ndOrderTwoNoisesP} together with \eqref{InteMalliavinDerivTwoNoisesPowerLess} (with $\Gamma_1$ replaced by $\mathrm{D}_{\delta,z_2}\Gamma_1$) and the assumption $\Gamma_1,\Gamma_2\in \mathcal{K}_{n}(\eta)$ provides a similar contribution as of \eqref{EquationFirstOrderDeriv}, which gives the first r. h. s. term of \eqref{EquationSecondOrderDeriv};
\item The second r. h. s. term of \eqref{Malliavin2ndOrderTwoNoisesP} can be treated the same way as the second r. h. s. term of \eqref{InteMalliavinDerivTwoNoises}, which gives the second r. h. s. term of \eqref{EquationSecondOrderDeriv};
\item Finally, we estimate the last r. h. s. term of \eqref{Malliavin2ndOrderTwoNoisesP}. To this end, we use the assumption $\Gamma_2\in \mathcal{K}_{n}(\eta)$ in the form of: there exists $C>0$ such that 
\begin{equation}\label{EstiProductTermThirdMDeriv}
\begin{aligned}
&\bigg\vert \int \dd z_1\cdots \dd z_{n+1}\, \Gamma_2(z_1,\cdots,z_{n+1})\big(\mathbb{X}_n\big(\mathrm{D}_{\delta,z_2}\Gamma_1\big)+\mathbb{X}_{n-1}(\Gamma_{\delta,z_2})\big)\big(\mathbb{X}_n\big(\mathrm{D}_{\delta,z_3}\Gamma_1\big)+\mathbb{X}_{n-1}(\Gamma_{\delta,z_3})\big)\bigg\vert\\
&\lesssim \int \dd z_1\, \dd z_2\, \dd z_3\,\eta(z_1)H_{C}(z_1-z_2)H_{C}(z_1-z_3)\big(\vert\mathbb{X}_n\big(\mathrm{D}_{\delta,z_2}\Gamma_1\big)\vert+\vert\mathbb{X}_{n-1}(\Gamma_{\delta,z_2})\vert\big)\big(\vert\mathbb{X}_n\big(\mathrm{D}_{\delta,z_3}\Gamma_1\big)\vert+\vert\mathbb{X}_{n-1}(\Gamma_{\delta,z_3})\vert\big)\\
&=\int \dd z_1\, \eta(z_1)\bigg(\Big\vert\mathbb{X}_n\Big(\underbrace{\int \dd z_2\, H_C(z_1-z_2)\mathrm{D}_{\delta,z_2}\Gamma_1}_{:=\Gamma^{\mathrm{2nd}}_{\delta,4}(z_1)}\Big)\Big\vert+\Big\vert\mathbb{X}_{n-1}\Big(\underbrace{\int \dd z_2\, H_C(z_1-z_2)\Gamma_{\delta,z_2}}_{:=\Gamma^{\mathrm{2nd}}_{\delta,5}(z_1)}\Big)\Big\vert\bigg)^2.
\end{aligned}
\end{equation}
Finally, $\Gamma^{\mathrm{2nd}}_{\delta,4}(z_1)\in \mathcal{K}_n(\eta)$ from the assumption $\Gamma_1\in \mathcal{K}_n(\eta)$ and $\Gamma^{\mathrm{2nd}}_{\delta,5}(z_1)\in \mathcal{K}_{n-1}(\eta\,\mathcal{C}_{\delta,C}(\cdot,z_1))$ similarly to \eqref{ShowInKappan1} using that $\Gamma_1,\Gamma_2\in \mathcal{K}_n(\eta)$ and Lemma \ref{LemmaConvoBound}.
\end{itemize}

\medskip

{\sc Step $3$. Proof of \eqref{EquationThirdOrderDeriv}. }First, we compute the third-order Malliavin derivative $\mathrm{D}_{\delta,z_2,z_3,z_4}\mathbb{X}^{2p-1}_n$ by using \eqref{InteMalliavinDerivTwoNoisesPowerLess}: on the one hand 
\begin{equation}\label{Malliavin3rdOrderTwoNoisesPowerLess}
\begin{aligned}
\mathrm{D}_{\delta,z_2,z_3,z_4}\mathbb{X}_n=&\mathrm{D}_{\delta,z_3,z_4}\mathbb{X}_n(\mathrm{D}_{\delta,z_2}\Gamma_1)+\mathrm{D}_{\delta,z_4}\mathbb{X}_{n-1}(\mathrm{D}_{\delta,z_3}\Gamma_{\delta,z_2})+\mathbb{X}_{n-2}(\mathrm{D}_{\delta,z_4}\Gamma_{\delta,z_2,z_3})\\
&+\sum_{\sigma\in \mathfrak{S}_{n-2}}\int \dd z'_1\, \dd z'_4\cdots \dd z'_{n+1}\, \Gamma_{\delta,z_2,z_3}(z'_1,z'_{\sigma(4)},\cdots,z'_{\sigma(n+1)})\rho_\delta\star\rho_\delta(z'_4-z_4)\mathop{\diamond}_{i=5}^{n+1}\xi_\delta(z'_i)\\
=&\mathrm{D}_{\delta,z_3,z_4}\mathbb{X}_n(\mathrm{D}_{\delta,z_2}\Gamma_1)+\mathrm{D}_{\delta,z_4}\mathbb{X}_{n-1}(\mathrm{D}_{\delta,z_3}\Gamma_{\delta,z_2})+\mathbb{X}_{n-2}(\mathrm{D}_{\delta,z_4}\Gamma_{\delta,z_2,z_3})+\mathbb{X}_{n-3}(\Gamma_{\delta,z_2,z_3,z_4})
\end{aligned}
\end{equation}
where we have set 
$$\Gamma_{\delta,z_2,z_3,z_4}:=\int \dd z'_4\, \sum_{\sigma\in \mathfrak{S}_{n-2}}\int \dd z'_1\, \dd z'_4\cdots \dd z'_{n+1}\, \Gamma_{\delta,z_2,z_3}(z'_1,z'_{\sigma(4)},\cdots,z'_{\sigma(n+1)})\rho_\delta\star\rho_\delta(z'_4-z_4),$$
and $\mathfrak{S}_{n-2}$ denotes the permutations of $\{4,\cdots,n+1\}$. On the other hand, expanding 
$$\mathrm{D}_{\delta,z_4}\mathbb{X}_{n-1}(\mathrm{D}_{\delta,z_3}\Gamma_{\delta,z_2})=\mathbb{X}_{n-1}(\mathrm{D}_{\delta,z_3,z_4}\Gamma_{\delta,z_2})+\mathrm{X}_{n-2}(\mathrm{D}_{\delta,z_3}\Gamma_{\delta,z_2,z_4}),$$
where we recall that $\Gamma_{\delta,z_2,z_4}$ is defined in \eqref{DefGammaTwoIndices}, and using \eqref{InteMalliavinDerivTwoNoisesPowerLess}, \eqref{Malliavin2ndOrderTwoNoisesPowerLess} and \eqref{Malliavin3rdOrderTwoNoisesPowerLess}, we obtain
\begin{equation}\label{Malliavin3rdOrderTwoNoises}
\begin{aligned}
&\mathrm{D}_{\delta,z_2,z_3,z_4}\mathbb{X}^{2p-1}_n\\
&=(2p-1)\mathbb{X}^{2(p-1)}_n\mathrm{D}_{\delta,z_2,z_3,z_4}\mathbb{X}_n+ 2(2p-1)(p-1)\mathbb{X}^{2p-4}_n\mathrm{D}_{\delta,z_2}\mathbb{X}_n\mathrm{D}_{\delta,z_3}\mathbb{X}_n\mathrm{D}_{\delta,z_4}\mathbb{X}_n\\
&\quad +2(2p-1)(p-1)(2p-3)\mathbb{X}^{2p-3}_n\sum_{\sigma\in \mathfrak{S}(\{2,3,4\})}\mathrm{D}_{\delta,z_{\sigma(2)}}\mathbb{X}_n\mathrm{D}_{\delta,z_{\sigma(3)},z_{\sigma(4)}}\mathbb{X}_n\\
&=(2p-1)\mathbb{X}^{2(p-1)}_{n}\Big(\mathrm{D}_{\delta,z_3,z_4}\mathbb{X}_n(\mathrm{D}_{\delta,z_2}\Gamma_1)+\mathbb{X}_{n-1}(\mathrm{D}_{\delta,z_3,z_4}\Gamma_1)+\mathbb{X}_{n-2}(\mathrm{D}_{\delta,z_3}\Gamma_{\delta,z_2,z_4}+\mathrm{D}_{\delta,z_4}\Gamma_{\delta,z_2,z_3})\Big)\\
&\quad + (2p-1)\mathbb{X}^{2(p-1)}_{n}\mathbb{X}_{n-3}(\Gamma_{\delta,z_2,z_3,z_4})\\
&\quad+2(2p-1)(p-1)(2p-3)\mathbb{X}^{2p-4}_n\big(\mathbb{X}_{n}(\mathrm{D}_{\delta,z_2}\Gamma_1)+\mathbb{X}_{n-1}(\Gamma_{\delta,z_2})\big)\big(\mathbb{X}_{n}(\mathrm{D}_{\delta,z_3}\Gamma_1)+\mathbb{X}_{n-1}(\Gamma_{\delta,z_3})\big)\big(\mathbb{X}_{n}(\mathrm{D}_{\delta,z_4}\Gamma_1)+\mathbb{X}_{n-1}(\Gamma_{\delta,z_4})\big)\\
&\quad+2(2p-1)(p-1)\mathbb{X}^{2p-3}_n\sum_{\sigma\in \mathfrak{S}(\{2,3,4\})}\\
&\quad\quad \Big(\mathbb{X}_n(\mathrm{D}_{\delta,z_{\sigma(2)}}\Gamma_1)+\mathbb{X}_{n-1}(\Gamma_{\delta,z_{\sigma(2)}})\Big)\Big(\mathbb{X}_{n}(\mathrm{D}_{\delta,z_{\sigma(3)},z_{\sigma(4)}}\Gamma_1)+\mathbb{X}_{n-1}(\mathrm{D}_{\delta,z_{\sigma(3)}}\Gamma_{\delta,z_{\sigma(4)}}+\mathrm{D}_{\delta,z_{\sigma(4)}}\Gamma_{\delta,z_{\sigma(3)}})+\mathbb{X}_{n-2}(\Gamma_{\delta,z_{\sigma(3)},z_{\sigma(4)}})\Big).
\end{aligned}
\end{equation}
We then treat each r. h. s. term of \eqref{Malliavin3rdOrderTwoNoises} separately:
\begin{itemize}
\item The first r. h. s. term together with \eqref{Malliavin2ndOrderTwoNoisesPowerLess} (with $\Gamma_1$ replaced by $\mathrm{D}_{\delta,z_2}\Gamma_1$) and the assumption $\Gamma_1,\Gamma_2\in \mathcal{K}_n(\eta)$ provides a similar contribution as of \eqref{EquationSecondOrderDeriv}, which gives the first and second r. h. s. term of \eqref{EquationThirdOrderDeriv};
\item The second r. h. s. term can be treated the same way as the second r. h. s. term of \eqref{InteMalliavinDerivTwoNoises}, which gives the third r. h. s. of \eqref{EquationThirdOrderDeriv};
\item The third r. h. s. term is treaty similarly as in \eqref{EstiProductTermThirdMDeriv}, which gives the forth r. h. s. term of \eqref{EquationThirdOrderDeriv};
\item Finally, we estimate the the last r. h. s. term of \eqref{Malliavin3rdOrderTwoNoises}. To this end, we use the assumption $\Gamma_2\in \mathcal{K}_{n}(\eta)$ in the form of: there exists $C>0$ such that, for any $\sigma\in \mathfrak{S}(\{2,3,4\})$
\begin{align*}
&\bigg\vert \int \dd z_1\cdots \dd z_{n+1}\, \Gamma_2(z_1,\cdots,z_{n+1})\\
&\Big(\mathbb{X}_n(\mathrm{D}_{\delta,z_{\sigma(2)}}\Gamma_1)+\mathbb{X}_{n-1}(\Gamma_{\delta,z_{\sigma(2)}})\Big)\\
&\Big(\mathbb{X}_{n}(\mathrm{D}_{\delta,z_{\sigma(3)},z_{\sigma(4)}}\Gamma_1)+\mathbb{X}_{n-1}(\mathrm{D}_{\delta,z_{\sigma(3)}}\Gamma_{\delta,z_{\sigma(4)}}+\mathrm{D}_{\delta,z_{\sigma(4)}}\Gamma_{\delta,z_{\sigma(3)}})+\mathbb{X}_{n-2}(\Gamma_{\delta,z_{\sigma(3)},z_{\sigma(4)}})\Big)\bigg\vert\\
&\lesssim \int \dd z_1\, \dd z_{\sigma(2)}\, \dd z_{\sigma(3)}\, \dd z_{\sigma(4)}\, \eta(z_1)\, H_C(z_1-z_{\sigma(2)})H_C(z_1-z_{\sigma(3)})H_C(z_1-z_{\sigma(4)})\Big(\mathbb{X}_n(\mathrm{D}_{\delta,z_{\sigma(2)}}\Gamma_1)+\mathbb{X}_{n-1}(\Gamma_{\delta,z_{\sigma(2)}})\Big)\\
&\quad \Big(\mathbb{X}_{n}(\mathrm{D}_{\delta,z_{\sigma(3)},z_{\sigma(4)}}\Gamma_1)+\mathbb{X}_{n-1}(\mathrm{D}_{\delta,z_{\sigma(3)}}\Gamma_{\delta,z_{\sigma(4)}}+\mathrm{D}_{\delta,z_{\sigma(4)}}\Gamma_{\delta,z_{\sigma(3)}})+\mathbb{X}_{n-2}(\Gamma_{\delta,z_{\sigma(3)},z_{\sigma(4)}})\Big)\\
&=\int \dd z_1\,\eta(z_1)\, \bigg(\mathbb{X}_n\bigg( \underbrace{\int \dd z_{\sigma(2)} H_C(z_1-z_{\sigma(2)})\mathrm{D}_{\delta,z_{\sigma(2)}}\Gamma_1}_{:=\Gamma^{\mathrm{3rd}}_{\delta,9}}\bigg)+\mathbb{X}_{n-1}\bigg(\underbrace{\int \dd z_{\sigma(2)} H_C(z_1-z_{\sigma(2)})\Gamma_{\delta,z_{\sigma(2)}}}_{:=\Gamma^{\mathrm{3rd}}_{\delta,10}}\bigg)\bigg)\\
&\quad \bigg(\mathbb{X}_{n}\bigg(\underbrace{\int \dd z_{\sigma(3)}\, \dd z_{\sigma(4)}\, H_C(z_1-z_{\sigma(3)})H_C(z_1-z_{\sigma(4)})\mathrm{D}_{\delta,z_{\sigma(3)},z_{\sigma(4)}}\Gamma_1}_{:=\Gamma^{\mathrm{3rd}}_{\delta,11}}\bigg)\\
&\quad +\mathbb{X}_{n-1}\bigg(\underbrace{\int \dd z_{\sigma(3)}\, \dd z_{\sigma(4)}\, H_C(z_1-z_{\sigma(3)})H_C(z_1-z_{\sigma(4)})\big(\mathrm{D}_{\delta,z_{\sigma(3)}}\Gamma_{\delta,z_{\sigma(4)}}+\mathrm{D}_{\delta,z_{\sigma(4)}}\Gamma_{\delta,z_{\sigma(3)}}\big)}_{:=\Gamma^{\mathrm{3rd}}_{\delta,12}}\bigg)\\
&\quad +\mathbb{X}_{n-2}\bigg(\underbrace{\int \dd z_{\sigma(3)}\, \dd z_{\sigma(4)}\, H_C(z_1-z_{\sigma(3)})H_C(z_1-z_{\sigma(4)})\Gamma_{\delta,z_{\sigma(3)},z_{\sigma(4)}}}_{:=\Gamma^{\mathrm{3rd}}_{\delta,13}}\bigg)\bigg).\end{align*}
\end{itemize}
Similarly as previously, we derive that $\Gamma^{\mathrm{3rd}}_{\delta,9},\Gamma^{\mathrm{3rd}}_{\delta,11}\in \mathcal{K}_n(\eta)$, $z_1\mapsto (\Gamma^{\mathrm{3rd}}_{\delta,10}(z_1),\Gamma^{\mathrm{3rd}}_{\delta,12}(z_1),\Gamma^{\mathrm{3rd}}_{\delta,13}(z_1))\in \mathcal{K}_{n-1}(\eta\,\mathcal{C}_\delta(\cdot,z_1))\times \mathcal{K}_{n-1}(\eta\, \mathcal{C}_\delta(\cdot,z_1))\times \mathcal{K}_{n-2}(\eta \,\mathcal{C}^2_\delta(\cdot,z_1))$.

\end{proof}
\subsection{Estitmate of the first-order quantity $\mathbb{X}_1$}
We show in this section the stochastic estimates of the first-order quantity $\mathbb{X}_1$. 
\begin{proposition}[Estimate of the $1^{\mathrm{st}}$-order quantity $\mathbb{X}_1$]\label{EstimateForOneNoise}
Let $\eta : \mathbb{R}^{1+2}\rightarrow (0,\infty)$ be measurable and deterministic such that $\mathrm{supp}\,\eta\subset [0,T)\times \mathbb{R}^2$ for some $T>0$. We consider a kernel $\Gamma\in \mathcal{K}_1(\eta)$ as defined in Definition \ref{DefClassKernel}. Then, the quantity
\begin{equation}\label{FirstOrderQuantityDef}
\mathbb{X}_1(\Gamma):=\int \dd z_1\, \dd z_2\,\Gamma(z_1,z_2)\xi_\delta(z_2),
\end{equation}
satisfies for any integer $p\geq 1$,
\begin{equation}\label{StochasticEstimateOneNoise}
\mathbb{E}^{\frac{1}{2p}}\Big[\big(\mathbb{X}_1(\Gamma)\big)^{2p}\Big]\lesssim_{T,p} \int \dd z\, \eta(z)+\bigg(\int \dd z\,\dd z'\, \eta(z)\eta(z')\big(1+\vert\log(\vert x-x'\vert)\vert\big)\bigg)^{\frac{1}{2}}.
\end{equation}
\end{proposition}
\begin{proof}
We show \eqref{StochasticEstimateOneNoise} by induction over $p\geq 1$, where the idea is to include all possible random kernels $\Gamma\in \mathcal{K}_1(\eta)$ in the induction hypothesis. For notational convenience we write, throughout the proof, $\mathbb{X}_1$ for $\mathbb{X}_1(\Gamma)$. Note that the assumptions on $\Gamma$ ensure that $\mathbb{X}_1$ satisfies the assumptions of Lemma \ref{IntWRTWick}. In the following, we use the notations of Definition \ref{QuantityNorm}.

\medskip

The strategy is to find an inequality estimating $\mathbb{X}^{2p}_1$ in terms of $\mathbb{X}^{2(p-1)}_1$ and then run the induction over $p\geq 1$. Our starting point is the following identity: for any $p\geq 1$
\begin{align*}
\mathbb{E}\big[\mathbb{X}_1^{2p}\big]=\int \dd z_1\,\dd z_2\,\mathbb{E}\Big[\mathbb{X}_1^{2p-1}\Gamma(z_1,z_2)\xi_\delta(z_2)\Big].
\end{align*}
We then use the Gaussian integration by parts \eqref{GaussianInteNnoises} with $n=1$ in the form of
\begin{equation}\label{FirstStepInductionOneNoise}
\begin{aligned}
\mathbb{E}\big[\mathbb{X}_1^{2p}\big]=\underbrace{\mathbb{E}\bigg[\mathbb{X}_1^{2p-1}\int \dd z_1\, \dd z_2\,\mathrm{D}_{\delta,z_2}\big(\Gamma(z_1,z_2)\big)\bigg]}_{:=\mathrm{I}_1}+\underbrace{\mathbb{E}\Big[\int \dd z_1\, \dd z_2\, \Gamma(z_1,z_2)\mathrm{D}_{\delta,z_2}\mathbb{X}^{2p-1}_1\Big]}_{:=\mathrm{I}_{2}}.
\end{aligned}
\end{equation}
We treat $\mathrm{I}_1$ using Young's inequality with exponents $(\frac{2p}{2p-1},2p)$, which gives applying the assumption \eqref{AssumptionKernelFunction1Noise} i. e.  $\Gamma\in \mathcal{K}_1(\eta)$: there exists $C_{T,p}>0$ such that
\begin{align*}
\mathrm{I}_1\leq \tfrac{1}{2} \mathbb{E}\big[\mathbb{X}_1^{2p}\big]+C_{T,p} \Big(\int \dd z\, \eta(z)\Big)^{2p},
\end{align*}
where we absorb the first term into the left-hand side of \eqref{FirstStepInductionOneNoise}, yielding to 
\begin{equation}\label{MomentBoundX1Reduced}
\mathbb{E}\big[\mathbb{X}_1^{2p}\big]\lesssim \Big(\int \dd z\, \eta(z)\Big)^{2p}+\mathrm{I}_{2}.
\end{equation}
We now turn to estimating $\mathrm{I_2}$.  To this end, we use \eqref{EquationFirstOrderDeriv} in Lemma \ref{LinkMalliavinDerive} in the form of: there exists $\Gamma^{\mathrm{1st}}_{\delta,1}\in \mathcal{K}_1(\eta\int \eta)$ and $\Gamma^{\mathrm{1st}}_{\delta,2}\in \mathcal{K}_0(\eta_1)$ with 
$$\eta_1(z'):=\eta(z')\int \dd z\, \eta(z)\, \big(1+|\log(|x'-x|)|\big)$$
such that 
\begin{equation}\label{MalliavinDeriveOneNoise}
\mathrm{I}_2=(2p-1)\mathbb{E}\Big[\mathbb{X}^{2(p-1)}_1\big(\mathbb{X}_1(\Gamma^{\mathrm{1st}}_{\delta,1})+\mathbb{X}_0(\Gamma^{\mathrm{1st}}_{\delta,2})\big)\Big],
\end{equation}
where from $\Gamma^{\mathrm{1st}}_{\delta,2}\in \mathcal{K}_0(\eta_1)$ 
$$|\mathbb{X}_0(\Gamma^{\mathrm{1st}}_{\delta,2})|\lesssim \int \dd z\, \eta_1(z)=\int \dd z\, \dd z'\, \eta(z)\eta(z')\big(1+\vert\log(\vert z-z'\vert)\vert\big).$$
The combination of \eqref{MomentBoundX1Reduced} and \eqref{MalliavinDeriveOneNoise} leads to
\begin{align*}
\mathbb{E}\big[\mathbb{X}_1^{2p}\big]\lesssim &\Big(\int \dd z\, \eta(z)\Big)^{2p}+\mathbb{E}\Big[\mathbb{X}^{2(p-1)}_1\Big]\int \dd z_1\, \dd z'_1\, \eta(z_1)\eta(z'_1)\big(1+\vert\log(\vert z_1-z'_1\vert)\vert\big)+\mathbb{E}\Big[\mathbb{X}_1^{2(p-1)}\mathbb{X}_1(\Gamma^{\mathrm{1st}}_{\delta,1})\Big].
\end{align*}
Finally, for $p\geq 2$, we use Young's inequality with exponents $(\frac{p}{p-1},p)$ to absorb $\mathbb{E}\big[\mathbb{X}^{2(p-1)}_1\big]$ into the l. h. s., leading to
\begin{equation}\label{EstimateForInducX1}
\mathbb{E}\big[\mathbb{X}_1^{2p}\big]\lesssim\Big(\int \dd z\, \eta(z)\Big)^{2p}+\Big(\int \dd z_1\, \dd z'_1\, \eta(z_1)\eta(z'_1)\big(1+\vert\log(\vert z_1-z'_1\vert)\vert\big)\Big)^p+\mathbb{E}\Big[\mathbb{X}_1^{2(p-1)}\mathbb{X}_1(\Gamma^{\mathrm{1st}}_{\delta,1})\Big],
\end{equation}
which still holds true for $p=1$. We now argue by induction over the moment's exponent $p\geq 1$: 
\begin{itemize}
\item[(i)]For the initialization, recall that $\Gamma^{\mathrm{1st}}_{\delta,1}\in \mathcal{K}_1(\eta\int \eta)$ and using the integration by parts formula in Lemma \ref{IntWRTWick}, it holds
$$\mathbb{E}\big[\mathbb{X}_1(\Gamma^{\mathrm{1st}}_{\delta,1})\big]=\int \dd z_1\, \dd z_2\, \mathbb{E}\big[\mathrm{D}_{\delta,z_2}\Gamma^{\mathrm{1st}}_{\delta,1}(z_1,z_2)\big]\lesssim \Big(\int \dd z\, \eta(z)\Big)^2,$$
which, with \eqref{EstimateForInducX1}, gives \eqref{StochasticEstimateOneNoise} for $p=1$.

\medskip

\item[(ii)]For the induction step, we simply use Young's inequality with exponents $(\frac{p}{p-1},p)$ to the third r. h. s. term of \eqref{EstimateForInducX1} in order to absorb $\mathbb{E}\big[\mathbb{X}^{2(p-1)}_1\big]$ into the l. h. s. in the form of
$$\mathbb{E}\big[\mathbb{X}_1^{2p}\big]\lesssim\Big(\int \dd z\, \eta(z)\Big)^{2p}+\Big(\int \dd z_1\, \dd z'_1\, \eta(z_1)\eta(z'_1)\big(1+\vert\log(\vert z_1-z'_1\vert)\vert\big)\Big)^p+\mathbb{E}\big[|\mathbb{X}_1(\Gamma^{\mathrm{1st}}_{\delta,1})|^p\big],$$
combined with the induction hypothesis to $\mathbb{E}\Big[|\mathbb{X}_1(\Gamma^{\mathrm{1st}}_{\delta,1})|^p\Big]$ using that $(\int\eta)^{-1}\Gamma^{\mathrm{1st}}_{\delta,1}\in \mathcal{K}_1(\eta)$.
\end{itemize}
\end{proof}

%
%
%
\subsection{Estimate of the second-order quantity $\mathbb{X}_2$}
We show in this section the stochastic estimates of the second-order quantity $\mathbb{X}_2$. 
%
%
%
%
\begin{proposition}[Estimate of the $2^{\mathrm{nd}}$-order quantity $\mathbb{X}_2$]\label{X2Estimates}
Let $\psi^\lambda:=\lambda^{-4}\psi(\frac{\cdot}{\lambda^2},\frac{\cdot}{\lambda})$ with $\text{supp } \psi\subset \cc_1$ and $\|\psi\|_{\cc^0(\mathbb{R}^{1+2})}\leq 1$. For any $\lambda\in (0,e^{-1}]$ and $|\tilde{z}|\leq \lambda$, we define 
\begin{equation}\label{AssumptionTestFunction2Noise}
\begin{aligned}
&\eta^\lambda_{\tilde{z}}(z):= \vert \psi^\lambda(z)\vert \big(1+\vert \log(\vert x-\tilde{x}\vert)\vert\big)\quad\text{for any $z\neq \tilde{z}\in \mathbb{R}^{1+2}$}.
\end{aligned}
\end{equation}
We finally consider a kernel $\Gamma^\lambda_{\tilde{z}}\in \mathcal{K}_{2}(\eta^\lambda_{\tilde{z}})$ symmetric w. r. t. its two lasts variables. Then the quantity
\begin{equation}\label{DefX2Quant}
\mathbb{X}_2(\Gamma^\lambda_{\tilde{z}}):=\int \dd z_1\,\dd z_2\, \dd z_3\, \Gamma^\lambda_{\tilde{z}}(z_1,z_2,z_3)\xi_\delta(z_2)\diamond\xi_\delta(z_3),
\end{equation}
satisfies for any integer $p\geq 1$ and $\lambda\in (0,e^{-1}]$,
\begin{equation}\label{RefinementMomentBound}
\mathbb{E}^\frac{1}{2p}\Big[\big(\mathbb{X}_2(\Gamma^\lambda_{\tilde{z}})\big)^{2p}\Big]\lesssim \mathrm{C}(\lambda,\tilde{x}),
\end{equation}
where
$$\mathrm{C}(\lambda,\tilde{x}):=|\log(\lambda)|\bigg(1+\fint_{\bb_\lambda}\dd x\,|\log(|x-\tilde{x}|)|^2\bigg)^\frac{1}{2}.$$
In particular, 
\begin{equation}\label{RefinementMomentBoundInte}
\mathbb{E}^\frac{1}{2p}\Big[\big(\mathbb{X}_2(\Gamma^\lambda)\big)^{2p}\Big]\lesssim |\log(\lambda)|^2\quad \text{for $\Gamma^\lambda:=\fint_{\bb_{\lambda}}\dd \tilde{z}\, \Gamma^\lambda_{\tilde{z}}$}.
\end{equation}
Furthermore, if the term $1+|\log(|x-\tilde{x}|)|$ is dropped in \eqref{AssumptionTestFunction2Noise} then \eqref{RefinementMomentBoundInte} holds with $|\log(\lambda)|$ on the r. h. s. 
\end{proposition}
\begin{proof}
We fix $|\tilde{z}|\leq \lambda$ and drop the dependence in $\tilde{z}$ in the notations for convenience. We show \eqref{RefinementMomentBound} by induction over $p\geq 1$, where the idea is to include all possible random kernels $\Gamma^\lambda\in \mathcal{K}_2(\eta^\lambda)$ in the induction hypothesis. For notational convenience we write, throughout the proof, $\mathbb{X}_2$ for $\mathbb{X}_2(\Gamma^{\lambda})$.  We also use the notations introduced in Lemma \ref{LemmaConvoBound} and Definition \ref{QuantityNorm}. Note that the assumptions on $\Gamma^\lambda$ ensure that $\mathbb{X}_2$ satisfies the assumptions of Lemma \ref{IntWRTWick}.  

\medskip

The strategy is to find an inequality estimating $\mathbb{E}\big[\mathbb{X}^{2p}_2\big]$ in terms of $\mathbb{E}\big[\mathbb{X}^{2(p-1)}_2\big]$ and $\mathbb{E}\big[\mathbb{X}^{2p-3}_2\big]$ and then run the induction over $p\geq 1$. 

\medskip

As in the proof of Proposition \ref{EstimateForOneNoise}, our starting point is the following identity: for any $p\geq 1$
\begin{align*}
\mathbb{E}\big[\mathbb{X}_2^{2p}\big]=\int \dd z_1\,\dd z_2\,\dd z_3\, \mathbb{E}\big[\mathbb{X}_2^{2p-1}\Gamma^\lambda(z_1,z_2,z_3)\xi_\delta(z_2)\diamond \xi_\delta(z_3)\big].
\end{align*}
We then use the Gaussian integration by parts \eqref{GaussianInteNnoises} with $n=2$ in the form of
\begin{equation}\label{AfterGaussianInteTwoNoise}
\begin{aligned}
\mathbb{E}\big[\mathbb{X}_2^{2p}\big]=&\underbrace{\mathbb{E}\bigg[\mathbb{X}_2^{2p-1}\int\dd z_1\, \dd z_2\, \dd z_3\,\mathrm{D}_{\delta,z_2,z_3}\Gamma^\lambda(z_1,z_2,z_3)\bigg]}_{:=\mathrm{I}_1}+2\underbrace{\mathbb{E}\bigg[\int \dd z_2\,\mathrm{D}_{\delta,z_2}\mathbb{X}_2^{2p-1}\int \dd z_1\, \dd z_3\,\mathrm{D}_{\delta,z_3}\Gamma^\lambda(z_1,z_2,z_3)\bigg]}_{:=\mathrm{I}_2}\\
&+\underbrace{\mathbb{E}\bigg[\int \dd z_2\, \dd z_3\,\mathrm{D}_{\delta,z_2,z_3}\mathbb{X}_2^{2p-1}\int \dd z_1\, \Gamma^\lambda(z_1,z_2,z_3)\bigg].}_{:=\mathrm{I}_3}
\end{aligned}
\end{equation}
We treat $\mathrm{I}_1$ using Young's inequality with exponents $(\frac{2p}{2p-1},2p)$ and the assumption \eqref{AssumptionKernelFunction1Noise} i. e. $\Gamma^\lambda\in \mathcal{K}_2(\eta^\lambda)$ in the form of: there exists $C_p>0$ such that
\begin{equation}\label{CounterpartsI1X3}
\begin{aligned}
\mathrm{I}_1\lesssim \mathbb{E}\big[\mathbb{X}_2^{2p-1}\big]\int \dd z\, \eta^\lambda_{\tilde{z}}(z)&\leq\tfrac{1}{2} \mathbb{E}\big[\mathbb{X}_2^{2p}\big]+C_p\Big(1+\fint_{\bb_\lambda}\dd x\,\vert \log(\vert x-\tilde{x}\vert)\vert\Big)^{2p}\\
&\leq \tfrac{1}{2} \mathbb{E}\big[\mathbb{X}_2^{2p}\big]+C_p (\mathrm{C}(\lambda,\tilde{x}))^{2p},
\end{aligned}
\end{equation}
where we then absorb the first r. h. s. term into the l. h. s. of \eqref{AfterGaussianInteTwoNoise}, yielding to
\begin{equation}\label{MomentBoundX2Reduced}
\mathbb{E}\big[\mathbb{X}_2^{2p}\big]\lesssim (\mathrm{C}(\lambda,\tilde{x}))^{2p}+\mathrm{I}_2+\mathrm{I}_3.
\end{equation}
We now turn to the estimate of $\mathrm{I}_2$. To this end, we use \eqref{EquationFirstOrderDeriv} in Lemma \ref{LinkMalliavinDerive} in the form of: there exists $\Gamma^{\lambda,\mathrm{1st}}_{\delta,1}\in \mathcal{K}_2(\eta^\lambda\int \eta^\lambda)$ and $\Gamma^{\lambda,\mathrm{1st}}_{\delta,2}\in \mathcal{K}_1(\eta^\lambda_{1})$, where for any $i\geq 1$
$$\eta^\lambda_{i}(z'):=\eta^{\lambda}(z')\int \dd z\, \eta^\lambda(z)\big(1+\vert \log(\vert x'-x\vert)\vert\big)^i,$$
such that
\begin{equation}\label{I2FORX2WLEMMA}
\mathrm{I}_2=(2p-1)\mathbb{E}\Big[\mathbb{X}^{2(p-1)}_2\big(\mathbb{X}_2(\Gamma^{\lambda,\mathrm{1st}}_{\delta,1})+\mathbb{X}_1(\Gamma^{\lambda,\mathrm{1st}}_{\delta,2})\big)\Big].
\end{equation}
Furthermore, by Proposition \ref{EstimateForOneNoise} and Lemma \ref{GeneLog} with $n=4$, we deduce for any $p\geq 1$
\begin{equation}\label{BoundX1WithProp11}
\begin{aligned}
\mathbb{E}^{\frac{1}{p}}\Big[\big\vert\mathbb{X}_1(\Gamma^{\lambda,\mathrm{1st}}_{\delta,2})\big\vert^p\Big]\lesssim_p \int \dd z\, \eta^\lambda_{1}(z)+\bigg(\int \dd z\, \dd z'\, \eta^\lambda_{1}(z)\eta^\lambda_{1}(z')\big(1+\vert\log(\vert x-x'\vert)\vert\big)\bigg)^{\frac{1}{2}}\lesssim  (\mathrm{C}(\lambda,\tilde{x}))^2.
\end{aligned}
\end{equation}

The combination of \eqref{MomentBoundX2Reduced}, \eqref{I2FORX2WLEMMA} and \eqref{BoundX1WithProp11} together with Young's inequality with exponents $(\frac{p}{p-1},p)$ to absorb part of the second r. h. s. term of \eqref{I2FORX2WLEMMA} into the l. h. s. of \eqref{MomentBoundX2Reduced} yields 
\begin{equation}\label{MomentBoundX2Reduced2}
\begin{aligned}
\mathbb{E}\big[\mathbb{X}_2^{2p}\big]\lesssim & (\mathrm{C}(\lambda,\tilde{x}))^{2p}+\mathbb{E}\Big[\mathbb{X}^{2(p-1)}_2\mathbb{X}_2(\Gamma^{\lambda,\mathrm{1st}}_{\delta,1})\Big]+\mathrm{I}_3.
\end{aligned}
\end{equation}
We finally turn to the estimate of $\mathrm{I}_3$. We use this time \eqref{EquationSecondOrderDeriv} in Lemma \ref{LinkMalliavinDerive}: there exist $\Gamma^{\lambda,\mathrm{2nd}}_{\delta,1}\in \mathcal{K}_2(\eta^\lambda\int \eta^\lambda)$, $\Gamma^{\lambda,\mathrm{2nd}}_{\delta,2}\in \mathcal{K}_{1}(\eta^\lambda_{1})$, $\Gamma^{\lambda,\mathrm{2nd}}_{\delta,3}\in \mathcal{K}_{0}(\eta^\lambda_{2})$ as well as $\mathbb{R}^{1+2}\ni z\mapsto (\Gamma^{\lambda,\mathrm{2nd}}_{\delta,4}(z),\Gamma^{\lambda,\mathrm{2nd}}_{\delta,5}(z))\in \mathcal{K}_{2}(\eta^\lambda)\times \mathcal{K}_{1}(\eta^\lambda \,\mathcal{C}_\delta(\cdot,z))$ such that
\begin{equation}\label{I3FORX2WLEMMA}
\begin{aligned}
\mathrm{I}_3\lesssim &\mathbb{E}\Big[\mathbb{X}^{2(p-1)}_2\big(\mathbb{X}_2(\Gamma^{\lambda,\mathrm{2nd}}_{\delta,1})+\mathbb{X}_{1}(\Gamma^{\lambda,\mathrm{2nd}}_{\delta,2})\big)+\mathbb{X}^{2(p-1)}_2|\mathbb{X}_{0}(\Gamma^{\lambda,\mathrm{2nd}}_{\delta,3})|\Big]\\
&+(p-1)\mathbb{E}\bigg[\vert\mathbb{X}_2\vert^{2p-3}\int \dd z\, \eta^\lambda(z)\, \big(\vert \mathbb{X}_2(\Gamma^{\lambda,\mathrm{2nd}}_{\delta,4}(z))\vert +\vert \mathbb{X}_{1}(\Gamma^{\lambda,\mathrm{2nd}}_{\delta,5}(z))\vert\big)^2\bigg].
\end{aligned}
\end{equation}
The first r. h. s. term of \eqref{I3FORX2WLEMMA} is controlled the same way as $\mathrm{I}_2$ in \eqref{I2FORX2WLEMMA}. For the second r. h. s. term of \eqref{I3FORX2WLEMMA} we apply Lemma \ref{GeneLog} for $n=2$ to get
\begin{equation}\label{X0FORX2I3}
\begin{aligned}
\vert \mathbb{X}_0(\Gamma^{\lambda,\mathrm{2nd}}_{\delta,3})\vert \lesssim \int \dd z\, \eta^{\lambda}_{2}(z)= \int \dd z\, \dd z'\, \eta^\lambda(z)\eta^\lambda(z')\big(1+\vert \log(\vert x-x'\vert)\vert\big)^2 \lesssim  & (\mathrm{C}(\lambda,\tilde{x}))^2.
\end{aligned}
\end{equation}
%
%
%
%
%
%
%
%
Finally, in the regime $p\geq \frac{3}{2}$, we bound the contribution coming from the product of the $\mathbb{X}_1$-type quantities in the third r. h. s. term of \eqref{I3FORX2WLEMMA} using Young's inequality with exponents $(\frac{2p}{2p-3},\frac{2}{3}p)$, $\mathbb{R}^{1+2}\ni z\mapsto \Gamma^{\lambda,\mathrm{2nd}}_{\delta,5}(z)\in \mathcal{K}_{1}(\eta^\lambda \,\mathcal{C}_\delta(\cdot,z))$, Lemma \ref{LemmaConvoBound}, Proposition \ref{EstimateForOneNoise} and Lemma \ref{GeneLog}: for any $\varepsilon>0$, there exists $C_\varepsilon>0$ such that 
\begin{align}
\mathbb{E}\bigg[\vert \mathbb{X}_2\vert^{2p-3}\int\dd z\, \eta^\lambda(z)\vert \mathbb{X}_1(\Gamma^{\lambda,\mathrm{2nd}}_{\delta,5}(z))\vert^2\bigg]\leq &\varepsilon \mathbb{E}[\mathbb{X}^{2p}_2]+C_\varepsilon \bigg(\int \dd z\, \eta^\lambda(z)\,\mathbb{E}^\frac{3}{2p}\Big[\big\vert \mathbb{X}_1(\Gamma^{\lambda,\mathrm{2nd}}_{\delta,5}(z))\big\vert^{\frac{4}{3}p}\Big]\bigg)^{\frac{2}{3}p}\nonumber\\
\stackrel{\eqref{StochasticEstimateOneNoise}}{\leq} & \varepsilon \mathbb{E}\big[\mathbb{X}^{2p}_2\big]+C_\varepsilon (\mathrm{C}(\lambda,\tilde{x}))^{2p}.\label{EstiI3ProdX1ForX2}
\end{align}
To conclude, the combination of \eqref{MomentBoundX2Reduced2}, \eqref{I3FORX2WLEMMA}, \eqref{X0FORX2I3} and \eqref{EstiI3ProdX1ForX2} yields for any $p\in \{1\}\cup [\frac{3}{2},\infty)$
\begin{equation}\label{EstimateForInducX2}
\begin{aligned}
\mathbb{E}\big[\mathbb{X}^{2p}_2\big]\lesssim \mathbb{E}\Big[\mathbb{X}^{2(p-1)}_2\mathbb{X}_2(\Gamma^{\lambda,\mathrm{1st}}_{\delta,1}+\Gamma^{\lambda,\mathrm{2nd}}_{\delta,1})\Big]+(p-1)\,\mathbb{E}\bigg[\vert \mathbb{X}_2\vert^{2p-3}\int \dd z\, \eta^\lambda(z)\vert \mathbb{X}_2(\Gamma^{\lambda,\mathrm{2nd}}_{\delta,4}(z))\vert^2\bigg]+(\mathrm{C}(\lambda,\tilde{x}))^{2p}.
\end{aligned}
\end{equation}
We now argue by induction over integers $p\geq 1$ and split the proof into two steps.

\medskip

{\sc Step 1. Initialization. }We note that \eqref{EstimateForInducX2} requires to initialize with $p\in \{1,2\}$. 

\medskip

For $p=1$, recalling that $\Gamma^{\lambda,\mathrm{1st}}_{\delta,1},\Gamma^{\lambda,\mathrm{2nd}}_{\delta,1}\in \mathcal{K}_2(\eta^\lambda\int \eta^\lambda)$, we have using the integration by parts formula \eqref{IntWRTWick}:
\begin{align*}
\mathbb{E}[\mathbb{X}^2_2] \lesssim \mathbb{E}[\mathbb{X}_2(\Gamma^{\lambda,\mathrm{1st}}_{\delta,1}+\Gamma^{\lambda,\mathrm{2nd}}_{\delta,1})] +(\mathrm{C}(\lambda,\tilde{x}))^{2}=\int \dd z_1\, \dd z_2\, \dd z_3\, \mathbb{E}\Big[\mathrm{D}_{\delta,z_2,z_3}\big(\Gamma^{\lambda,\mathrm{1st}}_{\delta,1}+\Gamma^{\lambda,\mathrm{2nd}}_{\delta,1}\big)\Big]+(\mathrm{C}(\lambda,\tilde{x}))^{2} \lesssim (\mathrm{C}(\lambda,\tilde{x}))^{2}.
\end{align*}
For $p=2$, we first estimate third order moments by applying \eqref{EstimateForInducX2} with $p=\frac{3}{2}$: recalling that $\mathbb{R}^{1+2}\ni z\mapsto \Gamma^{\lambda,\mathrm{2nd}}_{\delta,4}(z)\in \mathcal{K}_{2}(\eta^\lambda)$, we obtain using the bound on the second moments
\begin{align*}
\mathbb{E}\big[\mathbb{X}^{3}_2\big] &\lesssim \mathbb{E}\Big[\mathbb{X}_2\,\mathbb{X}_2(\Gamma^{\lambda,\mathrm{1st}}_{\delta,1}+\Gamma^{\lambda,\mathrm{2nd}}_{\delta,1})\Big]+\mathbb{E}\bigg[\int \dd z\, \eta^\lambda(z)\,\vert \mathbb{X}_2(\Gamma^{\lambda,\mathrm{2nd}}_{\delta,4}(z))\vert^2\bigg]+(\mathrm{C}(\lambda,\tilde{x}))^{3}\lesssim (\mathrm{C}(\lambda,\tilde{x}))^{3}.
\end{align*}
Finally, we can treat the case $p=2$ using Young's inequality with exponents $(2,2)$ and $(4,\frac{4}{3})$ as well as the bounds on the second and third moments: there exists $C>0$ (which may change from line to line) such that
\begin{align*}
\mathbb{E}[\mathbb{X}^4_2] &\lesssim \mathbb{E}[\mathbb{X}^2_2\, \mathbb{X}_2(\Gamma^{\lambda,\mathrm{1st}}_{\delta,1}+\Gamma^{\lambda,\mathrm{2nd}}_{\delta,1})]+\mathbb{E}\bigg[\vert \mathbb{X}_2\vert \int\dd z\, \eta^\lambda(z)\,\vert  \mathbb{X}_2(\Gamma^{\lambda,\mathrm{2nd}}_{\delta,4}(z))\vert^2\bigg]\\
&\leq \tfrac{1}{2}\mathbb{E}[\mathbb{X}^4_2]+C\mathbb{E}\big[(\mathbb{X}_2(\Gamma^{\lambda,\mathrm{1st}}_{\delta,1}+\Gamma^{\lambda,\mathrm{2nd}}_{\delta,1}))^2\big]+C\bigg(\int \dd z\, \eta^\lambda(z)\, \mathbb{E}^{\frac{3}{4}}\Big[\vert \mathbb{X}_2(\Gamma^{\lambda,\mathrm{2nd}}_{\delta,4}(z))\vert^{\frac{8}{3}}\Big]\bigg)^{\frac{4}{3}}\\
&\leq \tfrac{1}{2}\mathbb{E}[\mathbb{X}^4_2]+C(\mathrm{C}(\lambda,\tilde{x}))^{4},
\end{align*}
which gives the bound on the fourth moment.

\medskip

{\sc Step 2. Induction step. }We now assume that $p\geq 3$ and we show the induction step. We use Young's inequality in \eqref{EstimateForInducX2} with exponents $(p,\frac{p}{p-1})$ and $(\frac{2p}{2p-3},\frac{2}{3}p)$ in the form of: there exists $C>0$ such that
\begin{equation}\label{InducStepForX2}
\begin{aligned}
\mathbb{E}[\mathbb{X}^{2p}_2]\leq \tfrac{1}{2}\mathbb{E}[\mathbb{X}^{2p}_2]+C\mathbb{E}\big[\big(\mathbb{X}_2(\Gamma^{\lambda,\mathrm{1st}}_{\delta,1}+\Gamma^{\lambda,\mathrm{2nd}}_{\delta,1})\big)^p\big]+C\bigg(\int \dd z\, \eta^\lambda(z)\, \mathbb{E}^{\frac{3}{2p}}\Big[\vert \mathbb{X}_2(\Gamma^{\lambda,\mathrm{2nd}}_{\delta,4}(z))\vert^{\frac{4}{3}p}\Big]\bigg)^{\frac{2}{3}p}.
\end{aligned}
\end{equation}
Noticing that $\frac{4}{3}p\leq 2(p-1)$, we then use the induction hypothesis 
$$\mathbb{E}\Big[\big(\mathbb{X}_2(\Gamma^{\lambda,\mathrm{1st}}_{\delta,1}+\Gamma^{\lambda,\mathrm{2nd}}_{\delta,1})\big)^p\Big]\lesssim (\mathrm{C}(\lambda,\tilde{x}))^{p}\quad\text{and}\quad \mathbb{E}^{\frac{3}{2p}}\Big[\vert \mathbb{X}_2(\Gamma^{\lambda,\mathrm{2nd}}_{\delta,4}(z))\vert^{\frac{4}{3}p}\Big]\lesssim (\mathrm{C}(\lambda,\tilde{x}))^{2},$$
and conclude that 
$$\mathbb{E}[\mathbb{X}^{2p}_2]\lesssim (\mathrm{C}(\lambda,\tilde{x}))^{2p}.$$
\end{proof}
\subsection{Estimate of the third-order quantity $\mathbb{X}_3$}
We show in this section the stochastic estimates of the third-order quantity $\mathbb{X}_3$.  
\begin{proposition}[Estimate of the $3^{\text{rd}}$-order quantity $\mathbb{X}_3$]\label{EstiX3}
Let $\psi^\lambda:=\lambda^{-4}\psi(\frac{\cdot}{\lambda^2},\frac{\cdot}{\lambda})$ with $\mathrm{supp }\,\psi\subset \cc_1$ and $\|\psi\|_{\cc^0(\mathbb{R}^{1+2})}\leq 1$ and $\Gamma^\lambda\in \mathcal{K}_3(\psi^\lambda)$ symmetric w. r. t. its three lasts variables. Then, the quantity
$$\mathbb{X}_3(\Gamma^\lambda):=\int\dd z_1\, \dd z_2\, \dd z_3\, \dd z_4\, \Gamma^\lambda(z_1,z_2,z_3,z_4)\xi_\delta(z_2)\diamond \xi_\delta(z_3)\diamond \xi_\delta(z_4),$$
satisfies for any $p\geq 1$ and $\lambda\in (0,e^{-1}]$
\begin{equation}\label{MomentBounds3Noises}
\mathbb{E}^{\frac{1}{2p}}\Big[\big(\mathbb{X}_3(\Gamma^\lambda)\big)^{2p}\Big]\lesssim  \vert \log(\lambda)\vert^{\frac{3}{2}}.
\end{equation}
\end{proposition}
\begin{proof}
We show \eqref{RefinementMomentBound} by induction over $p\geq 1$, where the idea is to include all possible random kernels $\Gamma^\lambda\in \mathcal{K}_3(\psi^\lambda)$ in the induction hypothesis. For notational convenience we write, throughout the proof, $\mathbb{X}_3$ for $\mathbb{X}_2(\Gamma^{\lambda})$.  We also use the notations introduced in Lemma \ref{LemmaConvoBound} and Definition \ref{QuantityNorm}. Note that the assumptions on $\Gamma^\lambda$ ensure that $\mathbb{X}_3$ satisfies the assumptions of Lemma \ref{IntWRTWick}. 

\medskip

The strategy is to find an inequality estimating $\mathbb{E}\big[\mathbb{X}^{2p}_2\big]$ in terms of $\mathbb{E}\big[\mathbb{X}^{2(p-1)}_2\big]$, $\mathbb{E}\big[\mathbb{X}^{2p-3}_2\big]$ and $\mathbb{E}\big[\mathbb{X}^{2p-4}_2\big]$ and then run the induction over $p\geq 1$. 

\medskip

As in Proposition \ref{EstimateForOneNoise} and Proposition \ref{X2Estimates}, our starting point is the following identity
\begin{align*}
&\mathbb{E}\big[\mathbb{X}_3^{2p}\big]=\mathbb{E}\Big[\int\dd z_1\, \dd z_2\, \dd z_3\, \dd z_4\,\mathbb{X}^{2p-1}_3\Gamma^\lambda(z_1,z_2,z_3)\xi_\delta(z_2)\diamond \xi_\delta(z_3)\diamond \xi_\delta(z_4)\Big].
\end{align*}
We then use the Gaussian integration by parts \eqref{GaussianInteNnoises} with $n=3$ in the form of
\begin{equation}\label{Buckling3Noises}
\begin{aligned}
\mathbb{E}\big[\mathbb{X}_3^{2p}\big]=&\underbrace{\mathbb{E}\bigg[\mathbb{X}_3^{2p-1}\int \dd z_1\, \dd z_2\, \dd z_3\, \dd z_4\, \mathrm{D}_{\delta, z_2,z_3,z_4}\Gamma^\lambda (z_1,z_2,z_3,z_4)\bigg]}_{:=\mathrm{I}_1}\\
&+3\underbrace{\mathbb{E}\bigg[\int \dd z_2\, \mathrm{D}_{\delta,z_2}\mathbb{X}^{2p-1}_3\int \dd z_1\, \dd z_3\, \dd z_4\, \mathrm{D}_{\delta,z_3,z_4} \Gamma^\lambda(z_1,z_2,z_3,z_4)\bigg]}_{:=\mathrm{I}_2}\\
&+3\underbrace{\mathbb{E}\bigg[\int \dd z_2\, \dd z_3\,\mathrm{D}_{\delta,z_2,z_3} \mathbb{X}^{2p-1}_{3}\int \dd z_1\, \dd z_4\, \mathrm{D}_{\delta,z_4}\Gamma^\lambda(z_1,z_2,z_3,z_4)\bigg]}_{:=\mathrm{I}_3}\\
&+\underbrace{\mathbb{E}\bigg[\int \dd z_2\, \dd z_3\,\dd z_4\, \mathrm{D}_{\delta,z_2,z_3,z_4} \mathbb{X}^{2p-1}_{3}\int \dd z_1\, \Gamma^\lambda(z_1,z_2,z_3,z_4)\bigg]}_{:=\mathrm{I}_4}.
\end{aligned}
\end{equation}
Many contributions in \eqref{Buckling3Noises} are estimating the same way as their counterparts in Proposition \ref{X2Estimates}. Therefore, for those terms, we only point out the differences and the corresponding estimates.

\medskip

$\mathrm{I}_1$ is treated the same way as its counterpart in \eqref{CounterpartsI1X3} using Young's inequality and $\int \psi^\lambda \lesssim 1$, leading to 
\begin{equation*}\label{MomentBoundX3Reduced}
\mathbb{E}\big[\mathbb{X}^{2p}_3\big]\lesssim 1+\mathrm{I}_2+\mathrm{I}_3+\mathrm{I}_4.
\end{equation*}
We treat $\mathrm{I}_2$ using \eqref{EquationFirstOrderDeriv}: there exist $\Gamma^{\lambda,\mathrm{1st}}_{\delta,1}\in \mathcal{K}_3(\psi^\lambda\int \psi^\lambda)$ and $\Gamma^{\lambda,\mathrm{1st}}_{\delta,2}\in \mathcal{K}_2(\psi^\lambda_1)$, where for any $i\geq 1$
$$\psi^\lambda_i(z'):=\psi^\lambda(z') \int \dd z\, \psi^\lambda(z)\, \big(1+|\log(|x'-x|)\big)^i,$$
such that
\begin{equation*}\label{I2FORX3LEMMA}
\mathrm{I}_2=(2p-1)\mathbb{E}\Big[\mathbb{X}^{2(p-1)}_3\big(\mathbb{X}_3(\Gamma^{\lambda,\mathrm{1st}}_{\delta,1})+\mathbb{X}_2(\Gamma^{\lambda,\mathrm{1st}}_{\delta,2})\big)\Big].
\end{equation*}
Furthermore, by \eqref{RefinementMomentBoundInte} in Proposition \ref{X2Estimates}, we have for any $p\geq 1$
\begin{equation*}\label{BoundX2WithPropo11}
\mathbb{E}^\frac{1}{p}\Big[\big|\mathbb{X}_2(\Gamma^{\lambda,\mathrm{1st}}_{\delta,2})\big|^p\Big]\lesssim |\log(\lambda)|^2.
\end{equation*}
Proceeding as for \eqref{MomentBoundX2Reduced2}, we obtain 
\begin{equation}\label{EstiForIducX3FirstReduc}
\mathbb{E}\big[\mathbb{X}^{2p}_3\big]\lesssim |\log(\lambda)|^{2p}+\mathbb{E}\Big[\mathbb{X}^{2(p-1)}_3\mathbb{X}_3(\Gamma^{\mathrm{\lambda,1st}}_{\delta,1})\Big]+\mathrm{I}_3+\mathrm{I}_4.
\end{equation}
For $\mathrm{I}_3$, we use this time \eqref{EquationSecondOrderDeriv} in Lemma \ref{LinkMalliavinDerive}: there exist $\Gamma^{\lambda,\mathrm{2nd}}_{\delta,1}\in \mathcal{K}_3(\psi^\lambda\int \psi^\lambda)$, $\Gamma^{\lambda,\mathrm{2nd}}_{\delta,2}\in \mathcal{K}_{2}(\psi^\lambda_{1})$, $\Gamma^{\lambda,\mathrm{2nd}}_{\delta,3}\in \mathcal{K}_{1}(\psi^\lambda_{2})$ as well as $\mathbb{R}^{1+2}\ni z\mapsto (\Gamma^{\lambda,\mathrm{2nd}}_{\delta,4}(z),\Gamma^{\lambda,\mathrm{2nd}}_{\delta,5}(z))\in \mathcal{K}_{3}(\psi^\lambda)\times \mathcal{K}_{2}(\psi^\lambda \,\mathcal{C}_\delta(\cdot,z))$ such that
\begin{equation}\label{I3ForProofX3}
\begin{aligned}
\mathrm{I}_3\lesssim &\mathbb{E}\Big[\mathbb{X}^{2(p-1)}_3\big(\mathbb{X}_3(\Gamma^{\lambda,\mathrm{2nd}}_{\delta,1})+\mathbb{X}_{2}(\Gamma^{\lambda,\mathrm{2nd}}_{\delta,2})\big)+\mathbb{X}^{2(p-1)}_3|\mathbb{X}_{1}(\Gamma^{\lambda,\mathrm{2nd}}_{\delta,3})|\Big]\\
&+(p-1)\mathbb{E}\bigg[\vert\mathbb{X}_3\vert^{2p-3}\int \dd z\, \eta^\lambda(z)\, \big(\vert \mathbb{X}_3(\Gamma^{\lambda,\mathrm{2nd}}_{\delta,4}(z))\vert +\vert \mathbb{X}_{2}(\Gamma^{\lambda,\mathrm{2nd}}_{\delta,5}(z))\vert\big)^2\bigg].
\end{aligned}
\end{equation}
We then proceed as for \eqref{I3FORX2WLEMMA} where we use the estimates: 
\begin{itemize}
\item[(i)] From Proposition \ref{EstimateForOneNoise}
$$\mathbb{E}^\frac{1}{p}\Big[\big|\mathbb{X}_1(\Gamma^{\lambda,\mathrm{2nd}}_{\delta,3})\big|^p\Big]\lesssim \int \dd z\, \psi^\lambda_2+\bigg(\int \dd z\, \dd z'\, \psi^\lambda_2(z)\psi^\lambda_2(z')\big(1+|\log(|x-x'|)|\big)\bigg)^\frac{1}{2}\lesssim |\log(\lambda)|^{\frac{5}{2}},$$
\item[(ii)] From Proposition \ref{X2Estimates}
$$\bigg(\int \dd z\, \psi^\lambda(z)\, \mathbb{E}^{\frac{3}{2p}}\Big[\big|\mathbb{X}_2(\Gamma^{\lambda,\mathrm{2nd}}_{\delta,5}(z))\big|^{\frac{4}{3}p}\Big]\bigg)^{\frac{2}{3}p}\lesssim \bigg(\fint_{\bb_\lambda}\dd x\, (\mathrm{C}(\lambda,x))^2\bigg)^{\frac{2}{3}p}\lesssim |\log(\lambda)|^{2p},$$
\end{itemize}
which, from \eqref{EstiForIducX3FirstReduc}, leads to: for any $p\in \{1\}\cup [\frac{3}{2},\infty)$
\begin{equation}\label{EstiForIducX3SndReduc}
\mathbb{E}\big[\mathbb{X}^{2p}_3\big]\lesssim \mathbb{E}\Big[\mathbb{X}^{2(p-1)}_3\mathbb{X}_3(\Gamma^{\lambda,\mathrm{1st}}_{\delta,1}+\Gamma^{\lambda,\mathrm{2nd}}_{\delta,1})\Big]+(p-1)\mathbb{E}\bigg[\vert \mathbb{X}_3\vert^{2p-3}\int \dd z\, \psi^\lambda(z)\vert \mathbb{X}_3(\Gamma^{\lambda,\mathrm{2nd}}_{\delta,4}(z))\vert^2\bigg]+|\log(\lambda)|^{\frac{5}{2}p}+\mathrm{I}_4.
\end{equation}
We now turn to the estimate of $\mathrm{I}_4$. We use \eqref{EquationThirdOrderDeriv} in Lemma \ref{LinkMalliavinDerive}: there exist $\Gamma^{\lambda,\mathrm{3rd}}_{\delta,1}\in \mathcal{K}_3(\psi^\lambda\int \psi^\lambda)$, $\Gamma^{\lambda,\mathrm{3rd}}_{\delta,2}\in \mathcal{K}_{2}(\psi^\lambda_{1})$, $\Gamma^{\lambda,\mathrm{3rd}}_{\delta,3}\in \mathcal{K}_{1}(\psi^\lambda_{2})$, $\mathbb{R}^{1+2}\ni z\mapsto (\Gamma^{\lambda,\mathrm{3rd}}_{\delta,4}(z),\Gamma^{\lambda,\mathrm{3rd}}_{\delta,5}(z))\in \mathcal{K}_{3}(\psi^\lambda)\times \mathcal{K}_{2}(\eta \,\mathcal{C}_\delta(\cdot,z))$, $\Gamma^{\lambda,\mathrm{3rd}}_{\delta,6}\in \mathcal{K}_{0}(\psi^\lambda_3)$, $\mathbb{R}^{1+2}\ni z\mapsto (\Gamma^{\lambda,\mathrm{3rd}}_{\delta,7}(z),\Gamma^{\lambda,\mathrm{3rd}}_{\delta,8}(z))\in \mathcal{K}_{3}(\psi^\lambda)\times \mathcal{K}_{2}(\psi^\lambda\,\mathcal{C}_{\delta,C}(\cdot,z))$, $z\mapsto (\Gamma^{\lambda,\mathrm{3rd}}_{\delta,9}(z),\Gamma^{\lambda,\mathrm{3rd}}_{\delta,11}(z))\in \mathcal{K}_3(\psi^\lambda)\times \mathcal{K}_3(\psi^\lambda)$, $z\mapsto (\Gamma^{\lambda,\mathrm{3rd}}_{\delta,10}(z),\Gamma^{\lambda,\mathrm{3rd}}_{\delta,12}(z),\Gamma^{\lambda,\mathrm{3rd}}_{\delta,13}(z))\in \mathcal{K}_{2}(\psi^\lambda\,\mathcal{C}_{\delta,C}(\cdot,z))\times \mathcal{K}_{2}(\psi^\lambda\, \mathcal{C}_{\delta,C}(\cdot,z))\times \mathcal{K}_{1}(\psi^\lambda\,\mathcal{C}^2_{\delta,C}(\cdot,z))$ such that
\begin{equation}\label{EquationThirdOrderDerivProofX3}
\begin{aligned}
\mathrm{I}_4\lesssim \, & \mathbb{E}\bigg[\mathbb{X}^{2(p-1)}_3\big(\mathbb{X}_n(\Gamma^{\lambda,\mathrm{3rd}}_{\delta,1})+\mathbb{X}_{2}(\Gamma^{\lambda,\mathrm{3rd}}_{\delta,2})\big)+\mathbb{X}^{2(p-1)}_3|\mathbb{X}_{1}(\Gamma^{\lambda,\mathrm{3rd}}_{\delta,3})|\bigg]\\
&\quad +(p-1)\mathbb{E}\bigg[\vert\mathbb{X}_3\vert^{2p-3}\int \dd z\, \psi^\lambda(z)\, \big(\vert \mathbb{X}_3(\Gamma^{\lambda,\mathrm{3rd}}_{\delta,4}(z))\vert +\vert \mathbb{X}_{2}(\Gamma^{\lambda,\mathrm{3rd}}_{\delta,5}(z))\vert\big)^2\bigg]\\
&\quad + \mathbb{E}\bigg[\mathbb{X}^{2(p-1)}_3|\mathbb{X}_{0}(\Gamma^{\lambda,\mathrm{3rd}}_{\delta,6})|\bigg]\\
&\quad + (p-1)(2p-3)\mathbb{E}\bigg[\mathbb{X}^{2p-4}_3\int \dd z\, \psi^\lambda(z)\, \big(\vert \mathbb{X}_3(\Gamma^{\lambda,\mathrm{3rd}}_{\delta,7}(z))\vert +\vert \mathbb{X}_{2}(\Gamma^{\lambda,\mathrm{3rd}}_{\delta,8}(z))\vert\big)^3\bigg]\\
&\quad + (p-1)\mathbb{E}\bigg[|\mathbb{X}_3|^{2p-3}\int \dd z\, \psi^\lambda(z)\, \big(|\mathbb{X}_3(\Gamma^{\lambda,\mathrm{3rd}}_{\delta,9}(z))|+|\mathbb{X}_{2}(\Gamma^{\lambda,\mathrm{3rd}}_{\delta,10}(z))|\big)\\
&\qquad\qquad\qquad\qquad\qquad\qquad\qquad\qquad\big(|\mathbb{X}_3(\Gamma^{\lambda,\mathrm{3rd}}_{\delta,11}(z))|+|\mathbb{X}_{2}(\Gamma^{\lambda,\mathrm{3rd}}_{\delta,12}(z))|+|\mathbb{X}_{1}(\Gamma^{\lambda,\mathrm{3rd}}_{\delta,13}(z))|\big)\bigg].
\end{aligned}
\end{equation}
The two first r. h. s. terms are controlled the same way as $\mathrm{I}_3$ in \eqref{I3ForProofX3}. For the third r. h. s. term, we recall that $\Gamma^{\lambda,\mathrm{3rd}}_{\delta,6}\in \mathcal{K}_{0}(\psi^\lambda_3)$ and therefore get
\begin{equation}\label{EstimateX0ForI4X3}
|\mathbb{X}_0(\Gamma^{\lambda,\mathrm{3rd}}_{\delta,6})|\lesssim \int \dd z\, \psi^\lambda_3(z)=\int \dd z\, \dd z'\, \psi^\lambda(z)\, \psi^\lambda(z')\, \big(1+|\log(|x'-x|)|\big)^3\lesssim |\log(\lambda)|^{3}.
\end{equation}
For the fourth r.h.s. term, in the regime $p\geq 2$, we use Young's inequality with exponents $(\frac{2p}{2p-4},\frac{p}{2})$ and Proposition \ref{X2Estimates} in the form: for any $\varepsilon>0$, there exist $C_\varepsilon>0$ such that
\begin{equation}\label{EstimateForI4X32}
\begin{aligned}
\mathbb{E}\bigg[\vert\mathbb{X}_3(\Gamma_1)\vert^{2p-4}\int \dd z\, \psi^\lambda(z)\, \vert \mathbb{X}_{2}(\Gamma^{\lambda,\mathrm{3rd}}_{\delta,8}(z))\vert^3\bigg] &\leq \varepsilon \mathbb{E}\big[\mathbb{X}^{2p}_3\big]+C_\varepsilon\mathbb{E}\bigg[\Big(\int \dd z\, \psi^\lambda(z)\, |\mathbb{X}_2(\Gamma^{\lambda,\mathrm{3rd}}_{\delta,8}(z))|^3\Big)^{\frac{p}{2}}\bigg]\\
&\leq \varepsilon \mathbb{E}\big[\mathbb{X}^{2p}_3\big]+C_\varepsilon\Big(\int \dd z\, \psi^\lambda(z)\, \mathbb{E}^\frac{2}{p}\big[|\mathbb{X}_2(\Gamma^{\lambda,\mathrm{3rd}}_{\delta,8}(z))|^\frac{3p}{2}\big]\Big)^{\frac{p}{2}}\\
&\stackrel{\eqref{RefinementMomentBound}}{\lesssim} \varepsilon \mathbb{E}\big[\mathbb{X}^{2p}_3\big]+C_\varepsilon|\log(\lambda)|^{\frac{9}{4}p}.
\end{aligned}
\end{equation}
For the last r.h.s. term, in the regime $p\geq \frac{3}{2}$, we use Young's inequality with exponents $(\frac{2p}{2p-3}, \frac{2}{3}p)$ and Proposition \ref{EstimateForOneNoise} and Proposition \ref{X2Estimates} in the form of: for any $\varepsilon>0$, there exists $C_\varepsilon>0$ such that 
\begin{equation}\label{EstimateForI4X33}
\begin{aligned}
&\mathbb{E}\bigg[\int \dd z\, \psi^\lambda(z)\, |\mathbb{X}_2(\Gamma^{\lambda,\mathrm{3rd}}_{\delta,10}(z))|\big(|\mathbb{X}_2(\Gamma^{\lambda,\mathrm{3rd}}_{\delta,12}(z))|+|\mathbb{X}_1(\Gamma^{\lambda,\mathrm{3rd}}_{\delta,13}(z))|\big)\bigg]\\
&\leq \varepsilon \mathbb{E}\big[\mathbb{X}^{2p}_3\big]+C_\varepsilon\bigg(\int \dd z\, \psi^\lambda(z)\,\mathbb{E}^\frac{3}{2p}\Big[|\mathbb{X}_2(\Gamma^{\lambda,\mathrm{3rd}}_{\delta,10}(z))|^{\frac{2}{3}p}\big(|\mathbb{X}_2(\Gamma^{\lambda,\mathrm{3rd}}_{\delta,12}(z))|+|\mathbb{X}_1(\Gamma^{\lambda,\mathrm{3rd}}_{\delta,13}(z))|\big)^{\frac{2}{3}p}\Big]\bigg)^{\frac{2}{3}p}\\
&\stackrel{\eqref{StochasticEstimateOneNoise},\eqref{RefinementMomentBound}}{\leq}\varepsilon \mathbb{E}\big[\mathbb{X}^{2p}_3\big]+C_\varepsilon |\log(\lambda)|^{3p}.
\end{aligned}
\end{equation}
To conclude, the combination of \eqref{EstiForIducX3SndReduc}, \eqref{EquationThirdOrderDerivProofX3}, \eqref{EstimateX0ForI4X3}, \eqref{EstimateForI4X32} and \eqref{EstimateForI4X33} yields, for any $p\in \{1\}\cup [\frac{3}{2},\infty)$
\begin{equation}\label{FinalEstiInducX3}
\begin{aligned}
&\mathbb{E}\big[\mathbb{X}^{2p}_3\big]\\
&\lesssim \mathbb{E}\Big[\mathbb{X}^{2(p-1)}_3\mathbb{X}_3(\Gamma^{\lambda,\mathrm{1st}}_{\delta,1}+\Gamma^{\lambda,\mathrm{2nd}}_{\delta,1}+\Gamma^{\lambda,\mathrm{3rd}}_{\delta,1})\Big]+(p-1)\mathbb{E}\bigg[\vert \mathbb{X}_3\vert^{2p-3}\int \dd z\, \psi^\lambda(z)\big(\vert \mathbb{X}_3(\Gamma^{\lambda,\mathrm{2nd}}_{\delta,4}(z))\vert^2+\vert \mathbb{X}_3(\Gamma^{\lambda,\mathrm{3rd}}_{\delta,4}(z))\vert^2\big)\bigg]\\
&+(p-1)\mathbb{E}\bigg[|\mathbb{X}_3|^{2p-3}\int \dd z\, \psi^\lambda(z)\, |\mathbb{X}_{3}(\Gamma^{\lambda,\mathrm{3rd}}_{\delta,9}(z))||\mathbb{X}_{3}(\Gamma^{\lambda,\mathrm{3rd}}_{\delta,11}(z))|\bigg]\\
&+(p-1)\mathbb{E}\bigg[|\mathbb{X}_3|^{2p-3}\int \dd z\, \psi^\lambda(z)\, \big(|\mathbb{X}_{3}(\Gamma^{\lambda,\mathrm{3rd}}_{\delta,11}(z))||\mathbb{X}_{2}(\Gamma^{\lambda,\mathrm{3rd}}_{\delta,10}(z))|+|\mathbb{X}_{3}(\Gamma^{\lambda,\mathrm{3rd}}_{\delta,9}(z))|(|\mathbb{X}_{2}(\Gamma^{\lambda,\mathrm{3rd}}_{\delta,12}(z))|+|\mathbb{X}_{1}(\Gamma^{\lambda,\mathrm{3rd}}_{\delta,13}(z))|)\bigg]\\
&+(p-1)(2p-3)\underbrace{\mathbb{E}\bigg[\mathbb{X}^{2p-4}_3\int \dd z\, \psi^\lambda(z)\, |\mathbb{X}_3(\Gamma^{\lambda,\mathrm{3rd}}_{\delta,7}(z))|^3\bigg]}_{:=\mathrm{I}_5}\\
&+|\log(\lambda)|^{3p}.
\end{aligned}
\end{equation}
We now argue by induction over integers $p\geq 1$ and split the proof into two steps.

\medskip

{\sc Step 1. Initialization. }We note that \eqref{FinalEstiInducX3} requires to initialize with $p\in \{1,2,3\}$. The structure of the first four r.h.s. terms of \eqref{FinalEstiInducX3} is identical to that of the two first r.h.s terms of \eqref{EstimateForInducX2} and therefore are bounded the exact same way. In particular, we deduce \eqref{MomentBounds3Noises} for $p\in \{1,\frac{3}{2}\}$.

\medskip

In the following, we only explain how to handle $\mathrm{I}_5$ for $p\in \{2,3\}$. First, note that for $p=2$, $\mathrm{I}_5$ is controlled with the third order moment of $\mathbb{X}_3(\Gamma^{\lambda,\mathrm{3rd}}_{\delta,7}(z))$. Second, we apply \eqref{FinalEstiInducX3} with $p=\frac{9}{4}$ and use Young's inequality with exponents $(9,\frac{9}{8})$ to the effect of: for any $\varepsilon>0$, there exists $C_\varepsilon>0$ such that
\begin{align*}
\mathrm{I}_5=&\mathbb{E}\bigg[\mathbb{X}^\frac{1}{2}_3\int \dd z\, \psi^\lambda(z)\, |\mathbb{X}_3(\Gamma^{\lambda,\mathrm{3rd}}_{\delta,7}(z))|^3\bigg]\\
\leq & \varepsilon \mathbb{E}[\mathbb{X}^\frac{9}{2}]+C_\varepsilon \bigg(\int \dd z\, \psi^\lambda(z)\, \mathbb{E}^\frac{8}{9}\big[|\mathbb{X}_3(\Gamma^{\lambda,\mathrm{3rd}}_{\delta,7}(z))|^{\frac{27}{8}}\big]\bigg)^{\frac{9}{8}}.
\end{align*}
Noticing that $\frac{27}{8}\leq 4$, we deduce that 
$$\mathrm{I}_5\leq  \varepsilon \mathbb{E}[\mathbb{X}^\frac{9}{2}]+C_\varepsilon |\log(\lambda)|^{\frac{81}{16}},$$
and we can absorb the first r.h.s. into the l.h.s. of \eqref{FinalEstiInducX3}.

\medskip

To conclude, we can deduce the case $p=3$ using Young's inequality with exponents $(3,\frac{3}{2})$ to the effect of: for any $\varepsilon>0$, there exists $C_\varepsilon>0$ such that
\begin{align*}
\mathrm{I}_5 &\leq \varepsilon \mathbb{E}[\mathbb{X}^6_3\big]+C_\varepsilon\bigg(\int \dd z\, \psi^\lambda(z)\, \mathbb{E}^\frac{2}{3}\big[|\mathbb{X}_3(\Gamma^{\lambda,\mathrm{3rd}}_{\delta,7}(z))|^{\frac{9}{2}}\big]\bigg)^\frac{3}{2}\leq \varepsilon \mathbb{E}[\mathbb{X}^6_3\big]+C_\varepsilon|\log(\lambda)|^{\frac{27}{4}},
\end{align*}
and we can absorb the first r.h.s. into the l.h.s. of \eqref{FinalEstiInducX3}.

\medskip

{\sc Step 2. Induction step. }We assume that $p\geq 4$ and we show the induction step. Once again, the first four r.h.s. terms of \eqref{FinalEstiInducX3} are trated the same way as in \eqref{InducStepForX2}. For $\mathrm{I}_5$, we use Young's inequality with exponents $(\frac{2p}{2p-4},\frac{p}{2})$ in the form of: for any $\varepsilon>0$, there exists a constant $C_\varepsilon>0$ such that 
$$\mathrm{I}_5 \leq \varepsilon \mathbb{E}[\mathbb{X}^{2p}_3\big]+C_\varepsilon\bigg(\int \dd z\, \psi^\lambda(z)\, \mathbb{E}^\frac{2}{p}\big[|\mathbb{X}_3(\Gamma_{\delta,7}(z))|^{\frac{3}{2}p}\big]\bigg)^\frac{p}{2}.$$
Noticing that  $\frac{3}{2}p\leq 2(p-1)$, we then use the induction hypothesis 
$$\mathbb{E}^\frac{2}{p}\big[|\mathbb{X}_3(\Gamma_{\delta,7}(z))|^{\frac{3}{2}p}\big]\leq |\log(\lambda)|^{\frac{9}{4}p},$$
and conclude that 
$$\mathrm{I}_5 \leq \varepsilon \mathbb{E}[\mathbb{X}^{2p}_3\big]+C_\varepsilon|\log(\lambda)|^{\frac{9}{4}p},$$
where we can absorb the first r. h. s. term into the l. h. s. of \eqref{FinalEstiInducX3}.
\end{proof}
\subsection{Proof of Theorem \ref{MainResult}}
In this subsection, we provide the proof of Theorem \ref{MainResult}. At first sight, according to Proposition \ref{EstimateForOneNoise}, Proposition \ref{X2Estimates} and Proposition \ref{EstiX3}, we only have to show that $K : (z,z')\in \mathbb{R}^{1+2}\times \mathbb{R}^{1+2}\mapsto K_z(z-z')$ belongs to the class $\mathcal{K}$ defined in Definition \ref{DefClassKernel}. However, under the regularity condition on $m$ in Assumption \ref{AssumptionCoef}, this is in general not the case. We thus argue by approximation: let $\{\zeta_q\}_{q\geq 1}$ be an approximation of unity and set, for any $q\geq 1$, $m_q:=\zeta_q\star m$ and $K_q : (z,z')\in \mathbb{R}^{1+2}\times \mathbb{R}^{1+2} \mapsto K_{z,q}(z-z')$ defined in \eqref{DefKernelFrozen} with $g$ replaced by $g_q:=m_q \star \xi$. We build accordingly $\widehat{\lollipop}_{\delta,q}$, $\overline{\cherry}_{\delta,q}$ and $\overline{\chickenfoot}_{\delta,q}$ with $K_z$ replaced by $K_{z,q}$. Our task is then to show that $K_{q}$ belongs to $\mathcal{K}$ with \eqref{AssumptionKernelFunction1Noise} holding uniformly over $q$ (which implies that \eqref{StochasticEstimateOneNoise}, \eqref{RefinementMomentBound} and \eqref{MomentBounds3Noises} hold uniformly over $q$ as well) and for any $p<\infty$ 
\begin{equation}\label{PassToTheLimitRenor}
\begin{aligned}
&\mathbb{E}^{\frac{1}{p}}\big[(\widehat{\lollipop}_{\delta,q},\psi^\lambda)^p\big]\underset{q\uparrow\infty}{\longrightarrow} \mathbb{E}^{\frac{1}{p}}\big[(\widehat{\lollipop}_{\delta},\psi^\lambda)^p\big],\\
&\mathbb{E}^{\frac{1}{p}}\big[(\overline{\cherry}_{\delta,q},\psi^\lambda)^p\big]\underset{q\uparrow\infty}{\longrightarrow} \mathbb{E}^{\frac{1}{p}}\big[(\overline{\cherry}_{\delta},\psi^\lambda)^p\big] \quad\text{and}\quad \mathbb{E}^{\frac{1}{p}}\big[(\overline{\chickenfoot}_{\delta,q},\psi^\lambda)^p\big] \underset{q\uparrow\infty}{\longrightarrow} \mathbb{E}^{\frac{1}{p}}\big[(\overline{\chickenfoot}_{\delta},\psi^\lambda)^p\big].
\end{aligned}
\end{equation}
The fact that $K_q\in \mathcal{K}$ is a direct consequence of its definition: we directly have that $K_q=\bar{K}_q(\xi,\cdot,\cdot)$ with for any $T\in X$ and $z,z'\in \mathbb{R}^{1+2}$
\begin{equation}\label{KernelXDep}
\bar{K}_q(T,z,z')=\mathds{1}_{t\geq t'}\frac{1}{4\pi (t-t')\sqrt{\text{det}(A_{\text{sym}}(m_q\star T(z)))}}\exp\bigg(-\frac{(x-x')\cdot A^{-1}_{\text{sym}}(m_q\star T(z))(x-x')}{t-t'}\bigg).
\end{equation}
We observe that it has the following structure 
\begin{equation}\label{StructureKQ}
\bar{K}_q(T,z,z')=F_{z,z'}(m_q\star T(z))
\end{equation}
where $F_{z,z'} : \mathbb{R}\rightarrow \mathbb{R}$ is smooth and for any $n\geq 1$ and some $C_n>0$, $\sup_{\lambda\in \mathbb{R}}\vert F^{(n)}_{z,z'}(\lambda)\vert\lesssim \mathds{1}_{t\geq t'}G_{\sqrt{t-t'}}(\frac{1}{\sqrt{C_n}}(x-x'))$. The Fréchet differentiability and the estimate \eqref{ConstantDefClass} then follows by a simple computation using that 
\begin{equation}\label{EstimateConvoMT}
\vert m_q\star T(z)\vert\lesssim_q (1+\vert z\vert^2)^\sigma \|T\|_X,
\end{equation}
and this is where the integrability condition \eqref{AssIntM} kicks in. We obtain \eqref{EstimateConvoMT} by a direct calculation using the space-time Fourier transform. First, from the Cauchy-Schwarz inequality it holds\footnote{we recall that $S:=(\partial_t+1-\Delta)^{-\sigma}(1+\vert\cdot\vert^2)^{-\sigma}$ with $\sigma>\frac{1}{2}$ and 
$S^\star$ denotes its adjoint w. r. t. the $\LL^2(\mathbb{R}^{1+2})$-inner product}
\begin{equation}\label{BoundWithL2CONVO}
\vert m_q\star T(z)\vert\leq \|(S^\star)^{-1}m_q(z-\cdot)\|_{\LL^{2}(\mathbb{R}^{1+2})}\|T\|_X\lesssim \|(1+\vert\cdot\vert^2)^\sigma \mathcal{F}((1+\vert \cdot\vert^2)^\sigma m_q(z-\cdot))\|_{\LL^2(\mathbb{R}^{1+2})}\|T\|_X.
\end{equation}
Then, note that for any $\eta\in \mathbb{R}^{1+2}$
\begin{align*}
\mathcal{F}((1+\vert \cdot\vert^2)^\sigma m_q(z-\cdot))(\eta)&=\mathcal{F}\bigg((1+\vert \cdot\vert^2)^\sigma\int \dd y\, \zeta_q (z-\cdot-y)m(y)\bigg)(\eta)\\
&=\int \dd y\,\mathcal{F}\big((1+\vert \cdot\vert^2)^\sigma \zeta_q(z-\cdot-y)\big)(\eta)m(y).
\end{align*}
We then use that for any integer $k\geq 1$, we have by regularity
\begin{align*}
\big\vert \mathcal{F}\big((1+\vert \cdot\vert^2)^\sigma \zeta_q(z-\cdot-y)\big)(\eta)\big\vert &\lesssim (1+\vert \eta\vert^2)^{-k} \big\vert \mathcal{F}\big((\partial_t+1-\Delta)^k(1+\vert \cdot\vert^2)^\sigma \zeta_q(z-\cdot-y)\big)(\eta)\big\vert\\
&\lesssim_{k,q}  (1+\vert \eta\vert^2)^{-k} (1+\vert z\vert^2)^\sigma (1+\vert y\vert^2)^\sigma,
\end{align*}
so that 
$$\big\vert \mathcal{F}((1+\vert \cdot\vert^2)^\sigma m_q(z-\cdot))(\eta)\big\vert \lesssim_k (1+\vert \eta\vert^2)^{-k} (1+\vert z\vert^2)^\sigma \int \dd y\, (1+\vert y\vert^2)^\sigma \vert m(y)\vert\stackrel{\eqref{AssIntM}}{\lesssim} (1+\vert \eta\vert^2)^{-k} (1+\vert z\vert^2)^\sigma,$$
which concludes by plugging it in \eqref{BoundWithL2CONVO}.

\medskip

The second item of \eqref{AssumptionKernelFunction1Noise} holds uniformly over $q$ (the first item follows directly from \eqref{KernelXDep} and ellipticity of $A$) from a simple computation as well using \eqref{StructureKQ}: we have for any $z,z',y'_1,\cdots,y'_k\in \mathbb{R}^{1+2}$
\begin{equation}\label{MalliavinFrozenN}
\mathrm{D}_{y'_1}\cdots\mathrm{D}_{y'_k}\, K_{q}(z,z')=F^{(k)}_{z,z'}(g_q(z))\prod_{i=1}^k m_q (z-y'_i).
\end{equation}
Therefore, we deduce that
\begin{equation*}
\begin{aligned}
\int \dd y'_1\cdots\dd y'_k\, \rho_\delta(y'_1-z'_1)\cdots \rho_\delta(y'_k-z'_k)\, \mathrm{D}_{y'_1}\cdots\mathrm{D}_{y'_k} K_{q}(z,z')&=F^{(k)}_{z,z'}(g_q(z))\prod_{i=1}^k \rho_\delta \star m_q (z-z'_i)\\
&\lesssim_{\|m\|_{\LL^\infty(\mathbb{R}^{1+2})}} \mathds{1}_{t\geq t'}G_{\sqrt{t-t'}}\Big(\tfrac{1}{\sqrt{C_k}}(x-x')\Big),
\end{aligned}
\end{equation*}
which is uniform over $q$. 

\medskip

Finally, we show \eqref{PassToTheLimitRenor}. We only justify the second limit, the first and the thrid are proven the same way. We first have by the triangle inequality
\begin{align*}
\mathbb{E}^\frac{1}{p}\big[(\overline{\cherry}_{\delta,q}-\overline{\cherry}_{\delta},\psi^\lambda)^p\big]\leq & \mathbb{E}^\frac{1}{p}\bigg[\Big(\int \dd z\,\psi^\lambda(z)\int \dd z'\, \dd z''\, \big(K_{z,q}(z-z')-K_z(z-z')\big)K_{z,q}(z-z'')\xi_\delta(z')\diamond \xi_\delta(z'')\Big)^p\bigg]\\
&+\mathbb{E}^\frac{1}{p}\bigg[\Big(\int \dd z\,\psi^\lambda(z)\int \dd z'\, \dd z''\, K_z(z-z')\big(K_{z,q}(z-z'')-K_z(z-z'')\big)\xi_\delta(z')\diamond \xi_\delta(z'')\Big)^p\bigg].
\end{align*}
We then use in order \eqref{StructureKQ} in the form of $\vert K_{z,q}(z-z')-K_z(z-z')\vert\lesssim \mathds{1}_{t\geq t'} G_{\sqrt{t-t'}}(\frac{1}{\sqrt{C}}(x-x'))\vert g_q(z)-g(z)\vert$, $\xi\in \cc^{-2^{-}}$ in the form of $\mathbb{E}^{\frac{1}{p}}[(\xi_\delta(z'))^p]\lesssim \delta^{-2^{-}}$ and Gaussianity in the form of $\mathbb{E}^\frac{1}{2p}\big[\vert g_n(z)-g(z)\vert^{2p}\big]\lesssim_p \mathbb{E}^\frac{1}{2}\big[\vert g_q(z)-g(z)\vert^{2}\big]= \|m_q(z-\cdot)-m(z-\cdot)\|_{\LL^2(\mathbb{R}^{1+2})}$ to obtain 
\begin{align*}
\mathbb{E}^\frac{1}{p}\big[(\overline{\cherry}_{\delta,q}-\overline{\cherry}_{\delta},\psi^\lambda)^p\big]&\lesssim_\delta  \int \dd z\, \vert \psi^\lambda(z)\vert \,\mathbb{E}^\frac{1}{2p}\big[\vert g_q(z)-g(z)\vert^{2p}\big]\\
&\lesssim_{\delta,p} \int \dd z\, \vert \psi^\lambda(z)\,\vert \|m_q (z-\cdot)-m(z-\cdot)\|_{\LL^2(\mathbb{R}^{1+2})},
\end{align*}
which goes to $0$ as $q\uparrow\infty$ using $m\in \cc^\alpha(\mathbb{R}^{1+2})\cap \LL^2(\mathbb{R}^{1+2})$ and Lebesgue's convergence theorem.
\subsection{Estimate of the regular part}
In this subsection, we show that the moments of the regular parts $\{\mathcal{R}_\delta\}_{\delta>0}$ are uniformly controlled in $\delta>0$.
\begin{proposition}[Estimate of the regular part]\label{PropRegPart}
For any $\delta>0$, $\lambda\in (0,e^{-1}]$, $m\in \{1,2,3\}$, test functions $\psi^\lambda:=\lambda^{-4}\psi(\frac{\cdot}{\lambda^2},\frac{\cdot}{\lambda})$ with $\mathrm{supp}\,\psi\subset \cc_1$ and $\|\psi\|_{\cc^0(\mathbb{R}^{1+2})}\leq 1$, it holds
\begin{equation}\label{EstimateRegularPart}
\mathbb{E}^{\frac{1}{2p}}\bigg[\Big(\int \dd z\, \psi^\lambda(z)\Big(\int \dd z'\, R(z,z')\xi_\delta(z')\Big)^m\Big)^{2p}\bigg]\lesssim_p 1\quad\text{for any $p\geq 1$},
\end{equation}
where we recall that $R$ is defined in \eqref{RegPartDef}.
\end{proposition}
\begin{proof}
%
%
%
The proof of \eqref{EstimateRegularPart} follows the lines of the proofs of Proposition \ref{EstimateForOneNoise}, Proposition \ref{X2Estimates} and Proposition \ref{EstiX3}, with the following three modifications:
\medskip

\begin{itemize}
\item[(i)]To obtain the necessary regularity on $m$ to apply Lemma \ref{IntWRTWick}, we argue by approximation using a mollifier: let $\{\zeta_q\}_{q\geq 1}$ be an approximation of unity and set, for any $q\geq 1$, $m_q:=\zeta_q\star m$, $a_q:=\zeta_q\star A(m_q\star \xi)$ and define $R_q$ and $\Gamma_q$ as in \eqref{RegPartDef} with $a$ replaced by $a_q$. Under this mollification, the Malliavin differentiability of $R_q$ is determined by the one of $a_q$ through \eqref{DefKernelFrozen} and the PDE solved by $\Gamma_q$, therefore it follows easily that $R_q$ is infinitely Malliavin differentiable in the sense of Appendix \ref{AppendixGIBP}. 

\medskip
%


\item[(ii)]$R_q$ belongs to a version of the class $\mathcal{K}$ in Definition \ref{DefClassKernel}, more precisely:

\medskip

\begin{itemize}
\item[(a)]There exists $\bar{R}_q : X\times \mathbb{R}^{1+2}\times \mathbb{R}^{1+2}\rightarrow (0,\infty)$ such that $R_q=\bar{R}_q(\xi,\cdot,\cdot)$ and, for any $z,z'\in \mathbb{R}^{1+2}$, $\bar{R}_q(\cdot,z,z')$ satisfies the assumptions of Lemma \ref{IntWRTWick} with constants $C_q(z,z')$ (in \eqref{AssumptionInteOnF}) and $C_{n,q}(z,z')$ (in \eqref{AssumptionInteOnDNF}) which satisfy for some $\sigma_n,C_n>0$
\begin{equation*}\label{ConstantDefClass2}
\max\big\{C_q(z,z'),C_{n,q}(z,z')\big\}\lesssim_{n,q} \mathds{1}_{t\geq t'}(1+\vert z\vert^2+\vert z'\vert^2)^{\sigma_n} (t-t')G_{\sqrt{t-t'}}\Big(\tfrac{1}{\sqrt{C_n}}(x-x')\Big).
\end{equation*}
\item[(b)]$R_q$ and its Malliavin derivatives satisfy the following bounds (uniformly over $q$): for any $n\geq 1$, there exists $C,C_n>0$ such that for any $\delta>0$, $\alpha'<\alpha$, $z,z',z'_1,\cdots,z'_n\in \mathbb{R}^2$ and $p\geq 1$, it holds\footnote{for a kernel $K: \mathbb{R}^{1+2}\times \mathbb{R}^{1+2}$, we write, for any $i\in\{1,2\}$, $\nabla_i K$ its spatial derivative with respect to the $i^{th}$-slot} for some $\gamma'>0$ depending on $d$ and $\alpha$
\begin{equation}\label{BoundMalliavinDerivativeRest}
\begin{aligned}
& \mathbb{E}^{\frac{1}{p}}\big[\vert R_q(z,z')\vert^p\big]\lesssim_p \mathds{1}_{t\geq t'}(t-t')^{\frac{\alpha'}{2}}G_{\sqrt{t-t'}}\Big(\tfrac{1}{\sqrt{C}}(x-x')\Big),\\
&\quad \mathbb{E}^\frac{1}{p}\bigg[\Big\vert\mathrm{d}^n R_q(z,z').\bigotimes_{i=1}^n \rho_\delta(\cdot-z'_i)\Big\vert^p\bigg]\lesssim_p  \mathds{1}_{t\geq t'}\max\big\{1,(t-t')^{\gamma'}\big\}(t-t')^{\frac{\alpha'}{2}}G_{\sqrt{t-t'}}\Big(\tfrac{1}{\sqrt{C_n}}(x-x')\Big)\\
&\text{and}\quad \mathbb{E}^\frac{1}{p}\bigg[\Big\vert\mathrm{d}^n \nabla_2 R_q(z,z').\bigotimes_{i=1}^n \rho_\delta(\cdot-z'_i)\Big\vert^p\bigg]\lesssim_p  \mathds{1}_{t\geq t'}\max\big\{1,(t-t')^{\gamma'}\big\}(t-t')^{\frac{\alpha'-1}{2}}G_{\sqrt{t-t'}}\Big(\tfrac{1}{\sqrt{C_n}}(x-x')\Big).
\end{aligned}
\end{equation}
\end{itemize}
%
%
%
%
\item[(iii)]The test function $\psi^{\lambda}$ is now deterministic. 
\end{itemize}

\medskip

The main change is in the gain $(t-t')^{\frac{\alpha'}{2}}$ in the two first items of \eqref{BoundMalliavinDerivativeRest} which implies the uniform bound in $\lambda$ and $\delta$. In what follows, we split the proof into four steps. A first step is devoted to establish the necessary PDE characterization of $\bar{R}_q$ and its Malliavin derivatives. Then, the second and third steps are devoted to the proofs of $(a)$ and $(b)$ separately and finally, in a last step, we show how to pass to the limit as $q\uparrow\infty$.

\medskip

{\sc Step 1. A PDE characterization. }We derive in this step a PDE characterization of $\bar{R}_q$ and its Malliavin derivatives. For this, we introduce the following notations: 
\begin{itemize}
\item (Malliavin derivatives). For any $Y:=(y_1,\cdots,y_n)\in \mathbb{R}^{1+2}$ and $\mathrm{I}=\{i_1,\cdots,i_m\}\subset \{1,\cdots,n\}$, we introduce
$$\mathrm{D}^Y_{\mathrm{I}}:=\mathrm{D}_{y_{i_1}}\cdots\mathrm{D}_{y_{i_m}}\quad\text{and}\quad \mathrm{D}^Y_{\emptyset}:=\mathrm{Id}.$$
\item (Finite differences). For any $T,T'\in X$, we write $\delta F(T,T'):= F(T)-F(T')$ to denote the finite differences for a functional $F : X\rightarrow \mathbb{R}$.
\end{itemize}

\medskip

We claim that there exists $\bar{R}_q : X\times \mathbb{R}^{1+2}\times \mathbb{R}^{1+2}$ such that $R_q=\bar{R}_q(\xi,\cdot,\cdot)$ and for any $\mathrm{I}=\{i_1,\cdots,i_m\}\subseteq \{1,\cdots,n\}$, $y_{i_1}\neq\cdots\neq y_{i_m}$ and $T\in X$ it solves in the weak sense 
\begin{equation}\label{EqMalliavinDiffNoFinite}
\begin{aligned}
\partial_t \mathrm{D}^Y_{\mathrm{I}}\bar{R}_q(T,z,\cdot)+\nabla\cdot \bar{a}_q(T)\nabla \mathrm{D}^Y_{\mathrm{I}}\bar{R}_q(T,z,\cdot)=&\sum_{\mathrm{J}\subset \mathrm{I}}\nabla\cdot \mathrm{D}^Y_{\mathrm{J}^c}\bar{a}_q(T)\mathrm{D}^Y_{\mathrm{J}}\bar{R}_q(T,z,\cdot)\\
&+\sum_{\mathrm{J}\subseteq \mathrm{I}}\nabla\cdot \mathrm{D}^Y_{\mathrm{J}^c}\big(\bar{a}_q(T)-\bar{a}_q(T,z)\big)\nabla \mathrm{D}^Y_{\mathrm{J}} \bar{K}_q(T,z,\cdot)\quad \text{in $(0,t)\times \mathbb{R}^2$,}
\end{aligned}
\end{equation}
with $\mathrm{D}^Y_{\mathrm{I}}\bar{R}_q(T,z,(t,\cdot))\equiv 0$ and
\begin{equation}\label{EqMalliavinDiff}
\begin{aligned}
\partial_t \mathrm{D}^Y_{\mathrm{I}}\delta \bar{R}_q(T,T',z,\cdot)+\nabla\cdot \bar{a}_q(T)\mathrm{D}^Y_{\mathrm{I}}\nabla\delta \bar{R}_q(T,T',z,\cdot)=&\nabla\cdot \sum_{\mathrm{J}\subset \mathrm{I}}\mathrm{D}^Y_{\mathrm{J}^c}\bar{a}_q(T)\nabla \mathrm{D}^Y_{\mathrm{J}}\delta \bar{R}_q(T,T',z,\cdot)\\
&+\nabla\cdot\sum_{\mathrm{J}\subseteq \mathrm{I}}\mathrm{D}^Y_{\mathrm{J}^c}\big(\bar{a}_q(T)-\bar{a}_q(T,z)\big)\nabla \mathrm{D}^Y_{\mathrm{J}}\delta \bar{K}_{q}(T,T',z,\cdot)\\
&-\nabla\cdot \sum_{\mathrm{J}\subseteq \mathrm{I}} \mathrm{D}^Y_{\mathrm{J}^c}\delta \bar{a}_q(T,T')\nabla \mathrm{D}^Y_{\mathrm{J}}\bar{R}_q(T',z,\cdot)\\
&+\nabla\cdot \sum_{\mathrm{J}\subseteq\mathrm{I}} \mathrm{D}^Y_{\mathrm{J}^c}\big(\delta \bar{a}_q(T,T')-\delta \bar{a}_q(T,T',z)\big)\nabla \mathrm{D}^Y_{\mathrm{J}}\bar{K}_{q}(T',z,\cdot)\\
:=&\nabla\cdot\Lambda_{q,y_1,\cdots,y_n}\quad\text{in $(0,t)\times \mathbb{R}^2$}.
\end{aligned}
\end{equation}
with $\mathrm{D}^Y_{\mathrm{I}}\delta\bar{R}_q(T,T',z,(t,\cdot))\equiv 0$. 

\medskip

To obtain \eqref{EqMalliavinDiffNoFinite}, we appeal to the definition \eqref{RegPartDef} in the form of $R_q=\bar{R}_q(\xi,\cdot,\cdot):=\bar{\Gamma}_q(\xi,\cdot,\cdot)-\bar{K}_{q}(\xi,\cdot,\cdot)$ where for any $z\in \mathbb{R}^{1+2}$ and $T\in X$
$$\partial_t \bar{\Gamma}_q(T,z,\cdot)+\nabla\cdot \bar{a}_q(T)\nabla \bar{\Gamma}_q(T,z,\cdot)=0\quad\text{with $\bar{\Gamma}_q(T,z,(t,\cdot))=\delta(x-x')$},$$
and $\bar{a}_q(T):=\zeta_q\star A(m_q\star T)$; and $\bar{K}_q$ is defined in \eqref{KernelXDep} together with
$$\partial_t \bar{K}_q(T,z,\cdot)+\nabla\cdot \bar{a}_q(T,z)\nabla \bar{K}_q(T,z,\cdot)=0\quad\text{with $\bar{K}_q(T,z,(t,\cdot))=\delta(x-x')$}.$$
Taking the difference between the two equations yields
\begin{equation}\label{EquationResteR}
\partial_\tau \bar{R}_q(T,z,\cdot)+\nabla\cdot \bar{a}_q(T)\nabla \bar{R}_q(T,z,\cdot)=\nabla\cdot \big(\bar{a}_q(T)-\bar{a}_q(T,z)\big)\nabla \bar{K}_{q}(T,z,\cdot)\quad\text{in $(0,t)\times \mathbb{R}^2$}\quad\text{with }\bar{R}_q(T,z,(t,\cdot))\equiv 0.
\end{equation}
We then take Malliavin derivatives to get \eqref{EqMalliavinDiffNoFinite}. To conclude, we obtain \eqref{EqMalliavinDiff} by taking finite differences in \eqref{EqMalliavinDiffNoFinite}.

\medskip

{\sc Step 2. Proof of $(a)$. }We only give the argument for the second item of \eqref{AssumptionInteOnDNF} (the first item is proven the same way), where it amounts to prove that there exist $\sigma_n,C_n>0$ such that for any $\varphi_1,\cdots,\varphi_n\in X$ with $\|\varphi_i\|_X=1$ it holds
\begin{equation}\label{WhatToProveLemma1}
\Big\vert \dd^n \delta\bar{R}_q(T,T',z,z').\bigotimes_{i=1}^n \varphi_i\Big\vert \lesssim_{n,q}  \mathds{1}_{t\geq t'}(1+\vert z\vert^2+\vert z'\vert^2)^{\sigma _n}(t-t')G_{\sqrt{t-t'}}\Big(\tfrac{1}{\sqrt{C_n}}(x-x')\Big)\|T-T'\|_{X}.
\end{equation}
The estimate \eqref{WhatToProveLemma1} is based on three ingredients: an estimate on Malliavin derivatives of $\nabla \bar{R}_q(T',z,\cdot)$ given by for any $\mathrm{I}=\{i_1,\cdots,i_m\}\subset\{1,\cdots,n\}$ and some $C_m,\sigma_m>0$ we have
\begin{equation}\label{BoundRBarNoDif}
\bigg\vert\int \dd y_{i_1}\cdots\dd y_{i_m}\, \varphi(y_1)\cdots \varphi(y_{i_m})\, \nabla \mathrm{D}^Y_{i_1,\cdots,y_n} \bar{R}_{q}(T,z,z')\bigg\vert \lesssim \mathds{1}_{t\geq t'}(1+\vert z\vert^2+\vert z'\vert^2)^{\sigma_m} (t-t')^{\frac{1}{2}}G_{\sqrt{t-t'}}\Big(\tfrac{1}{\sqrt{C_m}}(x-x')\Big),
\end{equation}
and on estimates from parabolic regularity theory as well as an induction argument. The proof of \eqref{BoundRBarNoDif} follows the lines of the argument we present for \eqref{WhatToProveLemma1}, therefore we leave the details to the reader. We control the left-hand side of \eqref{WhatToProveLemma1} by an explicit Green's function representation via the kernel $\bar{\Gamma}_q(T)$ of $\partial_\tau+\nabla\cdot \bar{a}_q(T)\nabla$:
\begin{equation}\label{GreenRepreDivN}
 \dd^n \delta\bar{R}_q(T,T',z,z').\bigotimes_{i=1}^n \varphi_i=-\mathds{1}_{t\geq t'}\int \dd z''\, \mathds{1}_{[t',t]}(t'')\bigg(\int \dd y_1\cdots \dd y_n\, \varphi(y_1)\cdots \varphi(y_n)\, \Lambda_{q,y_1,\cdots,y_n}(z'')\bigg)\cdot \nabla_1 \bar{\Gamma}_q(T,z'',z').
\end{equation}
Furthermore, note that $\|\bar{a}_q(T)\|_{\cc^{0,1}(\mathbb{R}^{1+2})}\lesssim_q 1$ which implies that $\bar{\Gamma}_q$ satisfies the following standard heat-kernel bound 
\begin{equation}\label{EstimateGreenDepq}
\big\vert \nabla_1\bar{\Gamma}_q(T,z'',z')\big\vert\lesssim_q \mathds{1}_{t''\geq t'}(t''-t')^{-\frac{1}{2}} G_{\sqrt{t''-t'}}\Big(\tfrac{1}{\sqrt{C}}(x''-x')\Big).
\end{equation}
To estimate \eqref{GreenRepreDivN}, we control for any $z''\in \mathbb{R}^{1+2}$, each contributions of $\int \dd y_1\cdots \dd y_n\, \varphi(y_1)\cdots \varphi(y_n)\, \Lambda_{q,y_1,\cdots,y_n}(z'')$ from the right-hand sides of \eqref{EqMalliavinDiff} separately. We fix $\mathrm{J}\subset \{1,\cdots,n\}$. First, from the definition of $\bar{a}_q(T)$ we have
\begin{equation}\label{EstiRHS0Induc}
\Big\vert\int \prod_{i\in \mathrm{J}^c}\dd y_i\, \prod_{i\in \mathrm{J}^c}\varphi_i(y_i)\, \mathrm{D}^{Y}_{\mathrm{J}^c}\bar{a}_{q}(T,z'')\Big\vert\lesssim\zeta_q\star \prod_{i\in \mathrm{J}^c}\vert m_q\star \varphi_i\vert(z'')\stackrel{\eqref{EstimateConvoMT}}{\lesssim} (1+\vert z''\vert^2)^{\vert\mathrm{J}^c\vert\sigma};
\end{equation}
Second, note that similarly to \eqref{MalliavinFrozenN}
$$\nabla \mathrm{D}^Y_{\mathrm{J}}\bar{K}_q(T,z,z'')=\bar{F}^{(\vert \mathrm{J}\vert)}_{z,z''}(m_q\star T(z))\prod_{i\in \mathrm{J}}m_q(z-y_i),$$
where $\bar{F}_{z,z''}$ is smooth and for any $n\geq 1$ and some $C_n>0$, $\sup_{\lambda\in \mathbb{R}}\vert \bar{F}^{(n)}_{z,z''}(\lambda)\vert\lesssim \mathds{1}_{t\geq t'}(t-t'')^{-\frac{1}{2}}G_{\sqrt{t-t''}}(\frac{1}{\sqrt{C_n}}(x-x''))$. We deduce
\begin{equation}\label{EstiRHS1Induc}
\begin{aligned}
&\Big\vert\int \dd y_1\cdots\dd y_n\, \varphi_1(y_1)\cdots\varphi_n(y_n)  \mathrm{D}^Y_{\mathrm{J}^c}\big(\bar{a}_q(T,z'')-\bar{a}_q(T,z)\big)\nabla \mathrm{D}^Y_{\mathrm{J}}\delta \bar{K}_{q}(T,T',z,z'')\Big\vert\\
&\lesssim \Big\vert \zeta_q\star A^{(\vert \mathrm{J}^c\vert)}(m_q\star T)\prod_{i\in \mathrm{J}^c}m_q\star \varphi_i(z'')-\zeta_q\star A^{(\vert \mathrm{J}^c\vert)}(m_q\star T)\prod_{i\in \mathrm{J}^c}m_q\star \varphi_i(z)\Big\vert\\
&\quad\times\Big(\bar{F}^{(\vert \mathrm{J}\vert)}_{z,z''}(m_q\star T(z))-\bar{F}^{(\vert \mathrm{J}\vert)}_{z,z''}(m_q\star T'(z))\Big)\prod_{i\in \mathrm{J}}m_q\star\varphi_i(z)\\
&\stackrel{\eqref{EstimateConvoMT}}{\lesssim_q} \mathds{1}_{t\geq t''}(1+\vert z\vert^2+\vert z''\vert^2)^{\sigma+1} (t-t'')^{\frac{1}{2}}G_{\sqrt{t-t''}}\Big(\tfrac{1}{\sqrt{C_n}}(x-x'')\Big)\|T-T'\|_X;
\end{aligned}
\end{equation}
Thirdly, for some $\sigma_n>0$
\begin{equation}\label{EstiRHS2Induc}
\begin{aligned}
&\Big\vert\int \dd y_1\cdots\dd y_n\, \varphi_1(y_1)\cdots\varphi_n(y_n)  \mathrm{D}^Y_{\mathrm{J}^c}\delta \bar{a}_q(T,T',z'')\nabla \mathrm{D}^Y_{\mathrm{J}}\delta \bar{R}_{q}(T',z,z'')\Big\vert\\
&\stackrel{\eqref{BoundRBarNoDif}}{\lesssim} \mathds{1}_{t\geq t''}\zeta_q\star \vert A^{(\vert \mathrm{J}^c\vert)}(m_q\star T)-A^{(\vert \mathrm{J}^c\vert)}(m_q\star T')\vert \Big\vert \prod_{i\in \mathrm{J}^c}m_q\star \varphi_i\Big\vert(z'') (1+\vert z\vert^2+\vert z''\vert^2)^{\sigma_{\vert\mathrm{J}\vert}} (t-t'')^\frac{1}{2} G_{\sqrt{t-t''}}\Big(\tfrac{1}{\sqrt{C}}(x-x'')\Big)\\
&\lesssim_q \mathds{1}_{t\geq t''}(1+\vert z\vert^2+\vert z''\vert^2)^{\sigma_n}(t-t'')^\frac{1}{2}G_{\sqrt{t-t''}}\Big(\tfrac{1}{\sqrt{C}}(x-x'')\Big)\|T-T'\|_X,
\end{aligned}
\end{equation}
and finally we obtain similarly
\begin{equation}\label{EstiRHS3Induc}
\begin{aligned}
&\Big\vert\int \dd y_1\cdots\dd y_n\, \varphi_1(y_1)\cdots\varphi_n(y_n)  \mathrm{D}^Y_{\mathrm{J}^c}\big(\delta \bar{a}_q(T,T',z'')-\delta \bar{a}_q(T,T',z)\big)\nabla \mathrm{D}^Y_{\mathrm{J}}\delta \bar{K}_{q}(T',z,z'')\Big\vert\\
&\lesssim_q \mathds{1}_{t\geq t''}(1+\vert z\vert^2+\vert z''\vert^2)^{\sigma+1} (t-t'')^{\frac{1}{2}}G_{\sqrt{t-t''}}\Big(\tfrac{1}{\sqrt{C_n}}(x-x'')\Big)\|T-T'\|_X.
\end{aligned}
\end{equation}
Plugging \eqref{EstiRHS0Induc}, \eqref{EstiRHS1Induc}, \eqref{EstiRHS2Induc} and \eqref{EstiRHS3Induc}  in \eqref{GreenRepreDivN} and using \eqref{EstimateGreenDepq}, we obtain for some $\sigma_n>0$
\begin{align}
&\Big\vert  \dd^n \delta\bar{R}_q(T,T',z,z').\bigotimes_{i=1}^n \varphi_i\Big\vert\nonumber\\
&\lesssim_q \mathds{1}_{t\geq t'}\sum_{\mathrm{J}\subset \llbracket 1,n\rrbracket}(1+\vert z\vert^2+\vert z'\vert^2)^{\vert \mathrm{J}^c\vert \sigma}\int \dd z''\, \Big\vert \nabla_2 \dd^{\vert \mathrm{J}\vert}\delta \bar{R}_q(T,T',z,z'').\bigotimes_{i\in \mathrm{J}}\varphi_i\Big\vert(t''-t)^{-\frac{1}{2}}G_{\sqrt{t''-t'}}\Big(\tfrac{1}{\sqrt{C}}(x''-x')\Big)\nonumber\\
&\quad+\mathds{1}_{t\geq t'}(1+\vert z\vert^2+\vert z'\vert^2)^{\sigma_n}\int \dd z''\, (t-t'')^\frac{1}{2}G_{\sqrt{t-t''}}\Big(\tfrac{1}{\sqrt{C}}(x-x'')\Big)(t''-t')^{-\frac{1}{2}}G_{\sqrt{t''-t'}}\Big(\tfrac{1}{\sqrt{C}}(x''-x')\Big)\|T-T'\|_X\nonumber\\
&\lesssim_q \mathds{1}_{t\geq t'}\sum_{\mathrm{J}\subset \llbracket 1,n\rrbracket}(1+\vert z\vert^2+\vert z'\vert^2)^{\vert \mathrm{J}^c\vert \sigma}\int \dd z''\, \Big\vert \nabla_2\dd^{\vert \mathrm{J}\vert}\delta \bar{R}_q(T,T',z,z'').\bigotimes_{i\in \mathrm{J}}\varphi_i\Big\vert (t''-t)^{-\frac{1}{2}}G_{\sqrt{t''-t'}}\Big(\tfrac{1}{\sqrt{C}}(x''-x')\Big)\nonumber\\
&\quad+ \mathds{1}_{t\geq t'}(1+\vert z\vert^2+\vert z'\vert^2)^{\sigma_n}(t-t')G_{\sqrt{t-t'}}\Big(\tfrac{1}{\sqrt{C}}(x-x')\Big)\|T-T'\|_X.\label{FirstEstiForIndcution}
\end{align}
We further estimate $ \nabla_2\dd^{\vert \mathrm{J}\vert}\delta \bar{R}_q(T,T',z,z'').\bigotimes_{i\in \mathrm{J}}\varphi_i$ using Schauder' estimates of Proposition \ref{Schauder} (recalling $\|\bar{a}_q(T)\|_{\cc^{0,1}(\mathbb{R}^{1+2})}\lesssim_q 1$) applied to \eqref{EqMalliavinDiff}: for any $t\geq t'$ 
\begin{equation}\label{SchauderDNDerivative}
\begin{aligned}
\Big\vert \nabla_2 \dd^{\vert \mathrm{J}\vert} \delta \bar{R}_q(T,T',z,z'').\bigotimes_{i\in\mathrm{J}}\varphi_i\Big\vert \lesssim& (t-t'')^{-\frac{1}{2}}\bigg(\fint_{\cc_{\sqrt{t-t'}/2}(z')}\dd y\,\Big\vert \dd^{\vert \mathrm{J}\vert} \delta \bar{R}_q(T,T',z,z'').\bigotimes_{i\in\mathrm{J}}\varphi_i\Big\vert^2\bigg)^\frac{1}{2}\\
&+\sup_{\cc_{\sqrt{t-t'}/2}(z''))}\bigg\vert \int \prod_{i\in\mathrm{J}}\dd y_i\, \prod_{i\in\mathrm{J}} \varphi_i(y_i)\,\Lambda_{q,\{y_i\}_{i\in \mathrm{J}}}\bigg\vert\\
&+(t-t'')^{\frac{\alpha}{2}}\bigg\vert \int \prod_{i\in\mathrm{J}}\dd y_i\,\prod_{i\in\mathrm{J}}\varphi_i(y_i)\,\,\Lambda_{q,\{y_i\}_{i\in \mathrm{J}}}\bigg\vert_{\cc^\alpha(\cc_{\sqrt{t-t'}/2}(z''))},
\end{aligned}
\end{equation}
where the third right-hand side term is bounded similarly as in \eqref{EstiRHS1Induc}, \eqref{EstiRHS2Induc} and \eqref{EstiRHS3Induc}. To conclude, the combination of \eqref{FirstEstiForIndcution} and \eqref{SchauderDNDerivative} implies \eqref{WhatToProveLemma1} by induction.

\medskip

{\sc Step 3.  Proof of $(b)$. }The proof of $(b)$ follows the lines of Step $2$ with \eqref{EqMalliavinDiffNoFinite}, where the only difference is that we need to establish estimates uniformly in 
$q$ and $\delta$. This can be achieved with a few modifications. First, instead of using \eqref{EstimateGreenDepq}, we appeal to the regularity \eqref{RegCoefAHolder} of $a$ which has the effect of replacing \eqref{EstimateGreenDepq} with an annealed counterpart\footnote{where we use that the dependence of heat-kernel estimates depends polynomially on local H\"older norms of $a$} that holds uniformly in $q$: for any $p<\infty$ and some $\gamma>0$ depending on $\alpha$
$$\mathbb{E}^\frac{1}{p}\big[\vert \nabla_1 \Gamma_q(z'',z')\vert^p\big]\lesssim \mathds{1}_{t''\geq t'}\max\big\{1,(t''-t')^\gamma\big\}(t''-t')^{-\frac{1}{2}}G_{\sqrt{t''-t'}}\Big(\tfrac{1}{\sqrt{C}}(x''-x')\Big).$$
Second, we now bound uniformly $\sup_{y\in \mathbb{R}^{1+2}}\vert m_q\star \rho_\delta(y)\vert \lesssim \|m\|_{\LL^\infty(\mathbb{R}^{1+2})}$ which yields uniform bounds in $\delta$ and removes the factors like $(1+\vert z\vert^2+\vert z'\vert^2)^{\sigma_n}$ for $\sigma_n>0$. Thirdly, similar estimates to that of \eqref{EstiRHS1Induc}, \eqref{EstiRHS2Induc} and \eqref{EstiRHS3Induc} hold in expectation uniformly in $q$ using once more the regularity \eqref{RegCoefAHolder} of $a$, where $(t-t')^\frac{1}{2}$ is weakened to $(t-t')^{\frac{\alpha}{2}}$. Finally, we conclude by induction establishing similar inequalities to that of \eqref{FirstEstiForIndcution} and \eqref{SchauderDNDerivative}, where the later have an extra multiplicative polynomial factor depending on $1+(t-t'')^{\frac{\alpha}{2}}\vert a_q\vert_{\cc_{\sqrt{t-t''}/2}(z'')}$ (see Proposition \ref{Schauder}). Appealing to \eqref{RegCoefAHolder} and after taking expectation, it yields an extra fractor $(1+(t-t'')^{\frac{\alpha}{2}+\varepsilon})^\gamma$ for any $\varepsilon>0$ but provides uniform estimates in $q$. 

\medskip

{\sc Pass to the limit as $q\uparrow\infty$. }We now show that for any $p<\infty$ and $m\in \{1,2,3\}$
\begin{equation}\label{PassToLimitRest}
\mathbb{E}^{\frac{1}{2p}}\bigg[\Big(\int \dd z\, \psi^\lambda(z)\Big(\int \dd z'\, R_q(z,z')\xi_\delta(z')\Big)^m\Big)^{2p}\bigg]\underset{q\uparrow\infty}{\longrightarrow}\mathbb{E}^{\frac{1}{2p}}\bigg[\Big(\int \dd z\, \psi^\lambda(z)\Big(\int \dd z'\, R(z,z')\xi_\delta(z')\Big)^m\Big)^{2p}\bigg].
\end{equation}
We give the argument for $m=2$, the cases $m=1,3$ are obtained the same way. Noticing that 
\begin{align*}
&\mathbb{E}^{\frac{1}{2p}}\bigg[\bigg(\int \dd z\, \psi^\lambda(z)\Big(\Big(\int \dd z'\, R_q(z,z')\xi_\delta(z')\Big)^2-\Big(\int \dd z'\, R(z,z')\xi_\delta(z')\Big)^2\Big)\bigg)^{2p}\bigg]\\
&=\mathbb{E}^{\frac{1}{2p}}\bigg[\Big(\int \dd z\, \psi^\lambda(z)\Big(\int \dd z'\, (R_q(z,z')-R(z,z'))\xi_\delta(z')\Big)\Big(\int \dd z'\, (R_q(z,z')+R(z,z'))\xi_\delta(z')\Big)\Big)^{2p}\bigg],
\end{align*}
using for any $p\geq 1$, $\mathbb{E}^\frac{1}{p}[(\xi_\delta(z'))^p]\lesssim_p \delta^{-2^{-}}$ and from the first item of \eqref{BoundMalliavinDerivativeRest} (which also holds for $R$ the same way), by Lebesgue convergence theorem it suffices to show that for any $p\geq 1$ and $z,z'\in \mathbb{R}^{1+2}$
\begin{equation}\label{EstiDiffPassLimit}
\mathbb{E}^\frac{1}{p}\big[\vert R_q(z,z')-R(z,z')\vert^p\big]\underset{q\uparrow\infty}{\rightarrow}0.
\end{equation}
We first write from \eqref{EquationResteR} with $T=\xi$ the PDE solved by $R_q(z,\cdot)-R(z,\cdot)$:
\begin{equation}\label{EqDiffRqPassLimit}
\begin{aligned}
\partial_\tau (R_q-R)(z,\cdot)-\nabla\cdot a\nabla (R_q-R)(z,\cdot)=&\nabla\cdot (a-a_q)\nabla R_q(z,\cdot)+\nabla\cdot (a_q-a_q(z))\nabla (K_q(z,\cdot)-K(z,\cdot))\\
&+\nabla\cdot \big((a-a_q)-(a-a_q)(z)\big)\nabla K(z,\cdot)\\
:=&\nabla\cdot \Lambda_q(z,\cdot).
\end{aligned}
\end{equation}
This yields the Green's function representation:
\begin{equation}\label{GreenPassLimit}
(R_q-R)(z,z')=-\mathds{1}_{t\geq t'}\int \dd z'' \mathds{1}_{[t',t]}(t'')\Lambda_q(z,z'')\cdot \nabla_1 \Gamma_q(z'',z').
\end{equation}
We then estimate each right-hand side term of \eqref{EqDiffRqPassLimit}: first using \eqref{BoundMalliavinDerivativeRest} and $m\in \mathrm{B}^\alpha_{2,\infty}(\mathbb{R}^{1+2})$
\begin{equation}\label{Esti1PassLimitAMAq}
\begin{aligned}
&\mathbb{E}^\frac{1}{p}\big[\vert(a-a_q)(z'')\nabla_2 R_q(z,z'')\vert^p\big]\\
&\lesssim \mathds{1}_{t\geq t''}(t-t'')^{\frac{\alpha'}{2}}G_{\sqrt{t-t''}}\Big(\tfrac{1}{\sqrt{C}}(x-x'')\Big)\int \dd y\,\zeta_q(z''-y)\mathbb{E}^\frac{1}{2p}\big[\vert m\star \xi(z'')-m_q\star \xi(y)\vert^{2p}\big]\\
&\lesssim  \mathds{1}_{t\geq t''}(t-t'')^{\frac{\alpha'}{2}}G_{\sqrt{t-t''}}\Big(\tfrac{1}{\sqrt{C}}(x-x'')\Big)\int \dd y\,\zeta_q(z''-y)\|m(z''-\cdot)-m_q(y-\cdot)\|_{\LL^2(\mathbb{R}^{1+2})}\\
&\lesssim \mathds{1}_{t\geq t''}(t-t'')^{\frac{\alpha'}{2}}G_{\sqrt{t-t''}}\Big(\tfrac{1}{\sqrt{C}}(x-x'')\Big)\Big(\int \dd y\, \zeta_q(y)\vert y\vert^\alpha+\|m(z''-\cdot)-m_q(z''-\cdot)\|_{\LL^2(\mathbb{R}^{1+2})}\Big).
\end{aligned}
\end{equation}
Second, using \eqref{KernelXDep} with $T=\xi$ and $m\in \mathrm{B}^\alpha_{2,\infty}(\mathbb{R}^{1+2})$
\begin{align*}
&\mathbb{E}^\frac{1}{p}\big[\vert(a(z'')-a_q(z))(\nabla_2 K_q(z,z'')-\nabla_2 K(z,z''))\vert^p\big]\\
&\lesssim \mathds{1}_{t\geq t'}\vert z''-z\vert^\alpha (t-t'')^{-\frac{1}{2}}G_{\sqrt{t-t''}}\Big(\tfrac{1}{\sqrt{C}}(x-x'')\Big)\|m_q(z-\cdot)-m(z-\cdot)\|_{\LL^2(\mathbb{R}^{1+2})}.
\end{align*}
Thirdly similarly to \eqref{Esti1PassLimitAMAq}
\begin{align*}
&\mathbb{E}^\frac{1}{p}\big[\vert\big((a-a_q)(z'')-(a-a_q)(z)\big)\nabla_2 K(z,z'')\vert^p\big]\\
&\lesssim (t-t'')^{-\frac{1}{2}}G_{\sqrt{t-t''}}\Big(\tfrac{1}{\sqrt{C}}(x-x'')\Big)\Big(\int \dd y\, \zeta_q(y)\vert y\vert^\alpha+\|m(z''-\cdot)-m_q(z''-\cdot)\|_{\LL^2(\mathbb{R}^{1+2})}\Big)\\
&\quad\times\Big(\int \dd y\, \zeta_q(y)\vert y\vert^\alpha+\|m(z-\cdot)-m_q(z-\cdot)\|_{\LL^2(\mathbb{R}^{1+2})}\Big).
\end{align*}
We conclude by injecting the three last estimates into \eqref{GreenPassLimit} and by letting $q\uparrow\infty$ using Lebesgue's convergence theorem.
\end{proof}
\section{Malliavin calculus and Wick products}\label{AppendixGIBP}
In this section, we consider the space dimension to be $d\geq 2$. We interpret the space-time white noise $\xi$ as a random distribution, where its law $\mathbb{P}$ is a Gaussian measure with covariance given by the $\LL^2(\mathbb{R}^{1+d})$-inner product on a Banach space $X$ of Schwartz distributions such that $\LL^2(\mathbb{R}^{1+d})\subset X$ is compact in a Hilbert-Schmidt way (see \cite[Example 3.9.7.]{bogachev1998gaussian}). A choice is to consider the Hilbert space 
$$X:=\Big\{T\in \mathcal{S}'(\mathbb{R}^{1+d}): \int_{\mathbb{R}^{1+d}} \big\vert S\,T\big\vert^2<\infty\Big\},$$
where $S:=(\partial_t+1-\Delta)^{-\sigma}(1+\vert\cdot\vert^2)^{-\sigma}$ for a fixed exponent $\sigma>\frac{d}{4}$ and $\vert \cdot\vert$ denotes the parabolic norm, i. e. $\vert (t,x)\vert^2=\vert t\vert+\vert x\vert^2$ for any $(t,x)\in \mathbb{R}^{1+d}$. The operator $S$ is indeed Hilbert-Schmidt from $\LL^2(\mathbb{R}^{1+d})$ to $\LL^2(\mathbb{R}^{1+d})$: given an orthonormal basis $\{e_n\}_n$ of $\LL^2(\mathbb{R}^{1+d})$, we have by applying the space-time Fourier transform $\mathcal{F}$,
\begin{align*}
    \sum_{n\geq 0}\big(S e_n,S e_n\big)&=\sum_{n\geq 0}\int_{\mathbb{R}^{1+d}}\vert (\partial_t+1-\Delta)^{-\sigma}(1+\vert\cdot\vert^2)^{-\sigma}e_n\vert^2\\
    &=\sum_{n\geq 0}\int_{\mathbb{R}}\dd \eta \int_{\mathbb{R}^d} \dd\zeta\,\big(\vert \eta\vert^2+(1+\vert\zeta\vert^2)^2\big)^{-\sigma}\vert\mathcal{F}\big((1+\vert\cdot\vert^2)^{-\sigma}e_n\big)\vert^2\\
    &=\int_{\mathbb{R}}\dd \eta \int_{\mathbb{R}^d} \dd\zeta\,\big(\vert \eta\vert^2+(1+\vert\zeta\vert^2)^2\big)^{-\sigma}\int \dd t\,\int\dd x\,(1+\vert t\vert+\vert x\vert^2)^{-2\sigma}<\infty.
\end{align*}
%
%
%
In the following, we write $(\cdot,\cdot)_X$ for the inner product in $X$ and we recall that for any centred random variables $\{X_i\}_{i\leq n}$ with $n\geq 1$, the Wick product $\displaystyle\mathop{\diamond}_{i=1}^n X_i$ is defined recursively by 
\begin{equation}\label{WickProductDef}
\mathop{\diamond}_{i=1}^n X_i=X_n\mathop{\diamond}_{i=1}^{n-1} X_i-\sum_{j=1}^{n-1}\mathbb{E}\big[X_n X_j\big]\mathop{\diamond}_{\underset{\ell\neq j}{\ell=1}}^{n-1} X_\ell,
\end{equation}
where we set by convention $\diamond_{i=n}^m X_i=1$ for any $m<n$. We recall some basic properties of \eqref{WickProductDef}:
\begin{itemize}
\item[(i)]The Wick product is invariant by permutation, i.e. for any permutation $\sigma$ of $\{1,\cdots,n\}$ it holds
\begin{equation}\label{PermutationWick}
\mathop{\diamond}_{i=1}^n X_i=\mathop{\diamond}_{i=1}^n X_{\sigma(i)}.
\end{equation}
\item[(ii)]The Wick product satisfies the following differential property for any $\ell\in \{1,\cdots,n\}$
\begin{equation}\label{DiffEqWick}
\frac{\partial }{\partial X_\ell}\mathop{\diamond}_{i=1}^n X_i=\mathop{\diamond}_{\underset{i\neq \ell}{i=1}}^n X_i.
\end{equation}
\end{itemize}

We now prove Lemma \ref{IntWRTWick}. 
\begin{proof}[Proof of Lemma \ref{IntWRTWick}]
The proof follows from a finite-rank approximation on $X$ and we first define the projectors. We use that $\sqrt{S^\star S}$ is self-adjoint and Hilbert-Schmidt on $\LL^2(\mathbb{R}^{1+d})$, so there exists an orthonormal basis $\{h_j\}_{j\geq 1}$ of $\LL^2(\mathbb{R}^{1+d})$ and eigenvalues $\{\lambda_j\}_{j\geq 1}\subset \ell^2(\mathbb{N})$ such that  
\begin{equation}\label{DefEigenvalues}
\sqrt{S^\star S}\,h_j=\lambda_j\, h_j.
\end{equation}
We use the family $\{h_j,\lambda_j\}_{j\geq 1}$ to define the projectors on $X$: for any $N\geq 1$
$$P_N\xi:=\sum_{j=1}^N \big(\xi,\lambda^{-2}_j h_j\big)_{X}\,h_j\quad \text{for any $\xi\in X$}.$$
The family $\{P_N\}_{N\geq 1}$ are indeed finite-rank projectors as from \eqref{DefEigenvalues} and the definition of $(\cdot,\cdot)_X$,  $\{\lambda^{-1}_j h_j\}_{j\geq 1}$ is an orthonormal basis in $X$. Furthermore, note that from \eqref{DefEigenvalues}, it also holds
\begin{equation}\label{RestrictL2}
\big(\xi,\lambda^{-2}_j h_j\big)_{X}=\big(\xi,h_j\big)_{\LL^{2}(\mathbb{R}^{1+d})}\quad\text{for any $\xi\in \LL^2(\mathbb{R}^{1+d})$.}
\end{equation}
Finally we recall that in $\{h_j\}_{j\geq 1}$, the kernel of $\dd^n F$ is given by \eqref{FormulaKernelL2}.

\medskip
Now, we define the finite-rank approximation of the functional $F$: for any $G\in \mathbb{R}^N$ and $N\geq 1$,
$$F_N(G):=F\bigg(\sum_{j=1}^N G_j h_j\bigg).$$
We now turn to the proof of \eqref{GaussianInteNnoises}. We first note that in the formalism of Gaussian measures, it holds
$$\mathbb{E}\bigg[F(\xi)\, \mathop{\diamond}_{i=1}^n (\xi,\varphi_i)\bigg]=\mathbb{E}\bigg[F\, \mathop{\diamond}_{i=1}^n \ell_i\bigg],$$
where $\ell_i:=(\cdot,\varphi_i)_{\LL^2(\mathbb{R}^{1+d})}$. We then use \eqref{AssumptionInteOnF} and $P_N\underset{N\uparrow\infty}{\rightarrow}\mathrm{Id}$ in $X$, to obtain
$$\mathbb{E}\bigg[F\, \mathop{\diamond}_{i=1}^n \ell_i\bigg]=\lim_{N\uparrow\infty}\mathbb{E}\bigg[F_N\big(f_1,\cdots,f_N\big)\, \mathop{\diamond}_{i=1}^n \ell_i\bigg],$$
where $f_i:=\big(\cdot,\lambda^{-2}_i h_i\big)_{X}$. We then appeal to finite dimensional Gaussian calculus: setting $G_N:=\{g_i\}_{i\leq N}$ a Gaussian vector with law $(f_1,\cdots,f_N)\#\mathbb{P}$ and $\{X_i\}_{i\leq n}$ with law $(\ell_1,\cdots,\ell_n)\#\mathbb{P}$ (where for simplicity we still write $\mathbb{E}$ for the underlying ensemble average of $\{G_N,\{X_i\}_i\}$),  it holds by Stein's lemma \cite{Stein1981} for any $N\geq 1$
\begin{align*}
\mathbb{E}\bigg[F_N\big(f_1,\cdots,f_N\big)\, \mathop{\diamond}_{i=1}^n \ell_i\bigg]&=\mathbb{E}\bigg[F_N(G_N)\mathop{\diamond}_{i=1}^n X_i\bigg]\\
&\stackrel{\eqref{WickProductDef}}{=}\mathbb{E}\bigg[F_N(G_N)X_n\mathop{\diamond}_{i=1}^{n-1} X_i\bigg]-\sum_{j=1}^{n-1}\mathbb{E}\big[X_n X_j\big]\mathbb{E}\bigg[F_N(G_N)\mathop{\diamond}_{\underset{\ell\neq j}{\ell=1}}^{n-1} X_\ell\bigg]\\
&=\sum_{j=1}^N\mathbb{E}\bigg[\partial_j F_N(G_N)\mathop{\diamond}_{i=1}^{n-1} X_i\bigg]\mathbb{E}\big[g_j \, X_n\big]\\
&=\sum_{j=1}^N\mathbb{E}\bigg[\dd F\bigg(\sum_{i=1}^N g_i\, h_i\bigg).h_j\,\mathop{\diamond}_{i=1}^{n-1} X_i\bigg]\mathbb{E}\big[g_j \, X_n\big]\\
&=\sum_{j=1}^N\mathbb{E}\bigg[\dd F\bigg(\sum_{i=1}^N g_i\, h_i\bigg).h_j\,\mathop{\diamond}_{i=1}^{n-1} X_i\bigg]\big(h_j,\varphi_n\big)_{\LL^2(\mathbb{R}^{1+d})}\\
&=\int_{\mathbb{R}^{1+d}} \dd y\, \mathbb{E}\bigg[\sum_{j=1}^N\dd F\bigg(\sum_{i=1}^N g_i\, h_i\bigg).h_j\, h_j(y)\varphi_n(y)\,\mathop{\diamond}_{i=1}^{n-1} X_i\bigg].
\end{align*}
We then iterate and use \eqref{AssumptionInteOnDNF} to pass to the limit as $N\uparrow \infty$. 
\end{proof}
\section{Parabolic regularity theory}
We recall in this section standard parabolic theory. We first state standard Schauder's estimates (which can be found in \cite{krylov1996lectures}). 
\begin{proposition}[Schauder estimates]\label{Schauder}
Let $r>0$, $z\in \mathbb{R}^{1+d}$, $f\in \cc^{\alpha}(\mathbb{R}^{1+d})$ (for some $\alpha>0$) and $u$ be a weak solution of
$$\partial_\tau u+\nabla\cdot a\nabla u=\nabla\cdot f\quad\text{in $\cc_r(z)$}.$$
Then, there exists $\gamma>0$ depending on $d$ and $\alpha$ such that 
$$\sup_{\cc_{\frac{r}{2}(z)}}\vert \nabla u\vert\lesssim \Big(1+r^\alpha \vert a\vert_{\cc^\alpha(\cc_r(z))}\Big)^\gamma\bigg(r^{-1}\Big(\fint_{\cc_r(z)}\vert u\vert^2\Big)^\frac{1}{2}+\sup_{\cc_r(z)}\vert f\vert+r^\alpha\vert f\vert_{\cc^{\alpha}(\cc_r(z))}\bigg),$$
where $(1+r^\alpha \vert a\vert_{\cc^\alpha(\cc_r(z))})^\gamma$ can be removed if $ \vert a\vert_{\cc^\alpha(\mathbb{R}^{1+2})}<\infty$.
\end{proposition}
We then state and prove some localized energy estimates for parabolic equations that is used in Proposition \ref{PropRegPart}.
\section*{Acknowledgments}
The author warmly thanks Harprit Singh for introducing him to the field of singular SPDEs and for useful discussions and fruitful comments. The author would like to thank the Oberwolfach Research Fellows (OWRF) program for supporting a research
stay during which a part of this work was carried out, as well as Rhys Steele and Lucas Broux for valuable
discussions during that stay.
\bibliographystyle{plain}
\bibliography{references}
\end{document}